# INTERMITTENCY IN A CATALYTIC RANDOM MEDIUM[1]

By J. Gärtner and F. den Hollander

*Technische Universität Berlin and Leiden University*

In this paper, we study intermittency for the parabolic Anderson equation $\partial u/\partial t = \kappa\Delta u + \xi u$, where $u:\mathbb{Z}^d \times [0,\infty) \to \mathbb{R}$, $\kappa$ is the diffusion constant, $\Delta$ is the discrete Laplacian and $\xi:\mathbb{Z}^d \times [0,\infty) \to \mathbb{R}$ is a space-time random medium. We focus on the case where $\xi$ is $\gamma$ times the random medium that is obtained by running independent simple random walks with diffusion constant $\rho$ starting from a Poisson random field with intensity $\nu$. Throughout the paper, we assume that $\kappa, \gamma, \rho, \nu \in (0,\infty)$. The solution of the equation describes the evolution of a "reactant" $u$ under the influence of a "catalyst" $\xi$.

We consider the annealed Lyapunov exponents, that is, the exponential growth rates of the successive moments of $u$, and show that they display an interesting dependence on the dimension $d$ and on the parameters $\kappa, \gamma, \rho, \nu$, with qualitatively different intermittency behavior in $d = 1, 2$, in $d = 3$ and in $d \geq 4$. Special attention is given to the asymptotics of these Lyapunov exponents for $\kappa \downarrow 0$ and $\kappa \to \infty$.

## 1. Introduction and main results.

1.1. *Motivation.* The parabolic Anderson equation is the partial differential equation

$$\frac{\partial}{\partial t} u(x,t) = \kappa\Delta u(x,t) + \xi(x,t)u(x,t), \qquad x \in \mathbb{Z}^d, t \geq 0. \tag{1.1}$$

Here, the $u$-field is $\mathbb{R}$-valued, $\kappa \in (0,\infty)$ is the diffusion constant and $\Delta$ is the discrete Laplacian, acting on $u$ as

$$\Delta u(x,t) = \sum_{\substack{y\in\mathbb{Z}^d \\ \|y-x\|=1}} [u(y,t) - u(x,t)] \tag{1.2}$$

Received April 2004; revised July 2005.
[1]Supported in part by DFG-Schwerpunkt 1033.

(where $\|\cdot\|$ is the Euclidean norm), while

$$\xi = \{\xi(x,\cdot) : x \in \mathbb{Z}^d\} \tag{1.3}$$

is an $\mathbb{R}$-valued random field that evolves with time and drives the equation.

Equation (1.1) is the parabolic analogue of the Schrödinger equation in a random potential. It is a discrete heat equation with the $\xi$-field playing the role of a source or sink. One interpretation, coming from the study of population dynamics, is that $u(x,t)$ is the average number of particles at site $x$ at time $t$ when particles perform independent simple random walks at rate $\kappa$, split into two at rate $\xi(x,t)$ when $\xi(x,t) > 0$ (source term) and die at rate $-\xi(x,t)$ when $\xi(x,t) < 0$ (sink term). For more background on applications, the reader is referred to the monograph by Carmona and Molchanov ([4], Chapter I).

What makes (1.1) particularly interesting is that the two terms in the right-hand side *compete with each other*: the diffusion induced by $\Delta$ tends to make $u$ flat, while the branching induced by $\xi$ tends to make $u$ irregular. Consequently, in the population dynamics context, there is a competition between particles spreading out by diffusion and particles clumping around the areas where the sources are large.

A systematic study of the parabolic Anderson model for *time-independent* random fields $\xi$ has been carried out by Gärtner and Molchanov [18, 19, 20], Gärtner and den Hollander [12], Gärtner and König [14], Gärtner, König and Molchanov [16, 17] and Biskup and König [1, 2] (for a survey, see Gärtner and König [15]). The focus of these papers is on the study of the dominant spatial peaks in the $u$-field in the limit of large $t$, in particular, the height, the shape and the location of these peaks. Both the discrete model on $\mathbb{Z}^d$ (with i.i.d. $\xi$-fields) and the continuous model on $\mathbb{R}^d$ (with Gaussian and Poisson-like $\xi$-fields) have been investigated in the *quenched* setting (i.e., conditioned on $\xi$) as well as in the *annealed* setting (i.e., averaged over $\xi$).

Most of the theory currently available for *time-dependent* random fields $\xi$ is restricted to the situation where the components of the $\xi$-field are *uncorrelated in space and time*. Carmona and Molchanov ([4], Chapter III) have obtained an essentially complete qualitative description of the annealed Lyapunov exponents, that is, the exponential growth rates of the successive moments of $u(0,t)$ averaged w.r.t. $\xi$, for the case where the components of $\xi$ are independent Brownian noises. The quenched Lyapunov exponent, that is, the exponential growth rate of $u(0,t)$ conditioned on $\xi$, is harder to analyze. Carmona, Molchanov and Viens [5], Carmona, Koralov and Molchanov [3] and Cranston, Mountford and Shiga [6] have computed the asymptotics for $\kappa \downarrow 0$ of the quenched Lyapunov exponent for independent Brownian noises, which turns out to be singular. Cranston, Mountford and Shiga [7] have extended this result to independent Lévy noises. Further refinements for independent Brownian noises are obtained in Greven and den Hollander [21], including sharp bounds on the critical values of $\kappa$ where the annealed Lyapunov exponents change from positive to zero (resp. the quenched Lyapunov exponent changes from negative to zero), as well as a description of the equilibrium behavior when the quenched Lyapunov exponent is zero. These results are obtained from variational expressions for the Lyapunov exponents and are valid for general random walk transition kernels replacing $\Delta$.

In the present paper, we will be considering the situation where $\xi$ is given by

$$\xi(x,t)=\gamma\sum_k \delta_{Y_k(t)}(x) \tag{1.4}$$

with $\gamma \in (0,\infty)$ a coupling constant and

$$\{Y_k(\cdot): k \in \mathbb{N}\} \tag{1.5}$$

a collection of independent continuous-time simple random walks with diffusion constant $\rho \in (0,\infty)$ starting from a Poisson random field with intensity $\nu \in (0,\infty)$ (the index $k$ is an arbitrary numbering). As initial condition for (1.1), we take, for simplicity,

$$u(\cdot,0)\equiv 1. \tag{1.6}$$

We are interested in computing the annealed Lyapunov exponents of $u$ and studying their dependence on the parameters $\kappa$ and $\gamma,\rho,\nu$.

The population dynamics interpretation of (1.1) and (1.4)–(1.6) is as follows. Consider a spatially homogeneous system of two types of particles, $A$ (catalyst) and $B$ (reactant), performing independent continuous-time simple random walks such that:

(i) $B$-particles split into two at a rate that is $\gamma$ times the number of $A$-particles present at the same location;
(ii) $\rho$ and $\kappa$ are the diffusion constants of the $A$- and $B$-particles, respectively;
(iii) $\nu$ and 1 are the initial intensities of the $A$- and $B$-particles, respectively.

Then

$$u(x,t) = \text{the average number of } B\text{-particles at site } x \text{ at time } t \text{ conditioned on the evolution of the } A\text{-particles.} \tag{1.7}$$

Kesten and Sidoravicius [23] recently investigated this model with the addition of the following assumption:

(iv) $B$-particles die at rate $\delta \in (0,\infty)$.

The latter amounts to the transformation

$$u(x,t) \to u(x,t)e^{-\delta t}. \tag{1.8}$$

We describe their results in Section 1.4.

For a single moving catalyst, that is, $\xi(x,t) = \gamma\delta_{Y(t)}(x)$, the annealed Lyapunov exponents have recently been analyzed in Gärtner and Heydenreich [11]. The results are qualitatively different from ours and are, in fact, more closely related to those of Carmona and Molchanov [4] for white noise potentials.

1.2. *Catalytic and intermittent behavior.* Let $\langle\cdot\rangle$ denote expectation w.r.t. the $\xi$-field. For $p \in \mathbb{N}$ and $t > 0$, define

$$\Lambda_p(t) = \frac{1}{t}\log(e^{-\nu\gamma t}\langle u(0,t)^p\rangle^{1/p}). \tag{1.9}$$

This quantity monitors the effect of the randomness in the $\xi$-field on the growth of the $p$th moment. Indeed, if we would replace $\xi(x,t)$ in (1.1) by its average value $\langle\xi(x,t)\rangle = \nu\gamma$ [according to (1.4)], then the solution would be $u(\cdot,t) \equiv e^{\nu\gamma t}$, resulting in $\Lambda_p(\cdot) \equiv 0$.

The key quantities of interest in the present paper are the following *Lyapunov exponents*:

$$\begin{aligned} \widehat{\lambda}_p &= \lim_{t\to\infty}\frac{1}{t}\log\Lambda_p(t), \\ \lambda_p &= \lim_{t\to\infty}\Lambda_p(t). \end{aligned} \tag{1.10}$$

[Note that $\lambda_p$ is related to the moment Lyapunov exponent $\widetilde{\lambda}_p = \lim_{t\to\infty}\frac{1}{t}\times \log\langle u(0,t)^p\rangle$ via the relation $\lambda_p = \widetilde{\lambda}_p/p - \nu\gamma$.] The existence of the limits is not a priori evident and needs to be established. This will be done in Section 3 for $\widehat{\lambda}_p$ and in Section 4.1 for $\lambda_p$. From the Feynman–Kac representation for the moments of the solution of (1.1) and (1.4)–(1.6), given in Proposition 2.1 of Section 2.1, it will follow that $t \mapsto t\Lambda_p(t)$ is strictly positive and strictly increasing on $(0,\infty)$. Hence, $\widehat{\lambda}_p, \lambda_p \ge 0$. Further, we have $\Lambda_p(t) \ge \Lambda_{p-1}(t)$ by Hölder's inequality applied to the definition of $\Lambda_p(t)$. Hence, $\lambda_p \ge \lambda_{p-1}$. We will see in Section 4.3 that $\lambda_p > 0$.

Depending on the values of these Lyapunov exponents, we distinguish the following types of behavior.

DEFINITION 1.1. For $p \in \mathbb{N}$, we say that the solution is:

(a) strongly $p$-catalytic if $\widehat{\lambda}_p > 0$;
(b) weakly $p$-catalytic if $\widehat{\lambda}_p = 0$.

The solution being strongly catalytic means that the moments of the $u$-field grow much faster in the random medium $\xi$ than in the average medium $\langle\xi\rangle$, at a double-exponential rate. Weakly catalytic corresponds to a slower rate. Strongly catalytic behavior comes from an extreme form of clumping in the $\xi$-field.

DEFINITION 1.2. For $p \in \mathbb{N} \setminus \{1\}$, we say that the solution is:

(a) strongly $p$-intermittent if either $\lambda_p = \infty$ or $\lambda_p > \lambda_{p-1}$;
(b) weakly $p$-intermittent if $\lambda_p < \infty$ and $\lambda_p = \lambda_{p-1}$.

The solution being strongly $p$-intermittent means that the $1/p$th power of the $p$th moment of the $u$-field grows faster than the $1/(p-1)$th power of the $(p-1)$th moment, at an exponential rate. Weakly $p$-intermittent corresponds to a slower rate. Strongly intermittent behavior also comes from clumping in the $\xi$-field, but in a less extreme form than for strongly catalytic behavior. Note that strong $p$-intermittency implies strong $q$-intermittency for all $q > p$ (see Gärtner and Molchanov [18]). Also, note that our definition of weakly intermittent includes the possibility of no separation of the moments, usually called *nonintermittent*.

In the population dynamics context, both catalytic and intermittent behavior come from the $B$-particles clumping around the areas where the $A$-particles are clumping. It signals the appearance of rare high peaks in the $u$-field close to rare high peaks in the $\xi$-field. These peaks dominate the moments of the $u$-field (for more details, see [18], [26], Lecture 8, [22], Chapter 8, and [15]).

1.3. *Main theorems.* Let

$$\widehat{\varphi}(k) = \sum_{\substack{x\in\mathbb{Z}^d \\ \|x\|=1}} [1 - \cos(k \cdot x)], \qquad k \in [-\pi, \pi)^d. \tag{1.11}$$

For $\mu \geq 0$, define

$$R(\mu) = \frac{1}{(2\pi)^d} \int_{[-\pi,\pi)^d} \frac{dk}{\mu + \widehat{\varphi}(k)} \tag{1.12}$$

and put

$$r_d = \frac{1}{R(0)} \begin{cases} = 0, & \text{if } d = 1, 2, \\ > 0, & \text{if } d \geq 3. \end{cases} \tag{1.13}$$

Note that $R(\mu)$ is the Fourier representation of the kernel of the resolvent $(\mu - \Delta)^{-1}$ at 0; $R(0)$ equals the Green function at the origin of simple random walk on $\mathbb{Z}^d$ jumping at rate $2d$, that is, the Markov process generated by $\Delta$.

The following elementary and well-known fact is needed for Theorem 1.4(i) below (see Figure 1).

LEMMA 1.3. *For $r \in (0,\infty)$, let*

$$\mu(r) = \sup \mathrm{Sp}(\Delta + r\delta_0) \tag{1.14}$$

*denote the supremum of the spectrum of the operator $\Delta + r\delta_0$ in $\ell^2(\mathbb{Z}^d)$. Then*

(i) $\mathrm{Sp}(\Delta + r\delta_0) = [-4d, 0] \cup \{\mu(r)\}$ *with*

$$\mu(r) \begin{cases} = 0, & \text{if } 0 < r \le r_d, \\ > 0, & \text{if } r > r_d; \end{cases} \tag{1.15}$$

(ii) *for $r > r_d$, $\mu(r)$ is the unique solution of the equation $R(\mu) = 1/r$ and is an eigenvalue corresponding to a strictly positive eigenfunction;*

(iii) *on $(r_d, \infty)$, $r \mapsto \mu(r)/r$ is strictly increasing with $\lim_{r\to\infty} \mu(r)/r = 1$;*

(iv) *on $(0,\infty)$, $r \mapsto \mu(r)$ is convex.*

Our first theorem establishes the existence of the Lyapunov exponents $\widehat{\lambda}_p, \lambda_p$ and identifies $\widehat{\lambda}_p$.

THEOREM 1.4. *Let $p \in \mathbb{N}$.*

(i) *If $d \ge 1$, then the limit $\widehat{\lambda}_p$ exists, is finite and equals $\widehat{\lambda}_p = \rho\mu(p\gamma/\rho)$.*

(ii) *If $d \ge 3$ and $0 < p\gamma/\rho < r_d$, then the limit $\lambda_p$ exists and is finite.*

(iii) *If $d \ge 3$ and $p\gamma/\rho = r_d$, then the limit $\lambda_p$ exists and is infinite.*

Note from (1.15) that $\widehat{\lambda}_p > 0$ when either (a) $d = 1, 2$ or (b) $d \ge 3$ and $p\gamma/\rho > r_d$. Consequently, $\lambda_p = \infty$ in that regime.

Our second theorem addresses the $\kappa$-dependence of $\lambda_p = \lambda_p(\kappa)$ in the regime where it is finite. In order to state this theorem, we define, for $d = 3$,

$$\mathcal{P} = \sup_{\substack{f \in H^1(\mathbb{R}^3) \\ \|f\|_2 = 1}} [\|(-\Delta_{\mathbb{R}^3})^{-1/2} f^2\|_2^2 - \|\nabla_{\mathbb{R}^3} f\|_2^2] \in (0,\infty), \tag{1.16}$$

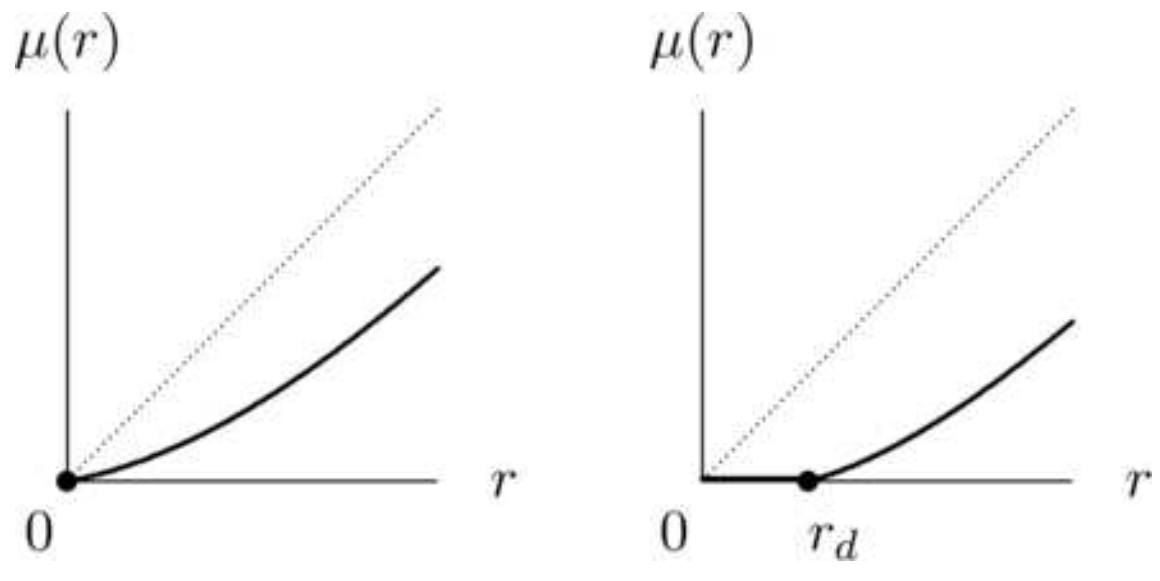

FIG. 1. $r \mapsto \mu(r)$ *for $d = 1, 2$, respectively, $d \ge 3$.*

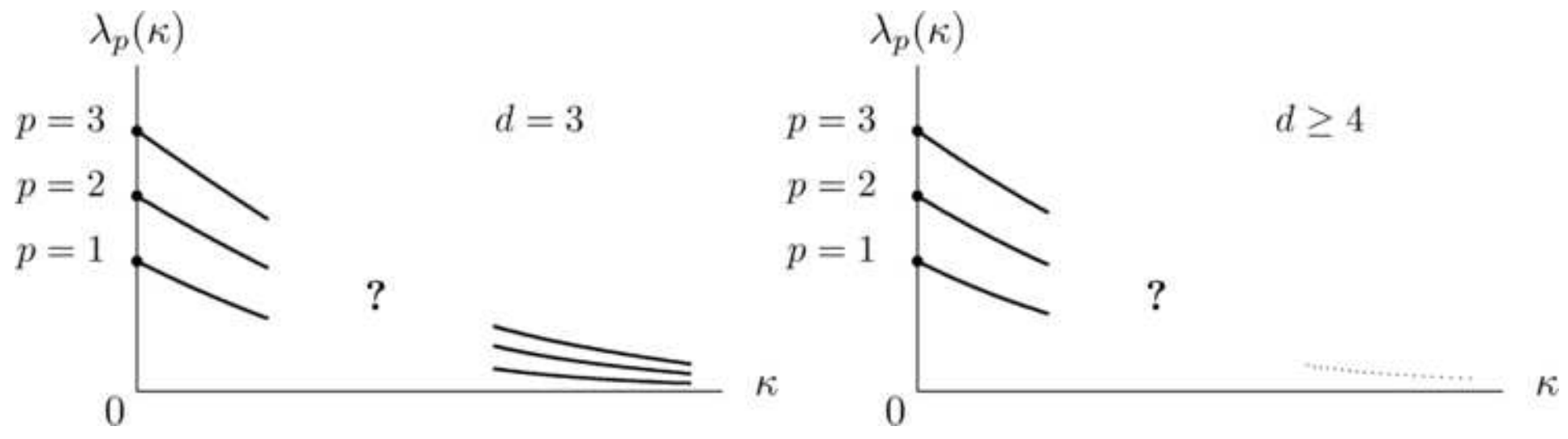

FIG. 2. *Qualitative picture of $\kappa \mapsto \lambda_p(\kappa)$. The dotted line represents the asymptotics for $d \geq 4$ given by (1.19).*

where $\nabla_{\mathbb{R}^3}$ and $\Delta_{\mathbb{R}^3}$ are the continuous (!) gradient and Laplacian, $\|\cdot\|_2$ is the $L^2$-norm, $H^1(\mathbb{R}^3) = \{f : \mathbb{R}^3 \to \mathbb{R} : f, \nabla_{\mathbb{R}^3} f \in L^2(\mathbb{R}^3)\}$ and

$$\|(-\Delta_{\mathbb{R}^3})^{-1/2} f^2\|_2^2 = \int_{\mathbb{R}^3} dx\, f^2(x) \int_{\mathbb{R}^3} dy\, f^2(y) \frac{1}{4\pi|x-y|}. \tag{1.17}$$

THEOREM 1.5. *Let $p \in \mathbb{N}$, $d \geq 3$ and $0 < p\gamma/\rho < r_d$.*

(i) *On $[0,\infty)$, $\kappa \to \lambda_p(\kappa)$ is strictly decreasing, continuous and convex.*

(ii)

$$\lim_{\kappa \downarrow 0} \lambda_p(\kappa) = \lambda_p(0) = \nu\gamma \frac{p\gamma/\rho}{r_d - p\gamma/\rho}. \tag{1.18}$$

(iii)

$$\lim_{\kappa \to \infty} \kappa \lambda_p(\kappa) = \frac{\nu\gamma^2}{r_d} + \mathbb{1}_{\{d=3\}} \left(\frac{\nu\gamma^2}{\rho} p\right)^2 \mathcal{P}. \tag{1.19}$$

Note that the asymptotics as $\kappa \to \infty$ are the same for all $p$ when $d \geq 4$; the correction term with $\mathcal{P}$ is present only when $d = 3$ (see Figure 2).

Summarizing, we have the following behavior:

COROLLARY 1.6. *Let $p \in \mathbb{N}$.*

(i) *The system is strongly $p$-catalytic if and only if either of the following holds:*

- *$d = 1, 2$;*
- *$d \geq 3$ and $p\gamma/\rho > r_d$.*

(ii) *The system is strongly $p$-intermittent if any of the following holds:*

- *$d = 1, 2$;*
- *$d \geq 3$ and $p\gamma/\rho \geq r_d$;*
- *$d \geq 3$, $0 < p\gamma/\rho < r_d$ and $\kappa$ is sufficiently small;*
- *$d = 3$, $0 < p\gamma/\rho < r_3$ and $\kappa$ is sufficiently large.*

1.4. *Discussion.* Theorems 1.4 and 1.5 show that there is a delicate interplay between the various parameters in the model.

Catalytic behavior is controlled by $\gamma/\rho$, the ratio of the strength and the speed of the catalyst $\xi$, and is sensitive to this ratio only when $d \geq 3$. For large ratios, the system is strongly catalytic; for small ratios, the system is weakly catalytic. The high peaks in the reactant $u$ develop around those sites where the catalyst $\xi$ piles up. The analysis behind Theorem 1.4(i) shows that strongly catalytic behavior corresponds to the high peaks in the $u$-field being concentrated at *single sites*, whereas weakly catalytic plus strongly intermittent behavior corresponds to the high peaks being *spread out over islands containing several sites* (weakly intermittent behavior corresponds to the presence of no relevant high peaks). It follows from Lemma 1.3 and Theorem 1.4(i) that $\rho \mapsto \widehat{\lambda}_p(\rho)$ is strictly decreasing in the strongly catalytic regime. Thus, as the catalyst $\xi$ moves faster, it is less effective. Moreover, $\lim_{\rho\downarrow 0} \widehat{\lambda}_p(\rho) = p\gamma$. Note that $\kappa$, the speed of the reactant $u$, plays no role, nor does $\nu$, the intensity of the catalyst $\xi$.

Intermittency has the following interpretation. Consider the situation where the system is strongly $p$-intermittent, that is, $\lambda_{p-1} < \lambda_p$. Pick any $\alpha \in (\lambda_{p-1}, \lambda_p)$. Then, on the one hand, the density of the point process

$$\Gamma_t = \{x \in \mathbb{Z}^d : u(x,t) > e^{\alpha t}\} \tag{1.20}$$

of high exceedances of the solution $u$ tends to zero exponentially fast as $t \to \infty$. On the other hand,

$$\langle u(0,t)^p \rangle \sim \langle u(0,t)^p \mathbb{1}_{\{u(0,t) > e^{\alpha t}\}} \rangle, \qquad t \to \infty, \tag{1.21}$$

and, therefore, by the ergodic theorem,

$$\frac{1}{|V_t|} \sum_{x \in V_t} u(x,t)^p \sim \frac{1}{|V_t|} \sum_{x \in V_t \cap \Gamma_t} u(x,t)^p, \qquad t \to \infty, \tag{1.22}$$

provided the centered boxes $V_t$ exhaust $\mathbb{Z}^d$ sufficiently fast. For details, we refer to Gärtner and König [15], Section 1.3. Thus, $p$-intermittency means that the $p$th moment of the solution is asymptotically "concentrated" on a thin set $\Gamma_t$ of high exceedances (which is expected to consist of "islands" that are located far from each other).

Intermittent behavior is sensitive to the parameters only when $d \geq 3$. Theorem 1.5(ii) shows that for small $\kappa$, the reactant $u$ has a range of high peaks that grow at different exponential rates and determine the successive moments, and so the system is strongly intermittent. For large $\kappa$, on the other hand, the behavior depends on the dimension. The large diffusion of the reactant $u$ prevents it from easily localizing around the high peaks where the catalyst $\xi$ piles up. As is clear from Theorem 1.5(iii), in $d = 3$, the system is strongly intermittent also for large $\kappa$, while in $d \geq 4$, it may or may not be. To decide this issue, we need finer asymptotics than those provided by (1.19). We conjecture the following.

Conjecture 1.7. *In $d=3$, the system is strongly p-intermittent for all $\kappa$.*

Conjecture 1.8. *For $d \geq d_0 \geq 4$, the system is weakly p-intermittent for $\kappa \geq \kappa_0(p)$.*

As promised at the end of Section 1.1, we discuss the results obtained by Kesten and Sidoravicius [23].

I. $d=1,2$: *For any choice of the parameters, the average number of B-particles per site tends to infinity at a rate faster than exponential.* This result is covered by our Theorem 1.4(i), because the inclusion of the death rate $\delta$ shifts $\lambda_1$ by $-\delta$ [recall (1.8)], but does not affect $\widehat{\lambda}_1$, while $\widehat{\lambda}_1 > 0$ in $d=1,2$ for any choice of the parameters.

II. $d \geq 3$: *For $\gamma$ sufficiently small and $\delta$ sufficiently large, the average number of B-particles per site tends to zero exponentially fast.* This result is covered by our Theorem 1.4(ii), because small $\gamma$ corresponds to the weakly catalytic regime for which $0 < \lambda_1 < \infty$ so that exponentially fast extinction occurs when $\delta > \lambda_1$.

III. $d \geq 1$: *For $\gamma$ sufficiently large, conditioned on the evolution of the A-particles, there is a phase transition: namely, for small $\delta$, the B-particles locally survive, while for large $\delta$ they become locally extinct.* This result is not linked to our theorems because we have no information on the quenched Lyapunov exponent.

The main focus of Kesten and Sidoravicius [23] is on survival versus extinction, while our focus is on moment asymptotics. Their approach does not lead to the identification of Lyapunov exponents, but it is more robust under an adaptation of the model than our approach, which is based on the Feynman–Kac representation in Section 2.1.

For related work on catalytic branching models, focusing in particular on continuum models with a singular catalyst in a measure-valued context, we refer to the overview papers by Dawson and Fleischmann [8] and Klenke [24]. Related references can also be found therein.

1.5. *Heuristics behind the asymptotics as $\kappa \to \infty$.* In this section, we summarize the main steps in the proof of Theorem 1.5(iii) in Sections 5–8. For simplicity, we restrict to the case $p=1$.

We will see that after a time scaling $t \mapsto t/\kappa$, the Feynman–Kac representation of the first moment (see Section 2.1) attains the form

$$\langle u(0,t/\kappa)\rangle = e^{\nu\gamma(t/\kappa)}\mathbb{E}_0^X\left(\exp\left[\frac{\nu\gamma}{\kappa}\int_0^t w^*(X(s),s)\,ds\right]\right), \tag{1.23}$$

where $X$ is simple random walk on $\mathbb{Z}^d$ (with generator $\Delta$) starting at the origin and $w^*$ denotes the solution of the random parabolic equation

$$\frac{\partial}{\partial t} w^* = \frac{\rho}{\kappa}\Delta w^* + \frac{\gamma}{\kappa}\delta_{X(t)}(1+w^*) \tag{1.24}$$

with zero initial condition. A serious complication is the long-range dependence of $w^*(\cdot,t)$ on the past $X(s)$, $s\in[0,t]$. For large $\kappa$, however, $w^*$ is small and, consequently, the $w^*$-term after the Kronecker symbol in the right-hand side of (1.24) is negligible. Therefore,

$$w^*(X(s),s) \approx \frac{\gamma}{\kappa}\int_0^s du\, p_{\rho/\kappa}(X(s)-X(u), s-u), \tag{1.25}$$

where $p_{\rho/\kappa}$ denotes the transition kernel of simple random walk with diffusion constant $\rho/\kappa$. Hence, the computation of

$$\lim_{\kappa\to\infty} \kappa\lambda_1(\kappa) = \lim_{\kappa\to\infty}\lim_{t\to\infty}\frac{\kappa^2}{t}\log(e^{-\nu\gamma(t/\kappa)}\langle u(0,t/\kappa)\rangle) \tag{1.26}$$

reduces to the asymptotic investigation of

$$\frac{\kappa^2}{t}\log \mathbb{E}_0^X\left(\exp\left[\frac{\nu\gamma^2}{\kappa^2}\int_0^t ds\int_s^t du\, p_{\rho/\kappa}(X(u)-X(s),u-s)\right]\right) \tag{1.27}$$

when first letting $t\to\infty$ and then $\kappa\to\infty$.

We split the inner integral into three parts by separately integrating over the time intervals $[s,s+\varepsilon\kappa^3]$, $[s+\varepsilon\kappa^3, s+K\kappa^3]$ and $[s+K\kappa^3,t]$, with $\varepsilon$ and $K$ being a small (resp., a large) constant. Through rough bounds, the third term turns out to be negligible. In $d\geq 4$, the same is true for the second term. We then show that a law of large numbers acts on the first term, that is, for large $\kappa$, the corresponding expression in the exponent may be replaced by its expectation. The lower bound is obvious from Jensen's inequality, but the proof of the upper bound turns out to be highly nontrivial. We have

$$\begin{aligned}&\frac{\kappa^2}{t}\mathbb{E}_0^X\left(\frac{\nu\gamma^2}{\kappa^2}\int_0^t ds\int_s^{s+\varepsilon\kappa^3} du\, p_{\rho/\kappa}(X(u)-X(s),u-s)\right)\\ &\qquad = \nu\gamma^2\int_0^{\varepsilon\kappa^3} du\, p_{1+\rho/\kappa}(0,u).\end{aligned} \tag{1.28}$$

As $\kappa\to\infty$, the integral in the right-hand side converges to $1/r_d$, the value of the Green function at 0 associated with $\Delta$. This yields assertion (1.19) for $d\geq 4$ and the first part of the desired expression for $d=3$.

In $d=3$, the first and second terms in the exponent of (1.27), as obtained via the above splitting, may be separated from each other with the help of Hölder's inequality (with a large exponent for the first factor and an exponent close to one for the second factor). Hence, for $d=3$, it only remains to consider the asymptotics of the second term as $t \to \infty$ and $\kappa \to \infty$ (in this order). After a Gaussian approximation of the transition kernel, this leads to the study of

$$\begin{aligned}(1.29)\qquad &\frac{\kappa^2}{t} \log \mathbb{E}_0^X \biggl( \exp \biggl[ \frac{\nu\gamma^2}{\kappa} \int_0^{t/\kappa^2} ds \int_{s+\varepsilon\kappa}^{s+K\kappa} du \\ &\qquad\qquad \times p_G\Bigl( X^\kappa(u) - X^\kappa(s), \frac{\rho}{\kappa}(u-s) \Bigr) \biggr] \biggr),\end{aligned}$$

where $p_G(x,t) = (4\pi t)^{-3/2} \exp[-\|x\|^2/4t]$ and $X^\kappa(\cdot) = X(\kappa^2 \cdot)/\kappa$ approaches Brownian motion as $\kappa \to \infty$. Next, observe that $(\rho/\kappa)(u-s)$ stays nearly constant when $u$ and $s$ with $u - s \geq \varepsilon\kappa$ vary over time intervals of length $\delta\kappa$ with $0 < \delta \ll \varepsilon$. But, as $\kappa \to \infty$, on each such time interval, we may apply the large deviation principle for the occupation time measure of $X^\kappa$. Then an application of the Laplace–Varadhan method yields that, for large $\kappa$ and $t \gg \kappa^3$, the expression in (1.29) behaves like

$$\begin{aligned}(1.30)\qquad &\frac{\kappa^2}{t} \sup_{\mu_{(\cdot)}} \biggl\{ \frac{\nu\gamma^2}{\kappa} \int_0^{t/\kappa^2} ds \int_{s+\varepsilon\kappa}^{s+K\kappa} du \int_{\mathbb{R}^3} \mu_s(dx) \int_{\mathbb{R}^3} \mu_u(dy) p_G\Bigl( y - x, \frac{\rho}{\kappa}(u-s) \Bigr) \\ &\qquad\qquad - \int_0^{t/\kappa^2} ds I(\mu_s) \biggr\},\end{aligned}$$

where $I$ denotes the large deviation rate function for the occupation time measure and the supremum is taken over (probability) measure-valued paths $\mu_{(\cdot)}$ on the time interval $[0, t/\kappa^2]$. It turns out that this supremum is attained for a time-independent path. Hence, (1.30) coincides with

$$(1.31)\qquad \sup_{\mu} \biggl\{ \frac{\nu\gamma^2}{\rho} \int_{\mathbb{R}^3} \mu(dx) \int_{\mathbb{R}^3} \mu(dy) \int_{\varepsilon}^{K} du\, p_G(y-x, u) - I(\mu) \biggr\}.$$

Finally, by letting $\varepsilon \to 0$ and $K \to \infty$, we see that the last integral approaches the Green function and the whole expression becomes

$$(1.32)\qquad \sup_{\mu} \biggl\{ \frac{\nu\gamma^2}{\rho} \int_{\mathbb{R}^3} \mu(dx) \int_{\mathbb{R}^3} \mu(dy) \frac{1}{4\pi|y-x|} - I(\mu) \biggr\}.$$

Since

$$(1.33)\quad I(\mu) = \begin{cases} \|\nabla_{\mathbb{R}^3} f\|_2^2, & \text{for } \mu(dx) = f^2(x)\, dx, \qquad f \in H^1(\mathbb{R}^3), \\ \infty, & \text{otherwise,} \end{cases}$$

(1.32) is easily seen to coincide with $(\nu\gamma^2/\rho)^2 \mathcal{P}$, where the variational expression for $\mathcal{P}$ is given by (1.16)–(1.17). In this way, we arrive at the second part of the expression in the right-hand side of (1.19) for $p = 1$ and $d = 3$, and we are done.

Interestingly, (1.16) is precisely the variational problem that arises in the so-called *polaron model*. Here, one takes Brownian motion $W$ on $\mathbb{R}^3$ with generator $\Delta_{\mathbb{R}^3}$, starting at the origin and, for $\alpha > 0$, considers the quantity

$$\begin{aligned}\Theta(t;\alpha) &= \frac{1}{\alpha^2 t}\log \mathbb{E}_0^W\left(\exp\left[\alpha\int_0^t ds\int_s^t du\frac{e^{-(u-s)}}{|W(u)-W(s)|}\right]\right)\\ &= \frac{1}{\alpha^2 t}\log \mathbb{E}_0^W\left(\exp\left[\frac{1}{\alpha^2}\int_0^{\alpha^2 t} ds\int_s^{\alpha^2 t} du\frac{e^{-(u-s)/\alpha^2}}{|W(u)-W(s)|}\right]\right).\end{aligned} \tag{1.34}$$

It was shown by Donsker and Varadhan [10] that

$$\theta(\alpha) = \lim_{t\to\infty}\Theta(t;\alpha), \qquad \alpha > 0, \tag{1.35}$$

exists and

$$\lim_{\alpha\to\infty}\theta(\alpha) = 4\sqrt{\pi}\mathcal{P}. \tag{1.36}$$

The expression obtained by substituting $\alpha^2 = \kappa/\rho$ and replacing $t$ by $\rho t/\kappa^3$ in the second line of (1.34) is qualitatively similar to (1.29). Although the two exponents are not the same, it turns out that they have the same large deviation behavior for $t\to\infty$ and $\kappa\to\infty$. Details can be found in Sections 5 and 7.

While Donsker and Varadhan use large deviations on the level of the process, we use large deviations on the level of the occupation time measure associated with the process.

It was shown by Lieb [25] that (1.16) has a unique maximizer modulo translations and that the centered maximizer is radially symmetric, radially nonincreasing, strictly positive and smooth.

1.6. *Future challenges.* One challenge is to understand the geometry and location of the high peaks in the $u$-field that determine the Lyapunov exponents in the weakly catalytic regime. These peaks (which are spread out over islands containing several sites) move and grow with time; the question is *how*.

Another challenge is to compute the quenched Lyapunov exponent, that is,

$$\lambda = \lim_{t\to\infty}\frac{1}{t}\log u(0,t), \qquad \xi\text{-}a.s., \tag{1.37}$$

and to study its dependence on the parameters.

Finally, the choice in (1.4) constitutes one of the simplest types of catalyst dynamics. What happens for other choices of the $\xi$-field, for example, when $\xi(x,t)$ is $\gamma$ times the occupation number at site $x$ at time $t$ of a system of particles performing a simple symmetric exclusion process in equilibrium (i.e., particles moving like symmetric random walks but not being allowed to occupy the same site)? This extension, which constitutes one of the simplest examples of a catalyst with interaction, will be addressed in Gärtner, den Hollander and Maillard [13]. Since particles cannot pile up in this model, there is no strongly catalytic regime (i.e., $\widehat{\lambda}_p = 0$). However, it turns out that the weakly catalytic regime again exhibits a delicate interplay of parameters controlling the intermittent behavior.

The asymptotic behavior for large $\kappa$ may be expected to be universal, that is, to some extent independent of the details of the dynamics of the catalysts. In fact, we will see evidence of this in [13].

1.7. *Outline.* We now outline the rest of this paper. In Section 2, we formulate some preparatory results, including a Feynman–Kac representation for the moments of the solution of (1.1) under (1.4)–(1.6), a certain concentration estimate, and the proof of Lemma 1.3. In Section 3, we prove Theorem 1.4(i) for $\widehat{\lambda}_p$. Section 4 contains the proof of Theorems 1.4(ii), (iii) and 1.5(i), (ii) for $\lambda_p = \lambda_p(\kappa)$ in three parts: existence, convexity and behavior for small $\kappa$. Sections 5–8, which take up over half of the paper, contain the proof of Theorem 1.5(iii): behavior for large $\kappa$.

## 2. Preparations.

Section 2.1 contains a Feynman–Kac representation for the moments of $u(0,t)$ that serves as the starting point of our analysis. Section 2.2 derives a certain concentration estimate that is needed for the proof of Theorem 1.4(i), while Section 2.3 contains the proof of Lemma 1.3.

2.1. *Feynman–Kac representation.* The formal starting point of our analysis of (1.1) is the following Feynman–Kac representation for the $p$th moment of the $u$-field.

PROPOSITION 2.1. *For any $p \in \mathbb{N}$,*

$$\langle u(0,t)^p\rangle = e^{p\nu\gamma t}\mathbb{E}_{0,\dots,0}^{X_1,\dots,X_p}\left(\exp\left[\nu\gamma\int_0^t \sum_{q=1}^p w(X_q(s),s)\,ds\right]\right), \tag{2.1}$$

*where $X_1,\dots,X_p$ are independent simple random walks on $\mathbb{Z}^d$ with step rate $2d\kappa$ starting from the origin. The expectation is taken with respect to these random walks and $w:\mathbb{Z}^d\times[0,\infty)\to\mathbb{R}$ is the solution of the Cauchy problem*

$$\begin{aligned}\frac{\partial}{\partial t}w(x,t) &= \rho\Delta w(x,t) + \gamma\left[\sum_{q=1}^p \delta_{X_q(t)}(x)\right]\{w(x,t)+1\},\\ w(\cdot,0) &\equiv 0.\end{aligned} \tag{2.2}$$

PROOF. We give the proof for $p=1$. Let $X,Y$ be independent copies of $X_1,Y_1$ [recall (1.5)]. By applying the Feynman–Kac formula to (1.1) and (1.6) and inserting (1.4), we have

$$\begin{aligned} u(0,t) &= \mathbb{E}_0^X\left(\exp\left[\int_0^t \xi(X(s),t-s)\,ds\right]\right) \\ &= \mathbb{E}_0^X\left(\prod_k \exp\left[\gamma\int_0^t \delta_{Y_k(t-s)}(X(s))\,ds\right]\right). \end{aligned} \tag{2.3}$$

Next, we take the expectation over the $\xi$-field. This is done by first taking the expectation over the trajectories $Y_k$, given the starting points $Y_k(0)$, and then taking the expectation over $Y_k(0)$ according to a Poisson random field with intensity $\nu$:

$$\begin{aligned} \langle u(0,t)\rangle &= \left\langle \mathbb{E}_0^X \prod_k \mathbb{E}_{Y_k(0)}^{Y_k}\left(\exp\left[\gamma\int_0^t \delta_{Y_k(t-s)}(X(s))\,ds\right]\right)\right\rangle \\ &= \mathbb{E}_0^X\left\langle \prod_k v(Y_k(0),t)\right\rangle \\ &= \mathbb{E}_0^X\left(\prod_{y\in\mathbb{Z}^d}\sum_{n\in\mathbb{N}_0}\frac{[\nu v(y,t)]^n}{n!}e^{-\nu}\right) \\ &= \mathbb{E}_0^X\left(\prod_{y\in\mathbb{Z}^d}\exp[\nu\{v(y,t)-1\}]\right) \end{aligned} \tag{2.4}$$

(where $\mathbb{N}_0=\mathbb{N}\cup\{0\}$) with

$$v(y,t) = \mathbb{E}_y^Y\left(\exp\left[\gamma\int_0^t \delta_{Y(t-s)}(X(s))\,ds\right]\right). \tag{2.5}$$

The latter is a functional of $X$ and is the solution of the Cauchy problem

$$\frac{\partial}{\partial t}v(x,t) = \rho\Delta v(x,t) + \gamma\delta_{X(t)}(x)v(x,t), \qquad v(\cdot,0)\equiv 1. \tag{2.6}$$

The last expectation in the right-hand side of (2.4) equals $\mathbb{E}_0^X(\exp[\nu\Sigma(t)])$ with $\Sigma(t)=\sum_{y\in\mathbb{Z}^d}\{v(y,t)-1\}$. But, from (2.6), we see that

$$\frac{d}{dt}\Sigma(t) = 0 + \gamma v(X(t),t), \qquad \Sigma(0)=0. \tag{2.7}$$

Hence, $\Sigma(t)=\gamma\int_0^t v(X(s),s)\,ds$. Now, put

$$w(x,t) = v(x,t) - 1 \tag{2.8}$$

to complete the proof. The extension to arbitrary $p$ is straightforwardly achieved by taking $p$ independent copies of the random walk $X$ (rather than one) and repeating the argument. □

It follows from (1.9) and Proposition 2.1 that

$$\Lambda_p(t) = \frac{1}{pt} \log \mathbb{E}_{0,\dots,0}^{X_1,\dots,X_p}\left(\exp\left[\nu\gamma \int_0^t \sum_{q=1}^p w(X_q(s),s)\,ds\right]\right). \tag{2.9}$$

This is the representation we will work with later. Note that

$$w = w_{X_1,\dots,X_p}, \tag{2.10}$$

that is, $w(\cdot,t)$ is to be solved as a function of the trajectories $X_1,\dots,X_p$ up to time $t$ (and of the parameters $p,\gamma,\rho$) and $\Lambda_p(t)$ is to be calculated after insertion of the solution into the Feynman–Kac representation (2.9). Thus, the study of $\Lambda_p(t)$ amounts to carrying out a *large deviation analysis for a time-inhomogeneous functional of $p$ random walks having long-time correlations.*

Note that

$$w(x,t) > 0 \qquad \forall\, x \in \mathbb{Z}^d,\ t > 0, \tag{2.11}$$

as can be seen from (2.2). Hence, $t \mapsto t\Lambda_p(t)$ is strictly positive and strictly increasing on $(0,\infty)$, as was claimed in Section 1.2.

2.2. *Concentration estimate.* The following estimate will be needed later on. It shows that the solution of (2.2) is maximal when $X_1,\dots,X_p$ stay at the origin.

PROPOSITION 2.2. *For any $p \in \mathbb{N}$ and $X_1,\dots,X_p$,*

$$w(x,t) \le \bar{w}(0,t) \qquad \forall\, x \in \mathbb{Z}^d,\ t \ge 0, \tag{2.12}$$

*where $\bar{w}: \mathbb{Z}^d \times [0,\infty) \to \mathbb{R}$ is the solution of the Cauchy problem*

$$\frac{\partial}{\partial t}\bar{w}(x,t) = \rho\Delta\bar{w}(x,t) + p\gamma\delta_0(x)\{\bar{w}(x,t)+1\}, \qquad \bar{w}(\cdot,0) \equiv 0. \tag{2.13}$$

PROOF. Recall (1.11). Abbreviate $\oint dk = (2\pi)^{-d}\int_{[-\pi,\pi)^d} dk$. Let

$$p_\rho(x,t) = \oint dk\, e^{-\rho t\widehat{\varphi}(k)} e^{-ik\cdot x}, \qquad x \in \mathbb{Z}^d,\ t \ge 0, \tag{2.14}$$

denote the Fourier representation of the transition kernel associated with $\rho\Delta$. From this representation, we see that

$$\max_{x\in\mathbb{Z}^d} p_\rho(x,t) = p_\rho(0,t) \qquad \forall\, t \ge 0. \tag{2.15}$$

The solution of (2.2) has the (implicit) representation

$$w(x,t) = \gamma \sum_{q=1}^{p} \int_0^t ds\, p_\rho(x - X_q(s), t-s)\{w(X_q(s), s) + 1\}. \tag{2.16}$$

Abbreviate

$$\widehat{\eta}(t) = \frac{1}{p} \sum_{q=1}^{p} w(X_q(t), t). \tag{2.17}$$

We first prove that

$$\widehat{\eta}(t) \le \bar{w}(0,t) \qquad \forall\, t \ge 0. \tag{2.18}$$

To that end, take $x = X_r(t)$, $r = 1, \dots, p$, in (2.16), sum over $r$ and use (2.15), to obtain

$$\widehat{\eta}(t) \le p\gamma \int_0^t ds\, p_\rho(0, t-s)\{\widehat{\eta}(s) + 1\}. \tag{2.19}$$

Define

$$h(t) = p\gamma p_\rho(0,t) \ge 0. \tag{2.20}$$

Then (2.19) can be rewritten as

$$\widehat{\eta} \le h * \{\widehat{\eta} + 1\}. \tag{2.21}$$

Next, put

$$\bar{\eta}(t) = \bar{w}(0,t). \tag{2.22}$$

Then the same formulas with $X_1(\cdot), \dots, X_p(\cdot) \equiv 0$ yield the relation

$$\bar{\eta} = h * \{\bar{\eta} + 1\}. \tag{2.23}$$

Thus, it remains to be shown that (2.21) and (2.23) imply (2.18), that is,

$$\widehat{\eta} \le \bar{\eta}. \tag{2.24}$$

This is achieved as follows.

Let $\delta = \bar{\eta} - \widehat{\eta}$. Then (2.21) and (2.23) give

$$\delta \ge h * \delta. \tag{2.25}$$

Iteration gives $\delta \ge h^{*n} * \delta$ and so, to prove (2.24), it suffices to show that $h^{*n}$ tends to zero as $n \to \infty$, uniformly on compact time intervals. To that end, put $h_T = \max_{t \in [0,T]} h(t)$. Then

$$0 \le h^{*n}(t) \le h_T \int_0^t h^{*(n-1)}(s)\, ds, \qquad t \in [0,T], \tag{2.26}$$

which, when iterated, gives

$$0 \le h^{*n}(t) \le h_T^n \frac{t^{n-1}}{(n-1)!}, \qquad t \in [0,T]. \tag{2.27}$$

Letting $n \to \infty$, we obtain the claimed assertion.

Finally, put

$$\eta(t) = \max_{x \in \mathbb{Z}^d} w(x,t), \qquad t \ge 0. \tag{2.28}$$

Then (2.15)–(2.17) and (2.24) give

$$\eta \le h * \{\widehat{\eta} + 1\} \le h * \{\bar{\eta} + 1\}. \tag{2.29}$$

Now, use (2.23) to get

$$\eta \le \bar{\eta}, \tag{2.30}$$

which, via (2.28), implies (2.12), as desired. □

PROPOSITION 2.3. *For any $p \in \mathbb{N}$, $t \mapsto \bar{w}(0,t)$ is nondecreasing and $\bar{w}(0) = \lim_{t\to\infty} \bar{w}(0,t)$ satisfies*

$$\bar{w}(0) = \begin{cases} \dfrac{p\gamma/\rho}{r_d - p\gamma/\rho}, & \text{if } 0 < \dfrac{p\gamma}{\rho} < r_d, \\ \infty, & \text{otherwise.} \end{cases} \tag{2.31}$$

PROOF. Returning to (2.22) and (2.23), and recalling (2.20), we have

$$\bar{w}(0,t) = p\gamma \int_0^t ds\, p_\rho(0,s)\{\bar{w}(0,t-s) + 1\}. \tag{2.32}$$

From this, we see that $t \mapsto \bar{w}(0,t)$ is nondecreasing. Using this fact in (2.32), we have

$$\bar{w}(0,t) \le p\gamma \left(\int_0^\infty ds p_\rho(0,s)\right)\{\bar{w}(0,t) + 1\} = \frac{p\gamma}{\rho}\frac{1}{r_d}\{\bar{w}(0,t) + 1\} \tag{2.33}$$

[recall (1.13)] and, hence,

$$\bar{w}(0,t) \le \mathrm{rhs}(2.31). \tag{2.34}$$

Taking the limit $t \to \infty$ in (2.32) and using monotone convergence, we get

$$\bar{w}(0) = p\gamma \left(\int_0^\infty du\, p_\rho(0,u)\right)\{\bar{w}(0) + 1\} = \frac{p\gamma}{\rho}\frac{1}{r_d}\{\bar{w}(0) + 1\}, \tag{2.35}$$

which implies the truth of the claimed assertion. □

2.3. *Proof of Lemma* 1.3. The proof is elementary.

(i)–(ii) For $r \in (0,\infty)$, let $\mathcal{H} = \Delta + r\delta_0$. This is a self-adjoint operator on $\ell^2(\mathbb{Z}^d)$. Let $\widehat{v}(k) = \sum_{x\in\mathbb{Z}^d} e^{ik\cdot x} v(x)$ denote the Fourier transform of $v \in \ell^2(\mathbb{Z}^d)$. The Fourier transform of $\mathcal{H}$ is the operator on $L^2([-\pi,\pi)^d)$ given by

$$(\widehat{\mathcal{H}}\widehat{v})(k) = -\widehat{\varphi}(k)\widehat{v}(k) + r \oint \widehat{v}(l)\, dl, \tag{2.36}$$

where we recall (1.11). Since $\mathrm{Sp}(\mathcal{H}) = \mathrm{Sp}(\widehat{\mathcal{H}})$, (1.14) reads as

$$\mu(r) = \sup \mathrm{Sp}(\widehat{\mathcal{H}}). \tag{2.37}$$

The spectrum of $\widehat{\mathcal{H}}$ consists of those $\lambda \in \mathbb{R}$ for which $\lambda - \widehat{\mathcal{H}}$ is not invertible. Consider, therefore, the equation

$$(\lambda - \widehat{\mathcal{H}}) f = g. \tag{2.38}$$

Substituting (2.36) into (2.38), we get

$$(\lambda + \widehat{\varphi}) f - r \oint f = g. \tag{2.39}$$

Now, the range of $\widehat{\varphi}$ is the interval $[0,4d]$. Thus, if $\lambda \in [-4d,0]$, then there exists $g \in L^2([-\pi,\pi)^d)$ for which (2.39), and hence (2.38), has no solution, that is,

$$\mathrm{Sp}(\widehat{\mathcal{H}}) \supset [-4d,0]. \tag{2.40}$$

Next, assume that $\lambda > 0$. Divide (2.38) by $\lambda + \widehat{\varphi}$ and integrate to get

$$[1 - rR(\lambda)] \oint f = \oint \frac{g}{\lambda + \widehat{\varphi}} \tag{2.41}$$

with $R$ as defined in (1.12). If $rR(\lambda) = 1$, then there is, again, no solution, that is,

$$rR(\lambda) = 1 \quad \Longrightarrow \quad \lambda \in \mathrm{Sp}(\widehat{\mathcal{H}}). \tag{2.42}$$

If, on the other hand, $rR(\lambda) \neq 1$, then (2.41) yields a unique solution

$$f = \frac{1}{\lambda + \widehat{\varphi}} \left( g + \frac{r}{1 - rR(\lambda)} \oint \frac{g}{\lambda + \widehat{\varphi}} \right), \tag{2.43}$$

which is in $L^2([-\pi,\pi)^d)$, that is,

$$rR(\lambda) \neq 1 \quad \Longrightarrow \quad \lambda \notin \mathrm{Sp}(\widehat{\mathcal{H}}). \tag{2.44}$$

The same argument shows that

$$(-\infty,-4d) \cap \mathrm{Sp}(\widehat{\mathcal{H}}) = \varnothing. \tag{2.45}$$

Combining (2.40), (2.42), (2.44) and (2.45), and noting that $rR(\lambda)=1$ has a unique solution $\lambda=\mu(r)>0$ if and only if $r>r_d$, we obtain assertions (i) and (ii). Note that if $r>r_d$, then

$$e=r(\mu(r)-\Delta)^{-1}\delta_0 \tag{2.46}$$

is a positive eigenfunction of $\mathcal{H}$ corresponding to the eigenvalue $\mu(r)$, normalized by $e(0)=1$ (rather than by $\|e\|_2=1$ with $\|\cdot\|_2$ the $\ell^2$-norm).

(iii) From (1.12), we have

$$\mu R(\mu)=\oint\frac{\mu}{\mu+\widehat{\varphi}}. \tag{2.47}$$

Differentiate this relation w.r.t. $\mu$ to obtain

$$[\mu R(\mu)]'=\oint\frac{\widehat{\varphi}}{(\mu+\widehat{\varphi})^2}>0. \tag{2.48}$$

Next, differentiate the relation $rR(\mu(r))=1$ w.r.t. $r$ and use the fact that $R'<0$ to obtain

$$\mu'(r)=-\frac{R(\mu(r))}{rR'(\mu(r))}>0. \tag{2.49}$$

From (2.48) and (2.49), we get

$$[\mu(r)/r]'=[\mu(r)R(\mu(r))]'=[\mu R(\mu)]'(r)\mu'(r)>0, \tag{2.50}$$

which proves the first part of assertion (iii). The second part of assertion (iii) follows from the estimate

$$0<\frac{1}{\mu}-R(\mu)=\oint\frac{\widehat{\varphi}}{\mu(\mu+\widehat{\varphi})}<\frac{1}{\mu^2}\oint\widehat{\varphi}=\frac{2d}{\mu^2} \tag{2.51}$$

after letting $\mu\to\infty$, corresponding to $r\to\infty$.

(iv) Differentiating (2.49) w.r.t. $r$, we obtain

$$\mu''(r)=\frac{R(\mu(r))}{r^2[R'(\mu(r))]^3}\{2[R'(\mu(r))]^2-R(\mu(r))R''(\mu(r))\}. \tag{2.52}$$

Using the integral representations of $R$, $R'$, $R''$ obtained from (2.47), we find that $R>0$ and $R'<0$ and, by an application of the Cauchy–Schwarz inequality, that the term between braces is $<0$. Hence $\mu''(r)>0$.

An alternative way of seeing (iii) and (iv) is via the Rayleigh–Ritz formula,

$$\mu(r)=\sup_{\substack{f\in\ell^2(\mathbb{Z}^d)\\ \|f\|_2=1}}\Bigg\{rf(0)-\tfrac{1}{2}\sum_{\substack{x,y\in\mathbb{Z}^d\\ \|x-y\|=1}}[f(x)-f(y)]^2\Bigg\}. \tag{2.53}$$

Indeed, this formula shows that $r\mapsto\mu(r)$ is a supremum of linear functions and is therefore convex. Moreover, it shows that $r\mapsto\mu(r)/r$ is nondecreasing and, since the supremum is attained when $r>r_d$, it, in fact, gives that $r\mapsto\mu(r)/r$ is strictly increasing on $(r_d,\infty)$ (and tends to 1 as $r\to\infty$).

**3. Proof of Theorem 1.4(i).** The proof uses spectral analysis.

3.1. *Upper and lower bounds.* Let $\mathcal{H} = \rho\Delta + p\gamma\delta_0$. This is a self-adjoint operator on $\ell^2(\mathbb{Z}^d)$. Equation (2.13) reads as

$$\frac{\partial}{\partial t}\bar{w} = \mathcal{H}\bar{w} + p\gamma\delta_0, \qquad \bar{w}(\cdot,0) \equiv 0. \tag{3.1}$$

By (1.14),

$$\sup \operatorname{Sp}(\mathcal{H}) = \rho\mu(p\gamma/\rho). \tag{3.2}$$

Suppose first that $\rho\mu(p\gamma/\rho) > 0$. Then, by Lemma 1.3, this is an eigenvalue of $\mathcal{H}$ corresponding to a strictly positive eigenfunction $e \in \ell^2(\mathbb{Z}^d)$ (normalized as $\|e\|_2 = 1$). From (2.9) and Proposition 2.2, we have

$$\begin{aligned} -2d\kappa + \nu\gamma \frac{1}{t}\int_0^t \bar{w}(0,s)\,ds &\le \Lambda_p(t;\kappa) \\ &\le \nu\gamma\frac{1}{t}\int_0^t \bar{w}(0,s)\,ds, \end{aligned} \tag{3.3}$$

where we use the fact that

$$\mathbb{P}_{0,\dots,0}^{X_1,\dots,X_p}(X_q(s) = 0 \quad \forall s \in [0,t] \quad \forall q = 1,\dots,p) = e^{-2d\kappa pt}. \tag{3.4}$$

From (3.1), we have

$$\bar{w}(\cdot,t) = p\gamma\int_0^t ds (e^{(t-s)\mathcal{H}}\delta_0)(\cdot). \tag{3.5}$$

Moreover, from the spectral representation of $e^{(t-s)\mathcal{H}}$ and (3.2), we have

$$e^{(t-s)\rho\mu(p\gamma/\rho)}\langle e,\delta_0\rangle \le \langle e^{(t-s)\mathcal{H}}\delta_0,\delta_0\rangle \le e^{(t-s)\rho\mu(p\gamma/\rho)}\|\delta_0\|_2^2. \tag{3.6}$$

Combining (3.3), (3.5) and (3.6), we arrive at

$$\widehat{\lambda}_p = \lim_{t\to\infty}\frac{1}{t}\log\Lambda_p(t;\kappa) = \rho\mu(p\gamma/\rho). \tag{3.7}$$

Next suppose that $\rho\mu(p\gamma/\rho) = 0$. Then the upper bound in (3.6) remains valid [despite the fact that no eigenfunction $e \in \ell^2(\mathbb{Z}^d)$ with eigenvalue 0 may exist] and so the limit equals zero.

**4. Proofs of Theorems 1.4(ii)–(iii) and 1.5(i)–(ii).** In Section 4.1, we prove Theorem 1.4(ii)–(iii) and in Sections 4.2–4.3 we prove Theorem 1.5(i)–(ii).

4.1. *Existence of* $\lambda_p$. We already know that $\lambda_p$ exists and is infinite in the strongly catalytic regime, that is, when $d=1,2$ or $d\geq 3$, $p\gamma/\rho > r_d$; see the remarks below Theorem 1.4(i). At the end of Section 4.3, we will see that the same is true at the boundary of the weakly catalytic regime, that is, when $d\geq 3$, $p\gamma/\rho = r_d$, as is claimed in Theorem 1.4(iii). The following lemma proves Theorem 1.4(ii):

LEMMA 4.1. *Let $d\geq 3$ and $p\in\mathbb{N}$. If $0<p\gamma/\rho<r_d$, then the limit $\lambda_p$ exists and is finite.*

PROOF. Fix $d\geq 3$ and $p\in\mathbb{N}$ and return to (2.3). We have

$$u(0,t)=\sum_{x\in\mathbb{Z}^d} S_x(t) \tag{4.1}$$

with

$$S_x(t)=\mathbb{E}_0^X\left(\exp\left[\int_0^t \xi(X(s),t-s)\,ds\right]\delta_x(X(t))\right). \tag{4.2}$$

Hence,

$$\begin{aligned}\langle u(0,t)^p\rangle &= \left\langle\left[\sum_{x\in\mathbb{Z}^d} S_x(t)\right]^p\right\rangle\\ &\leq \left\langle\left[\sum_{x\in Q_{t\log t}} S_x(t)\right]^p\right\rangle + p\left\langle\sum_{x_1\notin Q_{t\log t}} S_{x_1}(t)\left[\sum_{x\in\mathbb{Z}^d} S_x(t)\right]^{p-1}\right\rangle \qquad (4.3)\\ &= \left\langle\left[\sum_{x\in Q_{t\log t}} S_x(t)\right]^p\right\rangle + p\sum_{x_1\notin Q_{t\log t}}\ \sum_{x_2,\dots,x_p\in\mathbb{Z}^d}\left\langle\prod_{q=1}^p S_{x_q}(t)\right\rangle,\end{aligned}$$

where $Q_{t\log t}=[-t\log t,t\log t]^d\cap\mathbb{Z}^d$. By Jensen's inequality, the first term in the right-hand side of (4.3) is bounded above by

$$\begin{aligned}&|Q_{t\log t}|^{p-1}\left\langle\sum_{x\in Q_{t\log t}}[S_x(t)]^p\right\rangle\\ (4.4)\quad &= e^{p\nu\gamma t}|Q_{t\log t}|^{p-1}\\ &\quad\times\sum_{x\in Q_{t\log t}}\mathbb{E}_{0,\dots,0}^{X_1,\dots,X_p}\left(\exp\left[\nu\gamma\sum_{q=1}^p\int_0^t w(X_q(s),s)\,ds\right]\prod_{q=1}^p\delta_x(X_q(t))\right),\end{aligned}$$

where the last line follows the calculation in the proof of Proposition 2.1. The second term in the right-hand side of (4.3) is bounded above by

$$pe^{p\nu\gamma\{\bar{w}(0)+1\}t}\mathbb{P}_0^{X_1}(X_1(t)\notin Q_{t\log t}), \tag{4.5}$$

where we use the fact that $w(x,t) \le \bar{w}(0,t) \le \bar{w}(0)$ by Propositions 2.2 and 2.3, with $\bar{w}(0) < \infty$ strictly inside the weakly $p$-catalytic regime considered here. Now, define

$$\underline{\Lambda}_p(t) = \max_{x\in\mathbb{Z}^d} \frac{1}{pt} \log \mathbb{E}_{0,\dots,0}^{X_1,\dots,X_p}\left(\exp\left[\nu\gamma \sum_{q=1}^p \int_0^t w(X_q(s),s)\,ds\right] \times \prod_{q=1}^p \delta_x(X_q(t))\right). \tag{4.6}$$

Since the probability in (4.5) is superexponentially small (SES) in $t$, we see that a comparison of (2.9) and (4.6) yields the sandwich [combine (1.9) and (4.3)–(4.5)]

$$\begin{aligned} \underline{\Lambda}_p(t) &\le \Lambda_p(t) \\ &\le \frac{1}{pt}\log(|Q_{t\log t}|^p e^{pt\underline{\Lambda}_p(t)} + \mathrm{SES}), \end{aligned} \tag{4.7}$$

so that

$$\lim_{t\to\infty}[\Lambda_p(t) - \underline{\Lambda}_p(t)] = 0. \tag{4.8}$$

To prove existence of $\lambda_p$, it therefore suffices to prove existence of

$$\bar{\lambda}_p = \lim_{t\to\infty} \underline{\Lambda}_p(t), \tag{4.9}$$

after which we conclude that $\lambda_p = \bar{\lambda}_p$.

The proof of existence of (4.9) is achieved as follows. Write

$$w(x,s) = w_{X_1[0,t],\dots,X_p[0,t]}(x,s), \qquad s\in[0,t], \tag{4.10}$$

to exhibit the dependence of $w$ on the $p$ trajectories. We have, for any $s,t\ge 0$,

$$w_{X_1[0,s+t],\dots,X_p[0,s+t]}(x,u) \begin{cases} = w_{X_1[0,s],\dots,X_p[0,s]}(x,u), & \text{for } u\in[0,s], \\ \ge w_{X_1[s,s+t],\dots,X_p[s,s+t]}(x,u-s), & \text{for } u\in[s,s+t]. \end{cases} \tag{4.11}$$

Here, the inequality arises by resetting the initial condition to $\equiv 0$ at time $s$ and using the fact that the solution of (2.2) is monotone in the initial condition. It follows from (4.6) and (4.11) that

$$\begin{aligned} &p(s+t)\underline{\Lambda}_p(s+t) \\ &\quad\ge \max_{x,y\in\mathbb{Z}^d} \log \\ &\qquad \mathbb{E}_{0,\dots,0}^{X_1,\dots,X_p}\left(\exp\left[\nu\gamma\sum_{q=1}^p \int_0^{s+t} w_{X_1[0,s+t],\dots,X_p[0,s+t]}(X_q(u),u)\,du\right]\right. \end{aligned}$$

$$
\begin{aligned}
&\qquad\qquad \times \prod_{q=1}^{p} \delta_y(X_q(s)) \prod_{q=1}^{p} \delta_x(X_q(s+t))\Bigg) \\
&\geq \max_{x,y\in\mathbb{Z}^d} \Bigg\{ \log \mathbb{E}_{0,\ldots,0}^{X_1,\ldots,X_p} \Bigg( \exp\Bigg[\nu\gamma \sum_{q=1}^{p} \int_0^s w_{X_1[0,s],\ldots,X_p[0,s]}(X_q(u),u)\,du\Bigg] \\
&\qquad\qquad \times \prod_{q=1}^{p} \delta_y(X_q(s))\Bigg) \\
&\qquad + \log \mathbb{E}_{0,\ldots,0}^{X_1,\ldots,X_p} \Bigg( \exp\Bigg[\nu\gamma \sum_{q=1}^{p} \int_0^t w_{X_1[0,t],\ldots,X_p[0,t]}(X_q(u),u)\,du\Bigg] \\
&\qquad\qquad \times \prod_{q=1}^{p} \delta_{x-y}(X_q(t))\Bigg)\Bigg\} \\
&= ps\underline{\Lambda}_p(s) + pt\underline{\Lambda}_p(t),
\end{aligned} \tag{4.12}
$$

where we use the fact that $w_{y+X_1[0,t],\ldots,y+X_p[0,t]}(y+\cdot,u)$ does not depend on $y$. Thus, $t\mapsto t\underline{\Lambda}_p(t)$ is superadditive and so the limit in (4.9) indeed exists. It follows from Proposition 2.3 and (3.3) that $\lambda_p \leq p\nu\gamma\bar{w}(0)$, proving that $\lambda_p$ is finite strictly inside the weakly $p$-catalytic regime. □

4.2. *Convexity in $\kappa$.* We will write down a formal expansion of the expectation in the right-hand side of (2.9). From this expansion, it will immediately follow that $\Lambda_p(t)$ is a convex function of $\kappa$ for any $p$, $t$ and $\gamma$, $\rho$, $\nu$. After that, we can pass to the limit $t\to\infty$ to conclude that $\lambda_p = \lim_{t\to\infty}\Lambda_p(t)$ is also a convex function of $\kappa$.

PROPOSITION 4.2. *For any $p\in\mathbb{N}$,*

$$
\begin{aligned}
&\mathbb{E}_{0,\ldots,0}^{X_1,\ldots,X_p}\Bigg(\exp\Bigg[\nu\gamma\int_0^t \sum_{q=1}^{p} w(X_q(s),s)\,ds\Bigg]\Bigg) \\
&\quad = \sum_{n=0}^{\infty} (\nu\gamma)^n \Bigg(\prod_{m=1}^{n} \int_0^{s_{m-1}} ds_m\Bigg)\Bigg(\prod_{m=1}^{n}\sum_{r_m=1}^{p}\sum_{l_m=1}^{\infty}\Bigg)\gamma^{\sum_{m=1}^{n} l_m} \\
&\qquad\quad \times \Bigg(\prod_{\alpha=1}^{n}\prod_{\beta=1}^{l_\alpha}\int_0^{u_{\alpha,\beta-1}} du_{\alpha,\beta} \oint dk_{\alpha,\beta}\Bigg) \\
&\qquad\quad \times \exp\Bigg[-\rho\sum_{\alpha=1}^{n}\sum_{\beta=1}^{l_\alpha}(u_{\alpha,\beta-1}-u_{\alpha,\beta})\widehat{\varphi}(k_{\alpha,\beta})\Bigg]\Bigg(\prod_{\alpha=1}^{n}\prod_{\gamma=1}^{l_\alpha}\sum_{r_{\alpha,\gamma}=1}^{p}\Bigg)
\end{aligned} \tag{4.13}
$$

$$\times \exp\left[-\kappa \sum_{q=1}^{p} \int_0^t dv\, \widehat{\varphi}\left(\sum_{\alpha=1}^{n}\sum_{\beta=1}^{l_\alpha} k_{\alpha,\beta}\{\delta_{r_{\alpha,\beta},q}\mathbb{1}_{[0,u_{\alpha,\beta}]}(v) - \delta_{r_{\alpha,\beta-1},q}\mathbb{1}_{[0,u_{\alpha,\beta-1}]}(v)\}\right)\right],$$

*with the convention that* $s_0 = t$, $r_{\alpha,0} = r_\alpha$ *and* $u_{\alpha,0} = s_\alpha$, $\alpha \in \mathbb{N}$.

Proof. By Taylor expansion, we have

$$\mathbb{E}_{0,\dots,0}^{X_1,\dots,X_p}\left(\exp\left[\nu\gamma \int_0^t \sum_{q=1}^{p} w(X_q(s), s)\, ds\right]\right) = \sum_{n=0}^{\infty} (\nu\gamma)^n \left(\prod_{m=1}^{n} \int_0^{s_{m-1}} ds_m\right) \mathbb{E}_{0,\dots,0}^{X_1,\dots,X_p}\left(\prod_{m=1}^{n}\sum_{q=1}^{p} w(X_q(s_m), s_m)\right) \tag{4.14}$$

with $s_0 = t$. To compute the $n$-point correlation under the integral, we return to (2.16). By substituting (2.14) into (2.16) and iterating the resulting equation, we obtain the expansion

$$\begin{aligned} w(X_r(t), t) = \sum_{l=1}^{\infty} \gamma^l &\left(\prod_{\beta=1}^{l} \int_0^{u_{\beta-1}} du_\beta \oint dk_\beta\right) \\ &\times \exp\left[-\rho \sum_{\beta=1}^{l} (u_{\beta-1} - u_\beta)\widehat{\varphi}(k_\beta)\right] \\ &\times \left(\prod_{\gamma=1}^{l}\sum_{r_\gamma=1}^{p}\right) \exp\left\{i \sum_{\beta=1}^{l} k_\beta \cdot [X_{r_\beta}(u_\beta) - X_{r_{\beta-1}}(u_{\beta-1})]\right\} \end{aligned} \tag{4.15}$$

with $u_0 = t$ and $r_0 = r$. This expansion is convergent because the summand is bounded above by $(\gamma t p)^l / l!$. Using (4.15) in (4.14), we obtain

$$\begin{aligned} &\mathbb{E}_{0,\dots,0}^{X_1,\dots,X_p}\left(\prod_{m=1}^{n}\sum_{q=1}^{p} w(X_q(s_m), s_m)\right) \\ &\quad = \left(\prod_{m=1}^{n}\sum_{r_m=1}^{p}\sum_{l_m=1}^{\infty}\right) \gamma^{\sum_{m=1}^{n} l_m} \left(\prod_{\alpha=1}^{n}\prod_{\beta=1}^{l_\alpha} \int_0^{u_{\alpha,\beta-1}} du_{\alpha,\beta} \oint dk_{\alpha,\beta}\right) \\ &\qquad \times \exp\left[-\rho \sum_{\alpha=1}^{n}\sum_{\beta=1}^{l_\alpha} (u_{\alpha,\beta-1} - u_{\alpha,\beta})\widehat{\varphi}(k_{\alpha,\beta})\right] \left(\prod_{\alpha=1}^{n}\prod_{\gamma=1}^{l_\alpha}\sum_{r_{\alpha,\gamma}=1}^{p}\right) \\ &\qquad \times \mathbb{E}_{0,\dots,0}^{X_1,\dots,X_p}\left(\exp\left\{i \sum_{\alpha=1}^{n}\sum_{\beta=1}^{l_\alpha} k_{\alpha,\beta} \cdot [X_{r_{\alpha,\beta}}(u_{\alpha,\beta}) - X_{r_{\alpha,\beta-1}}(u_{\alpha,\beta-1})]\right\}\right) \end{aligned} \tag{4.16}$$

with $r_{\alpha,0}=r_\alpha$ and $u_{\alpha,0}=s_\alpha$, $\alpha=1,\ldots,n$. To complete the proof, it therefore suffices to show that

$$\mathbb{E}_{0,\ldots,0}^{X_1,\ldots,X_p}\left(\exp\left\{i\sum_{\alpha=1}^{n}\sum_{\beta=1}^{l_\alpha}k_{\alpha,\beta}\cdot[X_{r_{\alpha,\beta}}(u_{\alpha,\beta})-X_{r_{\alpha,\beta-1}}(u_{\alpha,\beta-1})]\right\}\right)$$

$$=\exp\left[-\kappa\sum_{q=1}^{p}\int_0^t dv\,\widehat{\varphi}\left(\sum_{\alpha=1}^{n}\sum_{\beta=1}^{l_\alpha}k_{\alpha,\beta}\{\delta_{r_{\alpha,\beta},q}\mathbb{1}_{[0,u_{\alpha,\beta}]}(v)\right.\right.$$
$$\left.\left.-\,\delta_{r_{\alpha,\beta-1},q}\mathbb{1}_{[0,u_{\alpha,\beta-1}]}(v)\}\right)\right].\tag{4.17}$$

By writing

$$\begin{aligned}X_{r_{\alpha,\beta}}&(u_{\alpha,\beta})-X_{r_{\alpha,\beta-1}}(u_{\alpha,\beta-1})\\&=\sum_{q=1}^{p}\{\delta_{r_{\alpha,\beta},q}X_q(u_{\alpha,\beta})-\delta_{r_{\alpha,\beta-1},q}X_q(u_{\alpha,\beta-1})\}\\&=\sum_{q=1}^{p}\int_0^t\{\delta_{r_{\alpha,\beta},q}\mathbb{1}_{[0,u_{\alpha,\beta}]}(v)-\delta_{r_{\alpha,\beta-1},q}\mathbb{1}_{[0,u_{\alpha,\beta-1}]}(v)\}\,dX_q(v)\end{aligned}\tag{4.18}$$

and noting that the increments $dX_q(v)$, $q=1,\ldots,p$, are independent, we see that (4.17) is a special case of the relation

$$\mathbb{E}_0^{X_q}\left(\exp\left[i\int_0^t f(v)\cdot dX_q(v)\right]\right)=\exp\left[-\kappa\int_0^t\widehat{\varphi}(f(v))\,dv\right],\qquad q=1,\ldots,p,\tag{4.19}$$

which holds for any $f:\mathbb{R}^d\to\mathbb{R}$ that is piecewise continuous and has bounded jumps. To see why (4.19) is true, we note that

$$\mathbb{E}_0^{X_q}(\exp[ik\cdot X_q(t)])=\sum_{x\in\mathbb{Z}^d}e^{ik\cdot x}p_\kappa(x,t)\tag{4.20}$$

with $p_\kappa$ denoting the transition kernel associated with $\kappa\Delta$. It follows from (2.14) that

$$\mathbb{E}_0^{X_q}(\exp[ik\cdot X_q(t)])=\exp[-\kappa t\widehat{\varphi}(k)].\tag{4.21}$$

From this relation, together with the fact that the increments of the process $X_q$ over disjoint time intervals are independent, we obtain (4.19). □

The expression in Proposition 4.2 is complicated, but the relevant point is that the right-hand side is a linear combination with nonnegative coefficients of functions that are negative exponentials in $\kappa$. Such a quantity is log-convex in $\kappa$, which tells us that $\Lambda_p(t)$ is convex in $\kappa$ [recall (2.9)]. Consequently, $\lambda_p=\lim_{t\to\infty}\Lambda_p(t)$ is also convex in $\kappa$.

4.3. *Small $\kappa$.* If $\kappa = 0$, then $X_1, \dots, X_p$ stay at the origin and so, from (2.9) and (2.12), we have that

$$\Lambda_p(t;0) = \nu\gamma \frac{1}{t} \int_0^t \bar{w}(0,s)\, ds. \tag{4.22}$$

Since $t \mapsto \bar{w}(0,t)$ is nondecreasing by Proposition 2.3, we have

$$\lambda_p(0) = \nu\gamma\bar{w}(0) \tag{4.23}$$

with $\bar{w}(0) = \lim_{t\to\infty} \bar{w}(0,t)$ given by (2.31). This proves the second equality in (1.18) in Theorem 1.5(ii). It follows from (3.3) and (4.22) that

$$\lambda_p(0) - 2d\kappa \le \lambda_p(\kappa) \le \lambda_p(0). \tag{4.24}$$

Hence, $\kappa \mapsto \lambda_p(\kappa)$ is continuous at 0 and bounded on $[0,\infty)$. This proves the first equality in (1.18) in Theorem 1.5(ii). Since $\kappa \mapsto \lambda_p(\kappa)$ is convex, as was shown in Section 4.2, it must be continuous and nonincreasing on $[0,\infty)$. Since it tends to zero like $1/\kappa$ as $\kappa \to \infty$ [as stated in Theorem 1.5(iii), which will be proven in Sections 5–8], it must be strictly positive and strictly decreasing on $[0,\infty)$. Thus, we have proven Theorem 1.5(i).

By Proposition 2.3 and (4.23), $\lambda_p(0) = \infty$ when $d \ge 3$, $p\gamma/\rho = r_d$. It therefore follows from (4.24) that $\lambda_p(\kappa) = \infty$. Thus, we have proven Theorem 1.4(iii). The proof of Theorem 1.4(ii) was already achieved with Lemma 4.1.

## 5. Proof of Theorem 1.5(iii).

The proof is long and technical. In Section 5.1, we introduce an appropriate scaling in $\kappa$. In Section 5.2, we formulate seven key lemmas that are the main ingredients of the proof. In Section 5.3, we prove Theorem 1.5(iii) subject to these lemmas. The proofs of the lemmas are deferred to Sections 6–8.

5.1. *Scaling.* To exhibit the dependence on the parameters, we henceforth write

$$\Lambda_p(T) = \Lambda_p(T;\kappa,\gamma,\rho,\nu), \tag{5.1}$$

where $\Lambda_p(T)$ is defined in (1.9). Substituting (2.16) into (2.9), we find that

$$\begin{aligned} &\Lambda_p(T;\kappa,\gamma,\rho,\nu) \\ &= \frac{1}{pT} \log \mathbb{E}_{0,\dots,0}^{X_1^\kappa,\dots,X_p^\kappa}\Bigg(\exp\Bigg[\nu\gamma^2 \sum_{k,l=1}^p \int_0^T ds \int_s^T dt \\ &\qquad \times p_\rho(X_l^\kappa(t) - X_k^\kappa(s), t-s)(1 + w(X_k^\kappa(s), s))\Bigg]\Bigg). \end{aligned} \tag{5.2}$$

In this formula, $X_1^\kappa, \dots, X_p^\kappa$ are independent simple random walks on $\mathbb{Z}^d$ with diffusion constant $\kappa$ (i.e., step rate $2d\kappa$), the expectation is over these random walks starting at 0, $p_\rho$ is the transition kernel associated with $\rho\Delta$ and $w$ denotes the solution of the Cauchy problem

$$\frac{\partial w}{\partial t} = \rho\Delta w + \gamma\left(\sum_{k=1}^{p} \delta_{X_k^\kappa(t)}\right)(1+w), \qquad w(\cdot,0) \equiv 0. \tag{5.3}$$

In Sections 2–4, the upper index $\kappa$ was suppressed. We introduce it here because we now want to remove the dependence of the random walks on $\kappa$. Indeed, in (5.2), we perform a time scaling $t \to t/\kappa$ in order to obtain

$$\Lambda_p(T;\kappa,\gamma,\rho,\nu) = \kappa\Lambda_p(\kappa T;1,\gamma/\kappa,\rho/\kappa,\nu). \tag{5.4}$$

Hence,

$$\Lambda_p(T;\kappa,\gamma,\rho,\nu) = \kappa\Lambda_p^*(\kappa T;\kappa,\gamma,\rho,\nu), \tag{5.5}$$

where

$$\begin{aligned} &\Lambda_p^*(T;\kappa,\gamma,\rho,\nu) \\ &= \frac{1}{pT}\log \mathbb{E}_{0,\dots,0}^{X_1,\dots,X_p}\left(\exp\left[\frac{\nu\gamma^2}{\kappa^2}\sum_{k,l=1}^{p}\int_0^T ds \int_s^T dt \right.\right. \\ &\qquad \left.\left. \times\, p_{\rho/\kappa}(X_l(t)-X_k(s),t-s)(1+w^*(X_k(s),s))\right]\right), \end{aligned} \tag{5.6}$$

$X_1,\dots,X_p$ are simple random walks on $\mathbb{Z}^d$ with diffusion constant 1 and $w^*$ solves

$$\frac{\partial w^*}{\partial t} = \frac{\rho}{\kappa}\Delta w^* + \frac{\gamma}{\kappa}\left(\sum_{k=1}^{p} \delta_{X_k(t)}\right)(1+w^*), \qquad w^*(\cdot,0) \equiv 0, \tag{5.7}$$

and satisfies $w^* \geq 0$.

The Lyapunov exponents in Theorem 1.4(ii)–(iii) are [recall (1.10)]

$$\lambda_p = \lambda_p(\kappa,\gamma,\rho,\nu) = \lim_{T\to\infty} \Lambda_p(T;\kappa,\gamma,\rho,\nu). \tag{5.8}$$

Because of (5.5), these are related to the rescaled Lyapunov exponents

$$\lambda_p^*(\kappa,\gamma,\rho,\nu) = \lim_{T\to\infty} \Lambda_p^*(T;\kappa,\gamma,\rho,\nu) \tag{5.9}$$

via

$$\lambda_p(\kappa,\gamma,\rho,\nu) = \kappa\lambda_p^*(\kappa,\gamma,\rho,\nu). \tag{5.10}$$

Also, note that (5.4) leads to the scaling

$$\lambda_p(\kappa,\gamma,\rho,\nu) = \kappa\lambda_p(1,\gamma/\kappa,\rho/\kappa,\nu). \tag{5.11}$$

We will frequently suppress the parameters $\gamma,\rho,\nu$ from the notation and write $\Lambda_p(T;\kappa)$, $\Lambda_p^*(T;\kappa)$ and $\lambda_p(\kappa)$, $\lambda_p^*(\kappa)$.

5.2. *Main ingredients of the proof.* The assertion of Theorem 1.5(iii) may now be restated as follows:

THEOREM 5.1. *Let $d \geq 3$, $p \in \mathbb{N}$ and*

$$0 < \frac{p\gamma}{\rho} < r_d. \tag{5.12}$$

(i) *For $d \geq 4$,*

$$\lim_{\kappa\to\infty} \kappa^2 \lambda_p^*(\kappa) = \frac{\nu\gamma^2}{r_d}. \tag{5.13}$$

(ii) *For $d = 3$,*

$$\lim_{\kappa\to\infty} \kappa^2 \lambda_p^*(\kappa) = \frac{\nu\gamma^2}{r_3} + \left(\frac{\nu\gamma^2}{\rho} p\right)^2 \mathcal{P} \tag{5.14}$$

*with $\mathcal{P}$ the constant defined in (1.16).*

The proof of Theorem 5.1 is based on seven lemmas, which are stated below and which provide lower and upper bounds for various parts contributing to (5.6). The guiding idea behind these lemmas is that the expectation in (5.6) can be moved to the exponential in the limit as $\kappa \to \infty$ uniformly in $T$, except for the part that produces the constant $\mathcal{P}$, which needs a large deviation analysis. This idea, though simple, is technically rather involved.

In the statement of the lemmas below, the following three auxiliary parameters appear:

$$0 < a < \infty, \qquad 0 < \varepsilon < K < \infty. \tag{5.15}$$

These parameters are needed to separate various time regimes. Four lemmas involve one random walk $(X)$, one lemma involves two random walks $(X, Y)$ and two lemmas involve $p$ random walks $(X_1, \dots, X_p)$. We use upper indices $-$ and $+$ for lim inf and lim sup, respectively.

5.2.1. *Lower bound.* The first lemma concerns the "diagonal term" ($0 \leq t - s \leq a\kappa^3$). Let

$$\begin{aligned} &\Lambda^-_{\mathrm{diag}}(T; a, \kappa) \\ &\quad = -\frac{1}{T} \log \mathbb{E}_0^X \left( \exp\left[ -\frac{\nu\gamma^2}{\kappa^2} \int_0^T ds \int_s^{s+a\kappa^3} dt\, p_{\rho/\kappa}(X(t) - X(s), t - s) \right] \right) \end{aligned} \tag{5.16}$$

and

$$\lambda^-_{\mathrm{diag}}(a, \kappa) = \liminf_{T\to\infty} \Lambda^-_{\mathrm{diag}}(T; a, \kappa). \tag{5.17}$$

LEMMA 5.2 (Lower bound for the diagonal term). *For* $d \geq 3$,

$$\liminf_{\kappa\to\infty} \kappa^2 \lambda^-_{\mathrm{diag}}(a,\kappa) \geq \frac{\nu\gamma^2}{r_d} \qquad \forall\, 0 < a < \infty. \tag{5.18}$$

The second lemma concerns the "variational term" ($\varepsilon\kappa^3 \leq t-s \leq K\kappa^3$), which involves $p$ random walks and which will turn out to be responsible for the second term in the right-hand side of (5.14). Let

$$\begin{aligned} &\Lambda_{\mathrm{var}}(T;\varepsilon,K,\kappa) \\ &\quad = \frac{1}{pT}\log \mathbb{E}^{X_1,\dots,X_p}_{0,\dots,0}\left(\exp\left[\frac{\nu\gamma^2}{\kappa^2}\sum_{k,l=1}^{p}\int_0^T ds \int_{s+\varepsilon\kappa^3}^{s+K\kappa^3} dt \right.\right. \\ &\qquad\qquad \left.\left. \times\, p_{\rho/\kappa}(X_l(t)-X_k(s),t-s)\right]\right) \end{aligned} \tag{5.19}$$

and

$$\lambda^-_{\mathrm{var}}(\varepsilon,K,\kappa) = \liminf_{T\to\infty} \Lambda_{\mathrm{var}}(T;\varepsilon,K,\kappa). \tag{5.20}$$

LEMMA 5.3 (Lower bound for the variational term). *For* $d = 3$,

$$\begin{aligned} &\liminf_{\kappa\to\infty} \kappa^2\lambda^-_{\mathrm{var}}(\varepsilon,K,\kappa) \geq \mathcal{P}_p(\varepsilon,K;\gamma,\rho,\nu) \\ &\qquad\qquad \forall\, 0<\varepsilon<K<\infty, \end{aligned} \tag{5.21}$$

*where*

$$\begin{aligned} &\mathcal{P}_p(\varepsilon,K;\gamma,\rho,\nu) \\ &\quad = \sup_{\substack{f\in H^1(\mathbb{R}^3)\\ \|f\|_2=1}} \left[\frac{\nu\gamma^2}{\rho} p \int_{\mathbb{R}^3} dx\, f^2(x) \int_{\mathbb{R}^3} dy\, f^2(y) \int_{\varepsilon\rho}^{K\rho} dt \right. \\ &\qquad\qquad \left. \times\, p_G(x-y,t) - \|\nabla_{\mathbb{R}^3} f\|_2^2\right] \end{aligned} \tag{5.22}$$

*with* $p_G(x,t) = (4\pi t)^{-3/2}\exp[-\|x\|^2/4t]$ *the Gaussian transition kernel associated with* $\Delta_{\mathbb{R}^3}$.

5.2.2. *Upper bound.* The third lemma is the counterpart of Lemma 5.2. Let

$$\begin{aligned} &\Lambda^+_{\mathrm{diag}}(T;a,\kappa) \\ &\quad = \frac{1}{T}\log \mathbb{E}^X_0\left(\exp\left[\frac{\nu\gamma^2}{\kappa^2}\int_0^T ds \int_s^{s+a\kappa^3} dt\, p_{\rho/\kappa}(X(t)-X(s),t-s)\right]\right) \end{aligned} \tag{5.23}$$

and

$$\lambda^+_{\mathrm{diag}}(a,\kappa) = \limsup_{T\to\infty} \Lambda^+_{\mathrm{diag}}(T;a,\kappa). \tag{5.24}$$

LEMMA 5.4 (Upper bound for the diagonal term).

(i) *For $d \geq 4$,*

$$\limsup_{\kappa\to\infty} \kappa^2 \lambda^+_{\mathrm{diag}}(a,\kappa) \leq \frac{\nu\gamma^2}{r_d} \qquad \forall\, 0 < a < \infty. \tag{5.25}$$

(ii) *For $d = 3$,*

$$\limsup_{a\downarrow 0} \limsup_{\kappa\to\infty} \kappa^2 \lambda^+_{\mathrm{diag}}(a,\kappa) \leq \frac{\nu\gamma^2}{r_3}. \tag{5.26}$$

The fourth lemma is the counterpart of Lemma 5.3. Let

$$\lambda^+_{\mathrm{var}}(\varepsilon,K,\kappa) = \limsup_{T\to\infty} \Lambda_{\mathrm{var}}(T;\varepsilon,K,\kappa). \tag{5.27}$$

LEMMA 5.5 (Upper bound for the variational term).

(i) *For $d \geq 4$,*

$$\lim_{\kappa\to\infty} \kappa^2 \lambda^+_{\mathrm{var}}(\varepsilon,K,\kappa) = 0 \qquad \forall\, 0 < \varepsilon < K < \infty. \tag{5.28}$$

(ii) *For $d = 3$,*

$$\limsup_{\kappa\to\infty} \kappa^2 \lambda^+_{\mathrm{var}}(\varepsilon,K,\kappa) \leq \mathcal{P}_p(\varepsilon,K;\gamma,\rho,\nu) \qquad \forall\, 0 < \varepsilon < K < \infty. \tag{5.29}$$

Three more lemmas deal with the upper bound, all of which turn out to involve terms that are negligible in the limit as $\kappa\to\infty$. The fifth lemma concerns the "off-diagonal" term ($t-s > a\kappa^3$). Let

$$\begin{aligned} &\Lambda_{\mathrm{off}}(T;a,\kappa) \\ &\quad = \frac{1}{T}\log \mathbb{E}_0^X\left(\exp\left[\frac{\nu\gamma^2}{\kappa^2}\int_0^T ds \int_{s+a\kappa^3}^\infty dt\, p_{\rho/\kappa}(X(t)-X(s),t-s)\right]\right) \end{aligned} \tag{5.30}$$

and

$$\lambda^+_{\mathrm{off}}(a,\kappa) = \limsup_{T\to\infty} \Lambda_{\mathrm{off}}(T;a,\kappa). \tag{5.31}$$

LEMMA 5.6 (Upper bound for the off-diagonal term).

(i) *For $d \geq 4$,*

$$\lim_{\kappa\to\infty} \kappa^2 \lambda^+_{\mathrm{off}}(a,\kappa) = 0 \qquad \forall\, 0 < a < \infty. \tag{5.32}$$

(ii) *For* $d = 3$,

$$\lim_{a\to\infty} \limsup_{\kappa\to\infty} \kappa^2 \lambda^+_{\mathrm{off}}(a,\kappa) = 0. \tag{5.33}$$

The sixth lemma concerns the "mixed" term and involves two random walks. Let

$$\begin{aligned}\Lambda_{\mathrm{mix}}(T;a,\kappa)\\ &= \frac{1}{T}\log \mathbb{E}^{X,Y}_{0,0}\left(\exp\left[\frac{\nu\gamma^2}{\kappa^2}\int_0^T ds \int_s^{s+a\kappa^3} dt\, p_{\rho/\kappa}(Y(t)-X(s),t-s)\right]\right)\end{aligned} \tag{5.34}$$

and

$$\lambda^+_{\mathrm{mix}}(a,\kappa) = \limsup_{T\to\infty} \Lambda_{\mathrm{mix}}(T;a,\kappa). \tag{5.35}$$

Lemma 5.7 (Upper bound for the mixed term).

(i) *For* $d \geq 4$,

$$\lim_{\kappa\to\infty} \kappa^2 \lambda^+_{\mathrm{mix}}(\infty,\kappa) = 0. \tag{5.36}$$

(ii) *For* $d = 3$,

$$\lim_{\kappa\to\infty} \kappa^2 \lambda^+_{\mathrm{mix}}(a,\kappa) = 0 \qquad \forall\, 0 < a < a_0 \textit{ with } a_0 \textit{ sufficiently small.} \tag{5.37}$$

The seventh lemma deals with a term that will be needed to treat the $w^*$-remainder in (5.6). Let

$$\begin{aligned}\Lambda_{\mathrm{rem}}(T;\kappa)\\ &= \frac{1}{T}\log \mathbb{E}^{X}_{0}\left(\exp\left[\frac{\nu\gamma^3}{\kappa^3}\int_0^T ds \left(\int_s^{\infty} dt\, p_{\rho/\kappa}(X(t)-X(s),t-s)\right)\right.\right.\\ &\qquad\qquad \left.\left.\times\left(\int_0^s du\, p_{\rho/\kappa}(X(s)-X(u),s-u)\right)\right]\right)\end{aligned} \tag{5.38}$$

and

$$\lambda^+_{\mathrm{rem}}(\kappa) = \limsup_{T\to\infty} \Lambda_{\mathrm{rem}}(T;\kappa). \tag{5.39}$$

[Note the extra factor $\gamma/\kappa$ in the exponent in the right-hand side of (5.38) compared to the previous definitions.]

Lemma 5.8 (Upper bound for the $w^*$-remainder). *For* $d \geq 3$,

$$\lim_{\kappa\to\infty} \kappa^2 \lambda^+_{\mathrm{rem}}(\kappa) = 0. \tag{5.40}$$

The proofs of Lemmas 5.2–5.8 are deferred to Sections 6–8.

5.3. *Proof of Theorem* 5.1. Recall that the solution of (5.7) admits the (implicit) integral representation [compare with (2.16)]

$$w^*(x,s) = \frac{\gamma}{\kappa}\sum_{l=1}^{p}\int_0^s du\, p_{\rho/\kappa}(x - X_l(u), s-u)(1 + w^*(X_l(u),u)). \tag{5.41}$$

Moreover, in the weakly catalytic regime given by (5.12), we have

$$w^*(x,s) \le \bar{w}(0) = C^* = \frac{p\gamma/\rho}{r_d - p\gamma/\rho} < \infty \qquad \forall x \in \mathbb{Z}^d,\ s \ge 0 \tag{5.42}$$

[recall (2.12), (2.16) and (2.31)]. Note that $C^*$ does not depend on $\kappa$.

For $d \ge 3$ and $a \ge 0$, abbreviate

$$G_a(0) = \int_a^\infty dt\, p(0,t). \tag{5.43}$$

We have $G_0(0) = R(0) = 1/r_d$ [recall (1.13)] and there exists a constant $c_d > 0$ such that

$$G_a(0) \le \frac{c_d}{r_d a^{(d-2)/2}}, \qquad a > 0. \tag{5.44}$$

5.3.1. *Lower bound.* Removing from (5.6) the terms with $w^*$, $t > s + K\kappa^3$ and $k \ne l$ for $t \le s + \varepsilon\kappa^3$, we get

$$\Lambda_p^*(T;\kappa,\gamma,\rho,\nu) \ge \frac{1}{pT}\log E_{0,\dots,0}^{X_1,\dots,X_p}(\exp[U + V - C]) \tag{5.45}$$

with

$$\begin{aligned} U &= \frac{\nu\gamma^2}{\kappa^2}\sum_{k=1}^{p}\int_0^T ds \int_s^{s+\varepsilon\kappa^3} dt\, p_{\rho/\kappa}(X_k(t) - X_k(s), t-s), \\ V &= \frac{\nu\gamma^2}{\kappa^2}\sum_{k,l=1}^{p}\int_0^T ds \int_{s+\varepsilon\kappa^3}^{s+K\kappa^3} dt\, p_{\rho/\kappa}(X_l(t) - X_k(s), t-s), \end{aligned} \tag{5.46}$$

where $C > 0$ is a constant that compensates for $t > T$ in (5.46). This constant may be chosen independently of $T$, as follows easily from rough estimates. By a reverse version of Hölder's inequality, we have

$$\begin{aligned} &E_{0,\dots,0}^{X_1,\dots,X_p}(\exp[U+V]) \\ &\quad \ge (E_{0,\dots,0}^{X_1,\dots,X_p}(\exp[-\zeta U]))^{-1/\zeta}(E_{0,\dots,0}^{X_1,\dots,X_p}(\exp[\theta V]))^{1/\theta}, \\ &\qquad\qquad \theta \in (0,1),\ \zeta = \frac{\theta}{1-\theta}. \end{aligned} \tag{5.47}$$

Hence, recalling (5.16) and (5.19), we obtain

$$\Lambda_p^*(T;\kappa,\gamma,\rho,\nu) \geq \frac{1}{\zeta}\Lambda_{\text{diag}}^-(T;\varepsilon,\kappa,\gamma,\rho,\zeta\nu) + \frac{1}{\theta}\Lambda_{\text{var}}(T;\varepsilon,K,\kappa,\gamma,\rho,\theta\nu). \tag{5.48}$$

By letting $T\to\infty$, recalling (5.9), letting $\kappa\to\infty$, using Lemmas 5.2 and 5.3 for the corresponding terms in the right-hand side and afterward letting $\theta\uparrow 1$, we arrive at

$$\liminf_{\kappa\to\infty}\kappa^2\lambda_p^*(\kappa) \geq \frac{\nu\gamma^2}{r_d}, \qquad \text{if } d\geq 4 \tag{5.49}$$

[drop the last term in (5.48)] and

$$\liminf_{\kappa\to\infty}\kappa^2\lambda_p^*(\kappa) \geq \frac{\nu\gamma^2}{r_3} + \mathcal{P}_p(\varepsilon,K;\gamma,\rho,\nu), \qquad \text{if } d=3 \tag{5.50}$$

[keep the last term in (5.48)]. In the latter, let $\varepsilon\downarrow 0$ and $K\to\infty$ and use the fact that, as is explained below,

$$\lim_{\varepsilon\downarrow 0, K\to\infty}\mathcal{P}_p(\varepsilon,K;\gamma,\rho,\nu) = \mathcal{P}_p(\gamma,\rho,\nu) \tag{5.51}$$

with

$$\mathcal{P}_p(\gamma,\rho,\nu) = \sup_{\substack{f\in H^1(\mathbb{R}^3)\\ \|f\|_2=1}}\left[\frac{\nu\gamma^2}{\rho}p\int_{\mathbb{R}^3}dx\,f^2(x)\int_{\mathbb{R}^3}dy\,f^2(y)\int_0^\infty dt \times p_G(x-y,t) - \|\nabla_{\mathbb{R}^3}f\|_2^2\right] \tag{5.52}$$

to obtain

$$\liminf_{\kappa\to\infty}\kappa^2\lambda_p^*(\kappa) \geq \frac{\nu\gamma^2}{r_3} + \mathcal{P}_p(\gamma,\rho,\nu), \qquad \text{if } d=3. \tag{5.53}$$

Finally, a straightforward scaling argument shows that

$$\mathcal{P}_p(\gamma,\rho,\nu) = \left(\frac{\nu\gamma^2}{\rho}p\right)^2\mathcal{P} \tag{5.54}$$

with $\mathcal{P}$ the constant defined in (1.16). This completes the proof of the lower bound in Theorem 5.1.

The fact that (5.51) holds is an immediate consequence of the fact that (1.16) and, hence, (5.52) has a maximizer $\bar{f}$, as shown by Lieb [25]. Indeed,

we have

$$
\begin{aligned}
0 &\le \mathcal{P}_p(\gamma,\rho,\nu) - \mathcal{P}_p(\varepsilon,K;\gamma,\rho,\nu) \\
&\le \frac{\nu\gamma^2}{\rho} p \int_{\mathbb{R}^3} dx\, \bar{f}^2(x) \int_{\mathbb{R}^3} dy\, \bar{f}^2(y) \\
&\quad \times \int_{(0,\varepsilon\rho)\cup(K\rho,\infty)} dt\, p_G(x-y,t)
\end{aligned}
\tag{5.55}
$$

and the right-hand side tends to zero as $\varepsilon \downarrow 0$ and $K \to \infty$ because the full integral is finite.

5.3.2. *Upper bound.* We begin by splitting the exponent in the right-hand side of (5.6) into various parts. The splitting is done with the various lemmas of Section 5.2.2 in mind and uses the parameters in (5.15) with $a = \varepsilon$ or $a = K$.

LEMMA 5.9. *For any $p \in \mathbb{N}$,*

$$
\begin{aligned}
&\sum_{k,l=1}^{p} \int_0^T ds \int_s^T dt\, p_{\rho/\kappa}(X_l(t) - X_k(s), t-s)(1 + w^*(X_k(s),s)) \\
&\quad \le \left(1 + \frac{D_\varepsilon}{\kappa^{d-2}}\right)(I + II + III) + \left(1 + \frac{D_\varepsilon}{\kappa^{d-2}} + 2(1+C^*)\frac{\gamma p}{r_d \rho}\right) IV \\
&\qquad + (1+C^*)\frac{\gamma}{\kappa} V,
\end{aligned}
\tag{5.56}
$$

*where $C^*$ is the constant in (5.42),*

$$
D_\varepsilon = \frac{(1+C^*) c_d \gamma p}{r_d \rho^{d/2} \varepsilon^{(d-2)/2}}
\tag{5.57}
$$

*with $c_d$ the constant in (5.44) and*

$$
\begin{aligned}
I &= \sum_{k=1}^{p} \int_0^T ds \int_s^{s+\varepsilon\kappa^3} dt\, p_{\rho/\kappa}(X_k(t) - X_k(s), t-s), \\
II &= \sum_{k,l=1}^{p} \int_0^T ds \int_{s+\varepsilon\kappa^3}^{s+K\kappa^3} dt\, p_{\rho/\kappa}(X_l(t) - X_k(s), t-s), \\
III &= \sum_{k,l=1}^{p} \int_0^T ds \int_{s+K\kappa^3}^{\infty} dt\, p_{\rho/\kappa}(X_l(t) - X_k(s), t-s), \\
IV &= \sum_{\substack{k,l=1 \\ k \ne l}}^{p} \int_0^T ds \int_s^{s+\varepsilon\kappa^3} dt\, p_{\rho/\kappa}(X_l(t) - X_k(s), t-s),
\end{aligned}
\tag{5.58}
$$

$$V = \sum_{k=1}^{p} \int_0^T ds \left( \int_0^s dr\, p_{\rho/\kappa}(X_k(s) - X_k(r), s - r) \right) \times \left( \int_s^\infty dt\, p_{\rho/\kappa}(X_k(t) - X_k(s), t - s) \right).$$

PROOF. For the term without $w^*$, we bound

$$\sum_{k,l=1}^{p} \int_0^T ds \int_s^T dt\, p_{\rho/\kappa}(X_l(t) - X_k(s), t - s) \le I + II + III + IV. \tag{5.59}$$

For the term with $w^*$, we bound, using (5.41) and (5.42),

$$\begin{aligned} &\sum_{k,l=1}^{p} \int_0^T ds \int_s^T dt\, p_{\rho/\kappa}(X_l(t) - X_k(s), t - s) w^*(X_k(s), s) \\ &\quad \le (1 + C^*) \frac{\gamma}{\kappa} \sum_{j,k,l=1}^{p} \int_0^T ds \left( \int_0^s dr\, p_{\rho/\kappa}(X_k(s) - X_j(r), s - r) \right) \\ &\qquad \times \left( \int_s^T dt\, p_{\rho/\kappa}(X_l(t) - X_k(s), t - s) \right). \end{aligned} \tag{5.60}$$

By (5.44),

$$\int_{\varepsilon\kappa^3}^\infty du\, p_{\rho/\kappa}(0, u) \le \frac{C_\varepsilon}{\kappa^{d-3}} \quad \text{with } C_\varepsilon = \frac{c_d}{r_d \rho^{d/2} \varepsilon^{(d-2)/2}}. \tag{5.61}$$

Hence, by (2.15),

$$\begin{aligned} &\int_0^{(s-\varepsilon\kappa^3)\vee 0} dr\, p_{\rho/\kappa}(X_k(s) - X_j(r), s - r) \\ &\qquad \le \int_{-\infty}^{s-\varepsilon\kappa^3} dr\, p_{\rho/\kappa}(0, s - r) \le \frac{C_\varepsilon}{\kappa^{d-3}}, \\ &\int_{(s+\varepsilon\kappa^3)\wedge T}^T dt\, p_{\rho/\kappa}(X_l(t) - X_k(s), t - s) \\ &\qquad \le \int_{s+\varepsilon\kappa^3}^\infty dt\, p_{\rho/\kappa}(0, t - s) \le \frac{C_\varepsilon}{\kappa^{d-3}}. \end{aligned} \tag{5.62}$$

Splitting the integrals in the two factors in the right-hand side of (5.60) into two parts, accordingly, and inserting (5.62), we find that

rhs (5.60)

$$
\begin{aligned}
&\le (1+C^*)\frac{\gamma}{\kappa}\sum_{j,k,l=1}^{p}\int_0^T ds\left(\int_{(s-\varepsilon\kappa^3)\vee 0}^{s} dr\, p_{\rho/\kappa}(X_k(s)-X_j(r),s-r)\right)\\
&\qquad\qquad \times\left(\int_s^{(s+\varepsilon\kappa^3)\wedge T} dt\, p_{\rho/\kappa}(X_l(t)-X_k(s),t-s)\right)\\
&\quad+\frac{D_\varepsilon}{\kappa^{d-2}}\sum_{k,l=1}^{p}\int_0^T ds\int_s^T dt\, p_{\rho/\kappa}(X_l(t)-X_k(s),t-s)
\end{aligned}
\tag{5.63}
$$

with $D_\varepsilon = 2(1+C^*)C_\varepsilon\gamma p$.

The second term in the right-hand side of (5.63) can be estimated using (5.59). For the first term, split the sum over the indices into $j=k=l$, $j\neq k$ and $k\neq l$. For $k\neq l$ ($j\neq k$), we estimate the first (second) inner integral by $\kappa/r_d\rho$. As a result, we obtain

$$
\begin{aligned}
\text{lhs (5.60)} &\le \frac{D_\varepsilon}{\kappa^{d-2}}(I+II+III+IV)\\
&\quad + 2(1+C^*)\frac{\gamma p}{r_d\rho}IV + (1+C^*)\frac{\gamma}{\kappa}V.
\end{aligned}
\tag{5.64}
$$

Combining (5.59) and (5.64), we arrive at the claimed assertion. □

Our next step is to apply Hölder's inequality to separate the various summands appearing in (5.58) so that we can apply to them the lemmas of Section 5.2.2. We will separate all summands except the ones in $II$, since the latter produces the variational problem in (5.22) and requires a cooperation of the $p$ random walks.

The total number of summands in (5.58) that are separated thus equals $q = p+1+p^2+p(p-1)+p = 2p^2+p+1$. Hence, substituting (5.58) into (5.56), substituting the resulting formula into (5.6) and applying Hölder's inequality

$$
\begin{aligned}
&E(e^{\sum_{r=1}^q S_r}) \le [E(e^{\theta S_1})]^{1/\theta}\prod_{r=2}^{q}[E(e^{\zeta S_r})]^{1/\zeta},\\
&\qquad\qquad \theta\in(1,\infty),\ \zeta=\frac{\theta}{\theta-1}(q-1),
\end{aligned}
\tag{5.65}
$$

to the expectation in the right-hand side of (5.6) (with $r=1$ reserved for $II$), we find that

$$
\begin{aligned}
&p\Lambda_p^*(T;\kappa,\gamma,\rho,\nu)\\
&\quad\le \frac{p}{\zeta}\Lambda_{\mathrm{diag}}^+\left(T;\varepsilon,\kappa,\gamma,\rho,\left(1+\frac{D_\varepsilon}{\kappa^{d-2}}\right)\zeta\nu\right)
\end{aligned}
$$

$$
(5.66)\qquad
\begin{aligned}
&+\frac{1}{\theta}p\Lambda_{\mathrm{var}}\left(T;\varepsilon,K,\kappa,\gamma,\rho,\left(1+\frac{D_\varepsilon}{\kappa^{d-2}}\right)\theta\nu\right)\\
&+\frac{p^2}{\zeta}\Lambda_{\mathrm{off}}\left(T;K,\kappa,\gamma,\rho,\left(1+\frac{D_\varepsilon}{\kappa^{d-2}}\right)\zeta\nu\right)\\
&+\frac{p(p-1)}{\zeta}\Lambda_{\mathrm{mix}}\left(T;\varepsilon,\kappa,\gamma,\rho,\left(1+\frac{D_\varepsilon}{\kappa^{d-2}}+2(1+C^*)\frac{\gamma p}{r_d\rho}\right)\zeta\nu\right)\\
&+\frac{p}{\zeta}\Lambda_{\mathrm{rem}}(T;\kappa,\gamma,\rho,(1+C^*)\zeta\nu).
\end{aligned}
$$

By letting $T\to\infty$, recalling (5.9), letting $\kappa\to\infty$, using Lemmas 5.4–5.8 for the corresponding terms in the right-hand side of (5.66) and afterward letting $\theta\downarrow 1$, we arrive at

$$
(5.67)\qquad \limsup_{\kappa\to\infty}\kappa^2\lambda_p^*(\kappa)\le\frac{\nu\gamma^2}{r_d},\qquad \text{if } d\ge 4,
$$

and after estimating $\mathcal{P}_p(\varepsilon,K;\gamma,\rho,\nu)\le\mathcal{P}_p(\gamma,\rho,\nu)$, using (5.26) with $a=\varepsilon$ and letting $\varepsilon\downarrow 0$, we arrive at

$$
(5.68)\qquad \limsup_{\kappa\to\infty}\kappa^2\lambda_p^*(\kappa)\le\frac{\nu\gamma^2}{r_3}+\mathcal{P}_p(\gamma,\rho,\nu),\qquad \text{if } d=3.
$$

For the second term in the right-hand side of (5.68), we may use (5.54). This completes the proof of the upper bound in Theorem 5.1.

**6. Proofs of Lemmas 5.2 and 5.4.** As we saw in Section 5.3, the "diagonal" contributions to the lower and the upper bound in the proof of Theorem 5.1 come from Lemmas 5.2 and 5.4, respectively. In this section, we prove these two lemmas. Let $p(x,t)$ denote the transition kernel associated with $\Delta$. Then $p_{\rho/\kappa}(x,t)=p(x,\frac{\rho}{\kappa}t)$.

6.1. *Proof of Lemma* 5.2.

PROOF OF LEMMA 5.2. Let $a,A>0$ be arbitrary. Estimate

$$
(6.1)\qquad
\begin{aligned}
&\int_0^T ds\int_s^{s+a\kappa^3}dt\,p_{\rho/\kappa}(X(t)-X(s),t-s)\\
&\quad\ge\left(\sum_{\substack{k=1\\ \text{even}}}^{\lfloor T/A\rfloor}+\sum_{\substack{k=1\\ \text{odd}}}^{\lfloor T/A\rfloor}\right)\int_{(k-1)A}^{kA}ds\int_s^{s+A}dt\,p_{\rho/\kappa}(X(t)-X(s),t-s)\\
&\qquad\qquad\qquad\forall\,\kappa\ge\kappa_0(a,A)=(A/a)^{1/3}.
\end{aligned}
$$

Note that the summands in each of the two sums are i.i.d. Hence, substituting (6.1) into (5.16) and applying the Cauchy–Schwarz inequality, we find that

$$\Lambda^-_{\mathrm{diag}}(T;a,\kappa) \geq -\frac{\lfloor T/A \rfloor}{2T} \log \mathbb{E}^X_0(\exp[-2W(A,\kappa)]) \tag{6.2}$$

with

$$W(A,\kappa) = \frac{\nu\gamma^2}{\kappa^2}\int_0^A ds \int_s^{s+A} dt\, p_{\rho/\kappa}(X(t)-X(s),t-s). \tag{6.3}$$

Next, note that by (2.15),

$$W(A,\kappa) \leq \frac{\nu\gamma^2}{\kappa^2}\int_0^A ds \int_s^{s+A} dt\, p_{\rho/\kappa}(0,t-s) \leq \frac{\nu\gamma^2}{\kappa^2} A \frac{\kappa}{\rho}\frac{1}{r_3}. \tag{6.4}$$

Since, for fixed $A$, the right-hand side tends to zero as $\kappa \to \infty$, it follows that

$$\begin{aligned} \mathbb{E}^X_0(\exp[-2W(A,\kappa)]) &\leq \exp[-2\theta\mathbb{E}^X_0(W(A,\kappa))] \\ &\qquad \forall\, \theta \in (0,1),\ \kappa \geq \kappa_1(\theta, A). \end{aligned} \tag{6.5}$$

Indeed, given $\theta \in (0,1)$, we can find an $\alpha(\theta) > 0$ such that $e^{-x} \leq 1 - \theta x$ for $0 \leq x \leq \alpha(\theta)$. Hence, for any random variable $\xi$ with $0 \leq \xi \leq \alpha(\theta)$, we have $\mathbb{E}(e^{-\xi}) \leq 1 - \theta\mathbb{E}(\xi) \leq e^{-\theta\mathbb{E}(\xi)}$.

Moreover, since

$$\begin{aligned} \mathbb{E}^X_0(p_{\rho/\kappa}(X(t)-X(s),t-s)) &= \mathbb{E}^X_0\left(p\left(X(t)-X(s), \frac{\rho}{\kappa}(t-s)\right)\right) \\ &= p\left(0, \left(1+\frac{\rho}{\kappa}\right)(t-s)\right), \end{aligned} \tag{6.6}$$

it follows from (6.3) that

$$\mathbb{E}^X_0(W(A,\kappa)) = \frac{\nu\gamma^2}{\kappa^2} A \int_0^A du\, p\left(0, \left(1+\frac{\rho}{\kappa}\right)u\right). \tag{6.7}$$

Inserting (6.5) and (6.7) into (6.2) and letting $T \to \infty$, we find that

$$\lambda^-_{\mathrm{diag}}(a,\kappa) \geq \theta \frac{\nu\gamma^2}{\kappa^2}\left(1+\frac{\rho}{\kappa}\right)^{-1} \int_0^{(1+\rho/\kappa)A} du\, p(0,u). \tag{6.8}$$

Hence,

$$\liminf_{\kappa\to\infty} \kappa^2 \lambda^-_{\mathrm{diag}}(a,\kappa) \geq \theta\nu\gamma^2 \int_0^A du\, p(0,u). \tag{6.9}$$

Now, let $A \to \infty$ and $\theta \uparrow 1$ to obtain the claimed assertion in (5.18). $\square$

6.2. *Proof of Lemma 5.4.* The proof of Lemma 5.4 relies on Lemma 6.1 below. For $a>0$, define

$$\Lambda_a(\gamma,\rho,\nu)=\limsup_{\kappa\to\infty}\frac{1}{a\kappa}\log\mathbb{E}_0^X\biggl(\exp\biggl[\frac{\nu\gamma^2}{\kappa^2}\int_0^{a\kappa^3}ds\int_s^\infty dt \times p_{\rho/\kappa}(X(t)-X(s),t-s)\biggr]\biggr). \tag{6.10}$$

Lemma 6.1. (a) *If $d\ge 4$, then*

$$\Lambda_a(\gamma,\rho,\nu)\le\frac{\nu\gamma^2}{r_d}\qquad\forall\,0<a<\infty. \tag{6.11}$$

(b) *If $d=3$, then*

$$\Lambda_a(\gamma,\rho,\nu)\le\frac{1+Ca^{1/4}(1/r_3)}{1-Ca^{1/4}(1+Ca^{1/4})(1/r_3)}\nu\gamma^2, \tag{6.12}$$

*provided $a>0$ is sufficiently small so that*

$$Ca^{1/4}(1+Ca^{1/4})\frac{1}{r_3}<1, \tag{6.13}$$

*where*

$$C=C(\gamma,\rho,\nu)=\biggl(\frac{2c_3\nu\gamma^2}{\sqrt{\rho}}\biggr)^{1/2}. \tag{6.14}$$

Before giving the proof of Lemma 6.1, we first prove Lemma 5.4.

Proof of lemma 5.4. Split the integral in the right-hand side of (5.23) as follows:

$$\begin{aligned}&\int_0^T ds\int_s^{s+a\kappa^3}dt\,p_{\rho/\kappa}(X(t)-X(s),t-s)\\&\quad\le\Biggl(\sum_{\substack{k=1\\ \text{even}}}^{\lceil T/a\kappa^3\rceil}+\sum_{\substack{k=1\\ \text{odd}}}^{\lceil T/a\kappa^3\rceil}\Biggr)\int_{(k-1)a\kappa^3}^{ka\kappa^3}ds\int_s^{s+a\kappa^3}dt\,p_{\rho/\kappa}(X(t)-X(s),t-s).\end{aligned} \tag{6.15}$$

Note that the summands in each of the two sums are i.i.d. Hence, substituting (6.15) into (5.23) and applying the Cauchy–Schwarz inequality, we find that

$$\begin{aligned}&\Lambda^+_{\text{diag}}(T;a,\kappa)\\&\quad\le\frac{\lceil T/a\kappa^3\rceil}{2T}\log\mathbb{E}_0^X\biggl(\exp\biggl[\frac{2\nu\gamma^2}{\kappa^2}\int_0^{a\kappa^3}ds\int_s^{s+a\kappa^3}dt\\&\qquad\times p_{\rho/\kappa}(X(t)-X(s),t-s)\biggr]\biggr).\end{aligned} \tag{6.16}$$

Letting $T \to \infty$, we arrive at

$$\lambda^+_{\mathrm{diag}}(a,\kappa) \le \frac{1}{2a\kappa^3} \log \mathbb{E}^X_0 \left( \exp\left[ \frac{2\nu\gamma^2}{\kappa^2} \int_0^{a\kappa^3} ds \int_s^{s+a\kappa^3} dt \right.\right.$$
$$\left.\left. \times\, p_{\rho/\kappa}(X(t)-X(s), t-s) \right] \right). \tag{6.17}$$

Assertion (5.25) follows from (6.17) after extending the second integral to infinity and applying Lemma 6.1(a) with $\nu$ replaced by $2\nu$. Assertion (5.26) follows similarly by applying Lemma 6.1(b). □

6.3. *Proof of Lemma* 6.1. The proof of Lemma 6.1 is based on two further lemmas. Recall (5.43).

Lemma 6.2. *For any $\alpha > 0$ and $M \in \mathbb{N}$,*

$$\mathbb{E}^X_0 \left( \exp\left[ \alpha \sum_{k=1}^M \int_0^\infty dt\, p_{\rho/\kappa}(Z_{k-1}(t) - Z_{k-1}(0), t) \right] \right)$$
$$\le \prod_{k=1}^M \max_{y_1,\ldots,y_{k-1}} \mathbb{E}^X_0 \left( \exp\left[ \alpha \sum_{l=0}^{k-1} \int_0^\infty dt\, p_{\rho/\kappa}\left( X(t) + y_l, \frac{l}{M} T + t \right) \right] \right), \tag{6.18}$$

*where $Z_k(t) = X(\frac{k}{M}T + t)$, $k \in \mathbb{N}_0$ and $y_0 = 0$.*

Lemma 6.3. *Let $d \ge 3$. For any $\alpha > 0$, $M \in \mathbb{N}$, $k \in \mathbb{N}_0$ and $y_0, \ldots, y_k \in \mathbb{Z}^d$,*

$$\mathbb{E}^X_0 \left( \exp\left[ \alpha \sum_{l=0}^{k} \int_0^\infty dt\, p_{\rho/\kappa}\left( X(t) + y_l, \frac{l}{M} T + t \right) \right] \right)$$
$$\le \exp\left[ \frac{\alpha \sum_{l=0}^k G_{(\rho T/\kappa M)l}(0)}{1 - \alpha \sum_{l=0}^k G_{(\rho T/\kappa M)l}(0)} \right], \tag{6.19}$$

*provided that $\alpha$ is sufficiently small so that*

$$\alpha \sum_{l=0}^k G_{(\rho T/\kappa M)l}(0) < 1. \tag{6.20}$$

Before giving the proof of Lemmas 6.2 and 6.3, we first prove Lemma 6.1.

Proof of Lemma 6.1. Let $M \in \mathbb{N}$ be arbitrary and abbreviate

$$Z_k(t) = X\left( \frac{k}{M} a\kappa^3 + t \right), \qquad k \in \mathbb{N}_0, \tag{6.21}$$

which is the same as that given below (6.18), with $T = a\kappa^3$. Then

$$\mathbb{E}_0^X\left(\exp\left[\frac{\nu\gamma^2}{\kappa^2}\int_0^{a\kappa^3} ds \int_s^\infty dt\, p_{\rho/\kappa}(X(t)-X(s),t-s)\right]\right)$$

(6.22)

$$= \mathbb{E}_0^X\left(\exp\left[\frac{\nu\gamma^2}{\kappa^2}\int_0^{a\kappa^3/M} ds \sum_{k=1}^{M}\int_s^\infty dt\, p_{\rho/\kappa}(Z_{k-1}(t)-Z_{k-1}(s),t-s)\right]\right).$$

After applying Jensen's inequality, we get

$$\text{rhs (6.22)} \le \frac{M}{a\kappa^3}\int_0^{a\kappa^3/M} ds$$

$$\times\, \mathbb{E}_0^X\left(\exp\left[\frac{\nu\gamma^2}{\kappa^2}\frac{a\kappa^3}{M}\sum_{k=1}^{M}\int_s^\infty dt\, p_{\rho/\kappa}(Z_{k-1}(t)-Z_{k-1}(s),t-s)\right]\right) \tag{6.23}$$

$$= \mathbb{E}_0^X\left(\exp\left[\nu\gamma^2\frac{a\kappa}{M}\sum_{k=1}^{M}\int_0^\infty dt\, p_{\rho/\kappa}(Z_{k-1}(t)-Z_{k-1}(0),t)\right]\right).$$

To the expression in the right-hand side, we may first apply Lemma 6.2 and then Lemma 6.3, both with $\alpha = \nu\gamma^2(a\kappa/M)$ and $T = a\kappa^3$. As a result, we obtain from (6.22) that

$$\frac{1}{a\kappa}\log \mathbb{E}_0^X\left(\exp\left[\frac{\nu\gamma^2}{\kappa^2}\int_0^{a\kappa^3} ds\int_s^\infty dt\, p_{\rho/\kappa}(X(t)-X(s),t-s)\right]\right)$$

$$\le \frac{1}{a\kappa}\sum_{k=1}^{M}\frac{\nu\gamma^2(a\kappa/M)\sum_{l=0}^{k-1}G_{\rho(a\kappa^2/M)l}(0)}{1-\nu\gamma^2(a\kappa/M)\sum_{l=0}^{k-1}G_{\rho(a\kappa^2/M)l}(0)} \tag{6.24}$$

$$\le \frac{\nu\gamma^2\sum_{l=0}^{M-1}G_{\rho(a\kappa^2/M)l}(0)}{1-\nu\gamma^2(a\kappa/M)\sum_{l=0}^{M-1}G_{\rho(a\kappa^2/M)l}(0)},$$

provided that

$$\nu\gamma^2\frac{a\kappa}{M}\sum_{l=0}^{M-1}G_{\rho(a\kappa^2/M)l}(0) < 1. \tag{6.25}$$

(a) Let $d \ge 4$. Then, by (5.44),

$$\sum_{l=0}^{M-1}G_{\rho(a\kappa^2/M)l}(0) \le G_0(0) + \left(\sum_{l=1}^{M-1}\frac{c_d}{\rho(a\kappa^2/M)l}\right)G_0(0)$$

(6.26)

$$\le \left(1+\widetilde{c}_d\frac{M\log M}{\rho a\kappa^2}\right)\frac{1}{r_d}$$

for some $\widetilde{c}_d > 0$ and all $M \in \mathbb{N}$. Now, choose

$$M = M(\kappa) = \lfloor \kappa^{3/2}\rfloor. \tag{6.27}$$

Then substituting (6.26) into (6.24) and letting $\kappa \to \infty$, we arrive at (6.11).

(b) Let $d=3$. Then, by (5.44),

$$\begin{aligned}\sum_{l=0}^{M-1} G_{\rho(a\kappa^2/M)l}(0) &\le G_0(0) + \left(\sum_{l=1}^{M-1} \frac{c_3}{\sqrt{\rho(a\kappa^2/M)l}}\right) G_0(0) \\ &\le \left(1 + \frac{2c_3}{\sqrt{\rho a}} \frac{M}{\kappa}\right) \frac{1}{r_3}\end{aligned} \tag{6.28}$$

for all $M \in \mathbb{N}$. Now, choose

$$M = M(\kappa) = \left\lfloor \left(\frac{\nu\gamma^2\sqrt{\rho}}{2c_3}\right)^{1/2} a^{3/4}\kappa \right\rfloor. \tag{6.29}$$

Then

$$\lim_{\kappa\to\infty} \frac{2c_3}{\sqrt{\rho a}} \frac{M(\kappa)}{\kappa} = Ca^{1/4} \tag{6.30}$$

and

$$\lim_{\kappa\to\infty} \nu\gamma^2 \frac{a\kappa}{M(\kappa)} = Ca^{1/4}, \tag{6.31}$$

where $C$ is given by (6.14). Substituting (6.28) into (6.24) and assumption (6.25), letting $\kappa \to \infty$, and taking into account (6.31) and (6.30), we arrive at (6.12) under assumption (6.13). $\square$

6.4. *Proofs of Lemmas 6.2 and 6.3.*

Proof of Lemma 6.2. We show that the function defined by

$$\begin{aligned}E(r) = \prod_{k=1}^{r} \max_{y_1,\dots,y_{k-1}} \mathbb{E}_0^X &\left(\exp\left[\alpha \sum_{l=0}^{k-1} \int_0^\infty dt\, p_{\rho/\kappa}\left(X(t)+y_l, \frac{l}{M}T+t\right)\right]\right) \\ \times \max_{z_1,\dots,z_r} \mathbb{E}_0^X &\left(\exp\left[\alpha \sum_{k=1}^{M-r} \int_0^\infty dt\, p_{\rho/\kappa}(Z_{k-1}(t)-Z_{k-1}(0),t)\right.\right. \\ &\left.\left. + \alpha \sum_{l=1}^{r} \int_0^\infty dt\, p_{\rho/\kappa}\left(X(t)+z_l, \frac{l}{M}T+t\right)\right]\right)\end{aligned} \tag{6.32}$$

for $r=0,\dots,M-1$ is nondecreasing in $r$. Then $E(0) \le E(M-1)$, which is the desired inequality. [Note that for $r=M-1$, the first term in the right-hand side of (6.32) corresponds to $k=1,\dots,M-1$ in the right-hand side of (6.18), the second term to $k=M$, $l=0$, and the third term to $k=M$, $l=1,\dots,M-1$.]

We now fix $r$ arbitrarily. We want to show that $E(r) \le E(r+1)$. To this end, we also fix $z_1,\dots,z_r$ arbitrarily. Separately handling the summand for $k=1$, splitting the integral over $(0,\infty)$ into integrals over $(0,T/M)$ and $(T/M,\infty)$, shifting time by $T/M$ for the latter and using the Markov property of $X$ at time $T/M$, we obtain

$$\begin{aligned}
&\mathbb{E}_0^X\Bigg(\exp\Bigg[\alpha\sum_{k=1}^{M-r}\int_0^\infty dt\, p_{\rho/\kappa}(Z_{k-1}(t)-Z_{k-1}(0),t)\\
&\qquad\qquad + \alpha\sum_{l=1}^{r}\int_0^\infty dt\, p_{\rho/\kappa}\Big(X(t)+z_l,\frac{l}{M}T+t\Big)\Bigg]\Bigg)\\
&\quad=\mathbb{E}_0^X\Bigg(\exp\Bigg[\alpha\int_0^{T/M} dt\, p_{\rho/\kappa}(X(t),t)\\
&\qquad\qquad + \alpha\int_0^\infty dt\, p_{\rho/\kappa}\Big(Z_1(t),\frac{1}{M}T+t\Big)\\
&\qquad\qquad + \alpha\sum_{k=1}^{M-(r+1)}\int_0^\infty dt\, p_{\rho/\kappa}(Z_k(t)-Z_k(0),t)\\
&\qquad\qquad + \alpha\sum_{l=1}^{r}\int_0^{T/M} dt\, p_{\rho/\kappa}\Big(X(t)+z_l,\frac{l}{M}T+t\Big)\\
&\qquad\qquad + \alpha\sum_{l=2}^{r+1}\int_0^\infty dt\, p_{\rho/\kappa}\Big(Z_1(t)+z_{l-1},\frac{l}{M}T+t\Big)\Bigg]\Bigg)
\end{aligned}\tag{6.33}$$

and

$$\begin{aligned}
&\text{rhs}(6.33)\\
&\quad\le\mathbb{E}_0^X\Bigg(\exp\Bigg[\alpha\int_0^{T/M} dt\, p_{\rho/\kappa}(X(t),t)\\
&\qquad\qquad + \alpha\sum_{l=1}^{r}\int_0^{T/M} dt\, p_{\rho/\kappa}\Big(X(t)+z_l,\frac{l}{M}T+t\Big)\Bigg]\Bigg)\\
&\qquad\times\max_{z_0}\mathbb{E}_0^X\Bigg(\exp\Bigg[\alpha\int_0^\infty dt\, p_{\rho/\kappa}\Big(X(t)+z_0,\frac{1}{M}T+t\Big)\\
&\qquad\qquad + \alpha\sum_{k=1}^{M-(r+1)}\int_0^\infty dt\, p_{\rho/\kappa}(Z_{k-1}(t)-Z_{k-1}(0),t)\\
&\qquad\qquad + \alpha\sum_{l=2}^{r+1}\int_0^\infty dt\, p_{\rho/\kappa}\Big(X(t)+z_0+z_{l-1},\frac{l}{M}T+t\Big)\Bigg]\Bigg).
\end{aligned}\tag{6.34}$$

In the last line, we have maximized over $Z_1(0) = X(T/M)$ after using the Markov property of $X$ at time $T/M$. Hence, combining (6.33) and (6.34), we get

$$\begin{aligned}
&\max_{z_1,\dots,z_r} \mathbb{E}_0^X\Biggl(\exp\Biggl[\alpha\sum_{k=1}^{M-r}\int_0^\infty dt\, p_{\rho/\kappa}(Z_{k-1}(t)-Z_{k-1}(0),t) \\
&\qquad\qquad + \alpha\sum_{l=1}^{r}\int_0^\infty dt\, p_{\rho/\kappa}\biggl(X(t)+z_l,\frac{l}{M}T+t\biggr)\Biggr]\Biggr) \\
&\quad\le \max_{y_1,\dots,y_r} \mathbb{E}_0^X\Biggl(\exp\Biggl[\alpha\sum_{l=0}^{r}\int_0^\infty dt\, p_{\rho/\kappa}\biggl(X(t)+y_l,\frac{l}{M}T+t\biggr)\Biggr]\Biggr) \\
&\qquad \times \max_{z_1,\dots,z_{r+1}} \mathbb{E}_0^X\Biggl(\exp\Biggl[\alpha\sum_{k=1}^{M-(r+1)}\int_0^\infty dt\, p_{\rho/\kappa}(Z_{k-1}(t)-Z_{k-1}(0),t) \\
&\qquad\qquad + \alpha\sum_{l=1}^{r+1}\int_0^\infty dt \\
&\qquad\qquad\qquad \times p_{\rho/\kappa}\biggl(X(t)+z_l,\frac{l}{M}T+t\biggr)\Biggr]\Biggr).
\end{aligned}\tag{6.35}$$

Here, we extend the first two integrals in the right-hand side of (6.34) from $T/M$ to infinity, use the fact that $y_0 = 0$ and replace $z_0$ by $z_1$ and $z_0 + z_{l-1}$ by $z_l$. Substituting (6.35) into (6.32), we get that $E(r) \le E(r+1)$, as desired. □

Proof of Lemma 6.3. A Taylor expansion of the exponential function yields

$$\begin{aligned}
&\mathbb{E}_0^X\Biggl(\exp\Biggl[\alpha\int_0^\infty dt\sum_{l=0}^{k} p_{\rho/\kappa}\biggl(X(t)+y_l,\frac{l}{M}T+t\biggr)\Biggr]\Biggr) \\
&\quad = \sum_{m=0}^{\infty}\alpha^m \mathbb{E}_0^X\Biggl(\prod_{j=1}^{m}\int_{t_{j-1}}^\infty dt_j\sum_{l=0}^{k} p_{\rho/\kappa}\biggl(X(t_j)+y_l,\frac{l}{M}T+t_j\biggr)\Biggr)
\end{aligned}\tag{6.36}$$

with $t_0 = 0$. A successive application of the Markov property at times $t_{m-1},\dots,t_1$ yields

$$\begin{aligned}
&\mathbb{E}_0^X\Biggl(\prod_{j=1}^{m}\int_{t_{j-1}}^\infty dt_j\sum_{l=0}^{k} p_{\rho/\kappa}\biggl(X(t_j)+y_l,\frac{l}{M}T+t_j\biggr)\Biggr) \\
&\quad = \mathbb{E}_0^X\Biggl(\prod_{j=1}^{m-1}\int_{t_{j-1}}^\infty dt_j\sum_{l=0}^{k} p_{\rho/\kappa}\biggl(X(t_j)+y_l,\frac{l}{M}T+t_j\biggr)\Biggr)
\end{aligned}$$

$$
(6.37)\quad \times \int_{t_{m-1}}^{\infty} dt_m \sum_{l=0}^{k} p\left(X(t_{m-1}) + y_l, \frac{\rho}{\kappa}\left(\frac{l}{M}T + t_m\right) + t_m - t_{m-1}\right)
$$

$$
\leq \mathbb{E}_0^X\left(\prod_{j=1}^{m-1}\int_{t_{j-1}}^{\infty} dt_j \sum_{l=0}^{k} p_{\rho/\kappa}\left(X(t_j) + y_l, \frac{l}{M}T + t_j\right)\right)\left(\sum_{l=0}^{k} G_{(\rho T/\kappa M)l}(0)\right)
$$

$$
\leq \cdots \leq \left(\sum_{l=0}^{k} G_{(\rho T/\kappa M)l}(0)\right)^m .
$$

In the first inequality, we have used the fact that $p(x,t) \leq p(0,t)$ and that $t \mapsto p(0,t)$ is nonincreasing. Substituting (6.37) into (6.36), summing the geometric series and using the inequality $1 + x \leq e^x$, $x \in \mathbb{R}$, we arrive at (6.19). $\square$

**7. Proofs of Lemmas 5.3 and 5.5.** As we saw in Section 5.3, the "variational" contributions to the lower and the upper bound in the proof of Theorem 5.1 come from Lemmas 5.3 and 5.5, respectively. In this section, we prove these two lemmas.

The proof of Lemma 5.5(i), which applies to $d \geq 4$, is easy. Indeed, in the right-hand side of (5.19), separate the $p^2$ summands with the help of Hölder's inequality [as in (5.66)]. The terms with $k = l$ are negligible for $\kappa \to \infty$, by Lemma 5.6(i) with $a = \varepsilon$, while the same is true for the terms with $k \neq l$, by Lemma 5.7(i). Lemmas 5.6 and 5.7 are proved in Section 8.

Thus, we may henceforth restrict our attention to $d = 3$.

7.1. *Space-time scaling.* We begin with a space-time scaling of the random walks. Let $\mathbb{Z}_\kappa^3 = \kappa^{-1}\mathbb{Z}^3$ and define

$$
(7.1)\quad \begin{aligned} X_k^{(\kappa)}(t) &= \kappa^{-1} X_k(\kappa^2 t), \qquad t \geq 0,\ k = 1, \ldots, p, \\ p^{(\kappa)}(x,t) &= \kappa^3 p(\kappa x, \kappa^2 t), \qquad x \in \mathbb{Z}_\kappa^3,\ t \geq 0. \end{aligned}
$$

Each $X_k^{(\kappa)}$ lives on $\mathbb{Z}_\kappa^3$, has generator

$$
(7.2)\quad (\Delta^{(\kappa)} f)(x) = \kappa^2 \sum_{\substack{y \in \mathbb{Z}_\kappa^3 \\ \|y-x\| = \kappa^{-1}}} [f(y) - f(x)], \qquad x \in \mathbb{Z}_\kappa^3,
$$

and has transition kernel whose density is $p^{(\kappa)}$ w.r.t. the discrete Lebesgue measure on $\mathbb{Z}_\kappa^3$, where each site carries weight $\kappa^{-3}$. As $\kappa \to \infty$, each $X_k^{(\kappa)}$ converges weakly to Brownian motion, which has as generator the continuous Laplacian $\Delta_{\mathbb{R}^3}$, and $p^{(\kappa)}$ converges weakly to $p_G$, the density of the transition kernel associated with Brownian motion w.r.t. the continuous Lebesgue measure on $\mathbb{R}^3$. The last convergence is uniform on compact sets, that is, for every compact set $C \subset \mathbb{R}^3 \times (0,\infty)$ and every $\theta \in (0,1)$, there exists $\kappa_0 = \kappa_0(C,\theta)$ such that

$$\theta p_G(x,t) \le p^{(\kappa)}(x,t) \le \frac{1}{\theta} p_G(x,t) \qquad \forall (x,t) \in C,\ \kappa \ge \kappa_0. \tag{7.3}$$

Further, note that

$$\min_{x\in\mathbb{R}^3} \frac{p_G(x,u_2)}{p_G(x,u_1)} = \frac{p_G(0,u_2)}{p_G(0,u_1)} = \left(\frac{u_1}{u_2}\right)^{3/2} \qquad \forall u_2 \ge u_1. \tag{7.4}$$

7.2. *Proof of Lemma 5.3.*

Proof of Lemma 5.3. Fix $0 < \varepsilon < K < \infty$, $\delta > 0$ small and $\theta \in (0,1)$. Abbreviate

$$\begin{aligned} L &= L(\delta,\varepsilon) = \lceil \varepsilon/\delta \rceil, \\ M &= M(\delta,K) = \lfloor K/\delta \rfloor, \\ N &= N(T;\delta,\kappa) = \lfloor T/\delta\kappa^3 \rfloor. \end{aligned} \tag{7.5}$$

Fix a large open cube $Q \subset \mathbb{R}^3$, centered at the origin. Later, we will take limits in the following order:

$$T \to \infty, \qquad \kappa \to \infty, \qquad \delta \downarrow 0, \qquad \theta \uparrow 1, \qquad Q \uparrow \mathbb{R}^3. \tag{7.6}$$

Let $C_Q$ be the event

$$\begin{aligned} C_Q &= C_Q(N,M,\delta,\kappa) \\ &= \{X_k^{(\kappa)}(t) \in Q\ \forall 0 \le t \le (N+M)\delta\kappa,\ k = 1,\ldots,p\}. \end{aligned} \tag{7.7}$$

Then from (5.19), (7.1) and the lower bound in (7.3), we get

$$\Lambda_{\mathrm{var}}(T;\varepsilon,K,\kappa) \ge \frac{1}{pT} \log \mathbb{E}_{0,\ldots,0}^{X_1^{(\kappa)},\ldots,X_p^{(\kappa)}} (\exp[U] \mathbb{1}_{C_Q}) \tag{7.8}$$

with

$$\begin{aligned} U &= \frac{\nu\gamma^2}{\kappa} \sum_{k,l=1}^{p} \int_0^{T/\kappa^2} ds \int_{s+\varepsilon\kappa}^{s+K\kappa} dt\, p^{(\kappa)}\Big(X_l^{(\kappa)}(t) - X_k^{(\kappa)}(s), \frac{\rho}{\kappa}(t-s)\Big) \\ &\ge \frac{\nu\gamma^2}{\kappa} \sum_{k,l=1}^{p} \int_0^{T/\kappa^2} ds \int_{s+\varepsilon\kappa}^{s+K\kappa} dt\, \theta p_G\Big(X_l^{(\kappa)}(t) - X_k^{(\kappa)}(s), \frac{\rho}{\kappa}(t-s)\Big) \end{aligned} \tag{7.9}$$

for $\kappa \ge \kappa_0(C,\theta)$ with $C = 2\bar{Q} \times [\varepsilon\rho, K\rho]$ ($\bar{Q}$ being the closure of $Q$). Moreover,

rhs (7.9)

$$\begin{aligned}
&\geq \frac{\nu\gamma^2}{\kappa}\sum_{k,l=1}^{p}\sum_{n=1}^{N}\int_{(n-1)\delta\kappa}^{n\delta\kappa} ds \int_{n\delta\kappa+\varepsilon\kappa}^{(n-1)\delta\kappa+K\kappa} dt \\
&\qquad \times \theta p_G\Big(X_l^{(\kappa)}(t) - X_k^{(\kappa)}(s), \frac{\rho}{\kappa}(t-s)\Big) \\
&\geq \frac{\nu\gamma^2}{\kappa}\sum_{k,l=1}^{p}\sum_{n=1}^{N}\sum_{m=L+1}^{M-1}\int_{(n-1)\delta\kappa}^{n\delta\kappa} ds \int_{(n+m-1)\delta\kappa}^{(n+m)\delta\kappa} dt \\
&\qquad \times \theta p_G\Big(X_l^{(\kappa)}(t) - X_k^{(\kappa)}(s), \frac{\rho}{\kappa}(t-s)\Big).
\end{aligned} \tag{7.10}$$

Next, note that $(m-1)\delta\rho \leq \frac{\rho}{\kappa}(t-s) \leq (m+1)\delta\rho$ for all $s,t$ in the domain of integration corresponding to $n,m$ and use (7.4) to obtain

$$\Lambda_{\text{var}}(T;\varepsilon,K,\kappa) \geq \frac{1}{pT}\log \mathbb{E}_{0,\dots,0}^{X_1^{(\kappa)},\dots,X_p^{(\kappa)}}(\exp[V]\mathbb{1}_{C_Q}) \tag{7.11}$$

with

$$\begin{aligned}
V = \frac{\nu\gamma^2}{\kappa}\sum_{k,l=1}^{p}\sum_{n=1}^{N}\sum_{m=L+1}^{M-1}\int_{(n-1)\delta\kappa}^{n\delta\kappa} ds \int_{(n+m-1)\delta\kappa}^{(n+m)\delta\kappa} dt\, \theta\Big(\frac{L}{L+2}\Big)^{3/2} \\
\times p_G(X_l^{(\kappa)}(t) - X_k^{(\kappa)}(s), (m-1)\delta\rho).
\end{aligned} \tag{7.12}$$

In this last expression, the time coordinate of the kernel is fixed for each $m$. Therefore, if we introduce the normalized occupation time measures

$$\begin{aligned}
\Xi_{k,r}^{(\kappa)}(A) = \frac{1}{\delta\kappa}\int_{(r-1)\delta\kappa}^{r\delta\kappa} ds\, \mathbb{1}_A(X_k^{(\kappa)}(s)), \\
k = 1,\dots,p,\ r = 1,\dots,N+M,\ A \subset \mathbb{R}^3 \text{ Borel},
\end{aligned} \tag{7.13}$$

then we may write

$$\begin{aligned}
V = \theta\Big(\frac{L}{L+2}\Big)^{3/2}\frac{\nu\gamma^2}{\rho}\delta\kappa \sum_{k,l=1}^{p}\sum_{n=1}^{N}\sum_{m=L+1}^{M-1}\int_Q \Xi_{k,n}^{(\kappa)}(dx)\int_Q \Xi_{l,n+m}^{(\kappa)}(dy) \\
\times \delta\rho p_G(y-x,(m-1)\delta\rho).
\end{aligned} \tag{7.14}$$

This representation puts us in a position where we can carry out a large deviation analysis, as follows.

For $\mu \in \mathcal{M}_1(Q)$, the set of probability measures on $Q$, let $\mathcal{U}_Q(\mu) \subset \mathcal{M}_1(Q)$ denote any weak open neighborhood of $\mu$ such that

$$\begin{aligned}\nu_1, \nu_2 \in \mathcal{U}_Q(\mu) \quad \Longrightarrow \quad & \int_Q \nu_1(dx) \int_Q \nu_2(dy) p_G(y-x,(m-1)\delta\rho) \\ & \geq \theta \int_Q \mu(dx) \int_Q \mu(dy) p_G(y-x,(m-1)\delta\rho) \\ & \qquad \forall\, m = L, \dots, M,\end{aligned} \tag{7.15}$$

and let $C_{Q,\mu}$ denote the event

$$C_{Q,\mu} = \{\Xi_{k,r}^{(\kappa)} \in \mathcal{U}_Q(\mu)\ \forall\, k = 1, \dots, p,\ r = 1, \dots, N+M\}. \tag{7.16}$$

Then, for any $\mu \in \mathcal{M}_1(Q)$, we may bound, via (7.11) and (7.14),

$$\begin{aligned}&\Lambda_{\mathrm{var}}(T;\varepsilon,K,\kappa) \\ &\quad \geq \frac{1}{pT} \log \mathbb{E}_{0,\dots,0}^{X_1^{(\kappa)},\dots,X_p^{(\kappa)}}(\exp[V] \mathbb{1}_{C_Q} \mathbb{1}_{C_{Q,\mu}}) \\ &\quad \geq \frac{1}{pT} \theta^2 \left(\frac{L}{L+2}\right)^{3/2} \frac{\nu\gamma^2}{\rho} p^2 N \\ &\qquad \times \delta\kappa \int_Q \mu(dx) \int_Q \mu(dy) \sum_{m=L+1}^{M-1} \delta\rho p_G(y-x,(m-1)\delta\rho) \\ &\qquad + \frac{1}{pT} \log \mathbb{P}_{0,\dots,0}^{X_1^{(\kappa)},\dots,X_p^{(\kappa)}}(C_Q \cap C_{Q,\mu}).\end{aligned} \tag{7.17}$$

By again appealing to (7.4), the sum in the first term in the right-hand side of (7.17) can be estimated as follows:

$$\begin{aligned}&\sum_{m=L+1}^{M-1} \delta\rho p_G(y-x,(m-1)\delta\rho) \\ &\qquad \geq \left(\frac{L}{L+2}\right)^{3/2} \int_{(L-1)\delta\rho}^{(M-2)\delta\rho} du\, p_G(y-x,u).\end{aligned} \tag{7.18}$$

As for the second term in the right-hand side of (7.17), by using the independence of the $p$ random walks as well as the Markov property at times $r\delta\kappa$ for $r = 1, \dots, N+M$, we may estimate (with $X^{(\kappa)} = X_1^{(\kappa)}$, $\Xi_r^{(\kappa)} = \Xi_{1,r}^{(\kappa)}$)

$$\begin{aligned}&\mathbb{P}_{0,\dots,0}^{X_1^{(\kappa)},\dots,X_p^{(\kappa)}}(C_Q \cap C_{Q,\mu}) \\ &\quad = [\mathbb{P}_0^{X^{(\kappa)}}(X^{(\kappa)}(t) \in Q\ \forall\, 0 \leq t \leq (N+M)\delta\kappa, \\ &\qquad\qquad \Xi_r^{(\kappa)} \in \mathcal{U}_Q(\mu)\ \forall\, r = 1, \dots, N+M)]^p\end{aligned} \tag{7.19}$$

$$\geq [\mathbb{P}_0^{X^{(\kappa)}}(X^{(\kappa)}(t)\in Q\ \forall\, 0\leq t\leq (N+M)\delta\kappa,$$
$$\Xi_r^{(\kappa)}\in\mathcal{U}_Q(\mu) \text{ and } X^{(\kappa)}(r\delta\kappa)\in\tfrac{1}{2}Q\ \forall\, r=1,\dots,N+M)]^p$$
$$\geq \Big[\min_{x\in\mathbb{Z}_\kappa^3\cap(1/2)Q}\mathbb{P}_x^{X^{(\kappa)}}(X^{(\kappa)}(t)\in Q\ \forall\, 0\leq t\leq\delta\kappa,$$
$$\Xi^{(\kappa)}\in\mathcal{U}_Q(\mu) \text{ and } X^{(\kappa)}(\delta\kappa)\in\tfrac{1}{2}Q)\Big]^{p(N+M)}.$$

The dependence on $N$ has now been pulled out of both terms in the right-hand side of (7.17) and so we can take the limit $T\to\infty$ to obtain from (5.20), (7.5) and (7.17)–(7.19) that

$$\begin{aligned}
&\kappa^2\lambda_{\mathrm{var}}^-(\varepsilon,K,\kappa)\\
&\quad\geq \theta^2\left(\frac{L}{L+2}\right)^3\frac{\nu\gamma^2}{\rho}p\int_Q\mu(dx)\int_Q\mu(dy)\int_{(L-1)\delta\rho}^{(M-2)\delta\rho}du\,p_G(y-x,u)\\
&\qquad+\frac{1}{\delta\kappa}\log\min_{x\in\mathbb{Z}_\kappa^3\cap(1/2)Q}\mathbb{P}_x^{X^{(\kappa)}}\Big(X^{(\kappa)}(t)\in Q\ \forall\, 0\leq t\leq\delta\kappa,\\
&\qquad\qquad\Xi^{(\kappa)}\in\mathcal{U}_Q(\mu) \text{ and } X^{(\kappa)}(\delta\kappa)\in\frac{1}{2}Q\Big)
\end{aligned}\tag{7.20}$$

for $\kappa\geq\kappa_0(C,\theta)$. The final step in the argument is the following large deviation bound:

Lemma 7.1. *For each $\mu\in\mathcal{M}_1(Q)$,*

$$\begin{aligned}
&\liminf_{\kappa\to\infty}\frac{1}{\delta\kappa}\log\min_{x\in\mathbb{Z}_\kappa^3\cap(1/2)Q}\mathbb{P}_x^{X^{(\kappa)}}\Big(X^{(\kappa)}(t)\in Q\ \forall\, 0\leq t\leq\delta\kappa,\\
&\qquad\Xi^{(\kappa)}\in\mathcal{U}_Q(\mu) \text{ and } X^{(\kappa)}(\delta\kappa)\in\frac{1}{2}Q\Big)\geq -S_Q(\mu)
\end{aligned}\tag{7.21}$$

*with $S_Q:\mathcal{M}_1(Q)\to[0,\infty]$ given by*

$$S_Q(\mu)=\begin{cases}\|\nabla_{\mathbb{R}^3}f\|_2^2, & \text{if } \mu\ll dx \text{ and } \sqrt{\dfrac{d\mu}{dx}}=f(x) \text{ with } f\in H_0^1(Q),\\ \infty, & \text{otherwise,}\end{cases}\tag{7.22}$$

*where $H_0^1(Q)$ is the completion of $C_c^\infty(Q)$ (the space of $C^\infty$-functions $f:Q\to\mathbb{R}$ with compact support) w.r.t. the $H^1$-norm $\|f\|_{H^1}=\|f\|_2+\|\nabla f\|_2$.*

The proof of Lemma 7.1 is deferred to Section 7.4. Letting $\kappa\to\infty$ in (7.20), using (7.21), letting $\delta\downarrow 0$, recalling (7.5), letting $\theta\uparrow 1$ and afterward taking the supremum over $\mu \in \mathcal{M}_1(Q)$, we arrive at

$$\begin{aligned} &\liminf_{\kappa\to\infty} \kappa^2 \lambda_{\mathrm{var}}^-(\varepsilon, K, \kappa) \\ &\geq \sup_{\substack{f\in H_0^1(Q)\\ \|f\|_2=1}} \biggl[\frac{\nu\gamma^2}{\rho} p \int_Q dx\, f^2(x) \int_Q dy\, f^2(y) \int_{\varepsilon\rho}^{K\rho} du\, p_G(y-x,u) \\ &\qquad\qquad - \|\nabla_{\mathbb{R}^3} f\|_2^2\biggr]. \end{aligned} \tag{7.23}$$

Finally, let $Q \uparrow \mathbb{R}^3$ and use a standard approximation argument to show that the variational expression in the right-hand side of (7.23) converges to

$$\begin{aligned} \sup_{\substack{f\in H^1(\mathbb{R}^3)\\ \|f\|_2=1}} &\biggl[\frac{\nu\gamma^2}{\rho} p \int_{\mathbb{R}^3} dx\, f^2(x) \int_{\mathbb{R}^3} dy\, f^2(y) \int_{\varepsilon\rho}^{K\rho} du\, p_G(y-x,u) \\ &\qquad - \|\nabla_{\mathbb{R}^3} f\|_2^2\biggr]. \end{aligned} \tag{7.24}$$

The latter is precisely $\mathcal{P}_p(\varepsilon, K; \gamma, \rho, \nu)$ as defined in (5.22), so we have completed the proof of Lemma 5.3. $\square$

7.3. *Proof of Lemma* 5.5. At the beginning of Section 7, we dealt with Lemma 5.5(i). Thus, we need only prove Lemma 5.5(ii).

Proof of Lemma 5.5(ii). Part of the argument runs parallel to Section 7.2. Fix $\varepsilon, K, \delta, \theta$ as before. Retain (7.5), but with $\lceil\cdot\rceil$ and $\lfloor\cdot\rfloor$ interchanged. Let $Q \subset \mathbb{R}^3$ be a large closed cube, centered at the origin. Later, we will again take limits in the order given in (7.6).

Let $l(Q)$ [resp. $l(Q^{(\kappa)})$] denote the side length of $Q$ [resp. $Q^{(\kappa)} = Q \cap \mathbb{Z}_\kappa^3$]. Let

$$\begin{aligned} &X_k^{(\kappa,Q)}(t), \qquad t \geq 0,\ k = 1, \ldots, p, \\ &p^{(\kappa,Q)}(x,t), \qquad x \in Q,\ t \geq 0, \end{aligned} \tag{7.25}$$

denote the $Q$-periodization of (7.1), that is,

$$\begin{aligned} X_k^{(\kappa,Q)}(t) &= X_k^{(\kappa)}(t) \ \mathrm{mod}(Q^{(\kappa)}), \\ p^{(\kappa,Q)}(x,t) &= \sum_{k\in\mathbb{Z}^3} p^{(\kappa)}\biggl(x + \frac{k}{\kappa} l(Q^{(\kappa)}), t\biggr). \end{aligned} \tag{7.26}$$

Similarly, let

$$p_G^{(Q)}(x,t), \qquad x \in Q,\ t \geq 0, \tag{7.27}$$

denote the $Q$-periodization of the Gaussian kernel, that is,

$$p_G^{(Q)}(x,t) = \sum_{k\in\mathbb{Z}^3} p_G(x+kl(Q),t). \tag{7.28}$$

From (5.19), (7.26) and the upper bound in (7.3) (which carries over to the $Q$-periodized kernels), we get

$$\Lambda_{\text{var}}(T;\varepsilon,K,\kappa) \le \frac{1}{pT}\log \mathbb{E}_{0,\dots,0}^{X_1^{(\kappa,Q)},\dots,X_p^{(\kappa,Q)}}(\exp[U]) \tag{7.29}$$

with

$$\begin{aligned} U &= \frac{\nu\gamma^2}{\kappa}\sum_{k,l=1}^{p}\int_0^{T/\kappa^2} ds \int_{s+\varepsilon\kappa}^{s+K\kappa} dt\, p^{(\kappa,Q)}\left(X_l^{(\kappa,Q)}(t) - X_k^{(\kappa,Q)}(s), \frac{\rho}{\kappa}(t-s)\right) \\ &\le \frac{\nu\gamma^2}{\kappa}\sum_{k,l=1}^{p}\int_0^{T/\kappa^2} ds \int_{s+\varepsilon\kappa}^{s+K\kappa} dt\, \frac{1}{\theta} p_G^{(Q)}\left(X_l^{(\kappa,Q)}(t) - X_k^{(\kappa,Q)}(s), \frac{\rho}{\kappa}(t-s)\right) \end{aligned} \tag{7.30}$$

for $\kappa \ge \kappa_0 = \kappa_0(C,\theta)$ with $C = 2Q\times[\varepsilon\rho, K\rho]$. Moreover,

$$\begin{aligned} &\text{rhs (7.30)} \\ &\quad \le \frac{\nu\gamma^2}{\kappa}\sum_{k,l=1}^{p}\sum_{n=1}^{N}\int_{(n-1)\delta\kappa}^{n\delta\kappa} ds \int_{(n-1)\delta\kappa+\varepsilon\kappa}^{n\delta\kappa+K\kappa} dt \\ &\qquad\qquad \times \frac{1}{\theta} p_G^{(Q)}\left(X_l^{(\kappa,Q)}(t) - X_k^{(\kappa,Q)}(s), \frac{\rho}{\kappa}(t-s)\right) \\ &\quad \le \frac{\nu\gamma^2}{\kappa}\sum_{k,l=1}^{p}\sum_{n=1}^{N}\sum_{m=L}^{M}\int_{(n-1)\delta\kappa}^{n\delta\kappa} ds \int_{(n+m-1)\delta\kappa}^{(n+m)\delta\kappa} dt \\ &\qquad\qquad \times \frac{1}{\theta} p_G^{(Q)}\left(X_l^{(\kappa,Q)}(t) - X_k^{(\kappa,Q)}(s), \frac{\rho}{\kappa}(t-s)\right). \end{aligned} \tag{7.31}$$

This is the analogue of (7.9) and (7.10).

Next, use (7.4) to obtain

$$\Lambda_{\text{var}}(T;\varepsilon,K,\kappa) \le \frac{1}{pT}\log \mathbb{E}_{0,\dots,0}^{X_1^{(\kappa,Q)},\dots,X_p^{(\kappa,Q)}}(\exp[V]), \tag{7.32}$$

with

$$V = \frac{\nu\gamma^2}{\kappa}\sum_{k,l=1}^{p}\sum_{n=1}^{N}\sum_{m=L}^{M}\int_{(n-1)\delta\kappa}^{n\delta\kappa} ds \int_{(n+m-1)\delta\kappa}^{(n+m)\delta\kappa} dt \frac{1}{\theta}\left(\frac{L+1}{L-1}\right)^{3/2}$$

$$\times p_G^{(Q)}(X_l^{(\kappa,Q)}(t) - X_k^{(\kappa,Q)}(s), (m+1)\delta\rho) \tag{7.33}$$

$$= \frac{1}{\theta}\left(\frac{L+1}{L-1}\right)^{3/2} \frac{\nu\gamma^2}{\rho} \delta\kappa \sum_{k,l=1}^{p} \sum_{n=1}^{N} \sum_{m=L}^{M} \int_Q \Xi_{k,n}^{(\kappa,Q)}(dx) \int_Q \Xi_{l,n+m}^{(\kappa,Q)}(dy)$$

$$\times \delta\rho p_G^{(Q)}(y-x,(m+1)\delta\rho),$$

which is the analogue of (7.12) and (7.14). Here,

$$\Xi_{k,r}^{(\kappa,Q)}(A) = \frac{1}{\delta\kappa} \int_{(r-1)\delta\kappa}^{r\delta\kappa} ds\, \mathbb{1}_A(X_k^{(\kappa,Q)}(s)),$$
$$k = 1,\dots,p, r = 1,\dots,N+M+1, A \subset Q \text{ Borel}, \tag{7.34}$$

is the analogue of (7.13).

For $\mu \in \mathcal{M}_1(Q)$, let $\mathcal{U}_Q(\mu) \subset \mathcal{M}_1(Q)$ be any weak neighborhood of $\mu$ such that

(1) for $\mu_1, \mu_2 \in \mathcal{M}_1(Q)$;

$$\nu_1 \in \mathcal{U}_Q(\mu_1),\ \nu_2 \in \mathcal{U}_Q(\mu_2) \quad \Longrightarrow \quad \int_Q \nu_1(dx) \int_Q \nu_2(dy) p_G^{(Q)}(y-x,u)$$
$$\le \frac{1}{\theta} \int_Q \mu_1(dx) \int_Q \mu_2(dy) p_G^{(Q)}(y-x,u) \tag{7.35}$$
$$\forall\, u \in [\varepsilon\rho, K\rho + 2\delta\rho];$$

(2) for $\mu \in \mathcal{M}_1(Q)$;

$$\inf_{\mu' \in \mathcal{U}_Q(\mu)} \widehat{S}_Q(\mu') \ge \theta \widehat{S}_Q(\mu). \tag{7.36}$$

Here, (7.35) is the analogue of (7.15), while $\widehat{S}_Q$ is the rate function defined in (7.45) below. The latter inequality can be achieved because $\mu \mapsto \widehat{S}_Q(\mu)$ is lower semi-continuous. Conditions (1) and (2) will be needed in the proof of Lemma 7.2 below (see Section 7.4).

Since $\mathcal{M}_1(Q)$ is compact, there exist finitely many $\mu_1, \dots, \mu_I \in \mathcal{M}_1(Q)$ (with $I$ not depending on $T, \kappa$) such that

$$\mathcal{M}_1(Q) \subset \bigcup_{i=1}^{I} \mathcal{U}_Q(\mu_i). \tag{7.37}$$

Let

$$\mathcal{J} = \{J : \{1,\dots,p\} \times \{1,\dots,N+M+1\} \to \{1,\dots,I\}\}. \tag{7.38}$$

For $J \in \mathcal{J}$, let $C_{Q,J}$ denote the event

$$C_{Q,J} = \{\Xi_{k,r}^{(\kappa,Q)} \in \mathcal{U}_Q(\mu_{J(k,r)})\ \forall\, k = 1,\dots,p,\ r = 1,\dots,N+M+1\}. \tag{7.39}$$

Then, because of (7.37), we may bound

$$\begin{aligned} &\Lambda_{\text{var}}(T;\varepsilon,K,\kappa) \\ &\quad \le \frac{1}{pT}\log\max_{J\in\mathcal{J}} \mathbb{E}_{0,\dots,0}^{X_1^{(\kappa,Q)},\dots,X_p^{(\kappa,Q)}}(\exp[V]\mathbb{1}_{C_{Q,J}}) + \frac{1}{pT}\log|\mathcal{J}|. \end{aligned} \tag{7.40}$$

On $C_{Q,J}$, we have, via (7.33) and (7.35),

$$\begin{aligned} V \le \frac{1}{\theta^2}\left(\frac{L+1}{L-1}\right)^{3/2}\frac{\nu\gamma^2}{\rho}\delta\kappa \sum_{k,l=1}^{p}\sum_{n=1}^{N}\sum_{m=L}^{M}\int_Q \mu_{J(k,n)}(dx)\int_Q \mu_{J(l,n+m)}(dy) \\ \times\, \delta\rho\, p_G^{(Q)}(y-x,(m+1)\delta\rho). \end{aligned} \tag{7.41}$$

Moreover, similarly as in (7.19),

$$\begin{aligned} &\mathbb{P}_{0,\dots,0}^{X_1^{(\kappa,Q)},\dots,X_p^{(\kappa,Q)}}(C_{Q,J}) \\ &\quad \le \prod_{k=1}^{p}\prod_{r=1}^{N+M+1}\max_{x\in\mathbb{Z}_\kappa^3\cap Q}\mathbb{P}_x^{X^{(\kappa,Q)}}(\Xi^{(\kappa,Q)}\in\mathcal{U}_Q(\mu_{J(k,r)})). \end{aligned} \tag{7.42}$$

Combining (7.40)–(7.42), it follows that

$$\begin{aligned} &\Lambda_{\text{var}}(T;\varepsilon,K,\kappa) \\ &\quad \le \frac{1}{pT}\max_{J\in\mathcal{J}}\Bigg[\frac{1}{\theta^2}\left(\frac{L+1}{L-1}\right)^{3/2}\frac{\nu\gamma^2}{\rho} \\ &\qquad\qquad \times \sum_{k,l=1}^{p}\sum_{n=1}^{N}\delta\kappa\sum_{m=L}^{M}\int_Q \mu_{J(k,n)}(dx)\int_Q \mu_{J(l,n+m)}(dy) \\ &\qquad\qquad\qquad \times\, \delta\rho p_G^{(Q)}(y-x,(m+1)\delta\rho) \\ &\qquad\qquad + \sum_{k=1}^{p}\sum_{r=1}^{N+M+1}\log\max_{x\in\mathbb{Z}_\kappa^3\cap Q}\mathbb{P}_x^{X^{(\kappa,Q)}}(\Xi^{(\kappa,Q)}\in\mathcal{U}_Q(\mu_{J(k,r)}))\Bigg] \\ &\qquad + \frac{1}{pT}\log|\mathcal{J}| \end{aligned} \tag{7.43}$$

for $\kappa \ge \kappa_0(C,\theta)$.

Below, we will need the following upper large deviation bound (with $\Xi^{(\kappa,Q)} = \Xi_{1,1}^{(\kappa,Q)}$) which is the reverse of Lemma 7.1:

LEMMA 7.2. *For each* $i \in \{1,\dots,I\}$,

$$\limsup_{\kappa\to\infty} \frac{1}{\delta\kappa} \log \max_{x\in\mathbb{Z}_\kappa^3\cap Q} \mathbb{P}_x^{X^{(\kappa,Q)}}(\Xi^{(\kappa,Q)} \in \mathcal{U}_Q(\mu_i)) \le -\theta \widehat{S}_Q(\mu_i) \tag{7.44}$$

*with* $\widehat{S}_Q$ *the* $Q$*-periodization of* $S_Q$*, that is,* $\widehat{S}_Q : \mathcal{M}_1(Q) \to [0,\infty]$ *is given by*

$$\widehat{S}_Q(\mu) = \begin{cases} \|\nabla_{\mathbb{R}^3} f\|_2^2, & \text{if } \mu \ll dx \\ & \text{and } \sqrt{\dfrac{d\mu}{dx}} = f(x) \text{ with } f \in H^1_{\mathrm{per}}(Q), \\ \infty, & \text{otherwise,} \end{cases} \tag{7.45}$$

*where* $H^1_{\mathrm{per}}(Q)$ *is the space of functions in* $H^1(Q)$ *with periodic boundary conditions.*

The proof of Lemma 7.2 is deferred to Section 7.4.

Next, define

$$\mu^J_{k,s} = \mu_{J(k,r)} \qquad \text{for } k = 1,\dots,p, \quad r = 1,\dots,N+M+1, (r-1)\delta\kappa \le s < r\delta\kappa. \tag{7.46}$$

The measure-valued paths $s \mapsto \mu^J_{k,s}$ are piecewise constant and take values in $\{\mu_1,\dots,\mu_I\}$. Once again using (7.4), we may revert back time from sums to integrals to obtain

$$\begin{aligned}
&\sum_{n=1}^N \delta\kappa \sum_{m=L}^M \int_Q \mu_{J(k,n)}(dx) \int_Q \mu_{J(l,n+m)}(dy) \delta\rho p_G^{(Q)}(y-x,(m+1)\delta\rho) \\
&\quad\le \frac{\rho}{\kappa}\left(\frac{L+3}{L+1}\right)^{3/2} \\
&\qquad \times \sum_{n=1}^N \int_{(n-1)\delta\kappa}^{n\delta\kappa} ds \sum_{m=L}^M \int_{(n+m-1)\delta\kappa}^{(n+m)\delta\kappa} dt \int_Q \mu^J_{k,s}(dx) \int_Q \mu^J_{l,t}(dy) \\
&\qquad \times p_G^{(Q)}\left(y-x, \frac{\rho}{\kappa}(t-s) + 2\delta\rho\right) \\
&\quad\le \frac{\rho}{\kappa}\left(\frac{L+1}{L-1}\right)^{3/2} \int_0^{N\delta\kappa} ds \int_{s+(L-1)\delta\kappa}^{s+(M+1)\delta\kappa} dt \\
&\qquad \times \int_Q \mu^J_{k,s}(dx) \int_Q \mu^J_{l,t}(dy) p_G^{(Q)}\left(y-x, \frac{\rho}{\kappa}(t-s) + 2\delta\rho\right) \\
&\quad\le \frac{\rho}{\kappa}\left(\frac{L+1}{L-1}\right)^{3/2} \int_0^{(N+M+1)\delta\kappa} ds \int_0^{(N+M+1)\delta\kappa} dt\, \mathbb{1}_{\{(L-1)\delta\kappa \le t-s \le (M+1)\delta\kappa\}} \\
&\qquad \times \int_Q \mu^J_{k,s}(dx) \int_Q \mu^J_{l,t}(dy) p_G^{(Q)}\left(y-x, \frac{\rho}{\kappa}(t-s) + 2\delta\rho\right)
\end{aligned} \tag{7.47}$$

and, according to Lemma 7.2,

$$\begin{aligned}
&\sum_{r=1}^{N+M+1} \log \max_{x\in\mathbb{Z}_\kappa^3\cap Q} \mathbb{P}_x^{X^{(\kappa,Q)}}(\Xi^{(\kappa,Q)} \in \mathcal{U}_Q(\mu_{J(k,r)})) \\
&\qquad = \int_0^{(N+M+1)\delta\kappa} ds \frac{1}{\delta\kappa} \log \max_{x\in\mathbb{Z}_\kappa^3\cap Q} \mathbb{P}_x^{X^{(\kappa,Q)}}(\Xi^{(\kappa,Q)} \in \mathcal{U}_Q(\mu_{k,s}^J)) \qquad (7.48)\\
&\qquad \le -\theta^2 \int_0^{(N+M+1)\delta\kappa} ds \widehat{S}_Q(\mu_{k,s}^J)
\end{aligned}$$

for $\kappa \ge \kappa_1(C,\theta) \ge \kappa_0(C,\theta)$. Inserting (7.47) and (7.48) into (7.43), we arrive at

$$\begin{aligned}
&\Lambda_{\mathrm{var}}(T;\varepsilon,K,\kappa) \\
&\quad \le \frac{1}{pT} \max_{J\in\mathcal{J}} \Bigg[ \frac{1}{\theta^2} \left(\frac{L+1}{L-1}\right)^3 \frac{\nu\gamma^2}{\kappa} \\
&\qquad\qquad \times \sum_{k,l=1}^{p} \int_0^{(N+M+1)\delta\kappa} ds \\
&\qquad\qquad\qquad \times \int_0^{(N+M+1)\delta\kappa} dt\, \mathbb{1}_{\{(L-1)\delta\kappa \le t-s \le (M+1)\delta\kappa\}} \\
&\qquad\qquad\qquad \times \int_Q \mu_{k,s}^J(dx) \qquad (7.49)\\
&\qquad\qquad\qquad \times \int_Q \mu_{l,t}^J(dy) p_G^{(Q)}\left(y-x, \frac{\rho}{\kappa}(t-s) + 2\delta\rho\right) \\
&\qquad\qquad\qquad\qquad - \theta^2 \sum_{k=1}^{p} \int_0^{(N+M+1)\delta\kappa} ds \widehat{S}_Q(\mu_{k,s}^J) \Bigg] \\
&\qquad + \frac{1}{pT} \log|\mathcal{J}|
\end{aligned}$$

for $\kappa \ge \kappa_1(C,\theta)$.

At this point we can perform a time-diagonalization.

LEMMA 7.3. *For every $A > 0$ and $\mu_{k,s} \in \mathcal{M}_1(Q)$ with $k = 1,\dots,p$, $0 \le s \le (N+M+1)\delta\kappa$,*

$$\begin{aligned}
&\frac{A}{\kappa} \int_0^{(N+M+1)\delta\kappa} ds \int_s^{(N+M+1)\delta\kappa} dt\, \mathbb{1}_{\{(L-1)\delta\kappa \le t-s \le (M+1)\delta\kappa\}} \\
&\qquad \times \sum_{k,l=1}^{p} \int_Q \mu_{k,s}(dx) \int_Q \mu_{l,t}(dy) p_G^{(Q)}\left(y-x, \frac{\rho}{\kappa}(t-s) + 2\delta\rho\right) \qquad (7.50)
\end{aligned}$$

$$
\begin{aligned}
&-\sum_{k=1}^{p}\int_0^{(N+M+1)\delta\kappa} ds\,\widehat{S}_Q(\mu_{k,s})\\
&\le p(N+M+1)\delta\kappa\\
&\quad\times \sup_{\nu\in\mathcal{M}_1(Q)}\biggl[\frac{A}{\kappa}p\int_Q\nu(dx)\int_Q\nu(dy)\int_{(L-1)\delta\kappa}^{(M+1)\delta\kappa}du\\
&\qquad\times p_G^{(Q)}\biggl(y-x,\frac{\rho}{\kappa}u+2\delta\rho\biggr)-\widehat{S}_Q(\nu)\biggr].
\end{aligned}
$$

The proof of Lemma 7.3 is given below. Inserting (7.50) with $A=\theta^{-4}(\frac{L+1}{L-1})^3\nu\gamma^2$ into (7.49), inserting (7.45), letting $T\to\infty$ and recalling (5.27), we obtain

$$
\begin{aligned}
&\kappa^2\lambda^+_{\mathrm{var}}(\varepsilon,K,\kappa)\\
&\quad\le\theta^2\sup_{\substack{f\in H^1_{\mathrm{per}}(Q)\\ \|f\|_2=1}}\biggl[\frac{1}{\theta^4}\biggl(\frac{L+1}{L-1}\biggr)^3\frac{\nu\gamma^2}{\rho}p\int_Q dx\,f^2(x)\int_Q dy\,f^2(y)\\
&\qquad\qquad\times\int_{(L-1)\delta\rho}^{(M+1)\delta\rho}du\,p_G^{(Q)}(y-x,u+2\delta\rho)-\|\nabla_{\mathbb{R}^3}f\|_2^2\biggr]\\
&\qquad+\frac{1}{\delta\kappa}\log I,
\end{aligned}
\tag{7.51}
$$

where we note that $\log|\mathcal{J}|=p(N+M)\log I$ and recall the last line of (7.5). Now, let $\kappa\to\infty$, $\delta\downarrow 0$ [yielding $L\to\infty$, $(L-1)\delta\to\varepsilon$ and $(M+1)\delta\to K$] and $\theta\uparrow 1$, to obtain

$$
\begin{aligned}
&\limsup_{\kappa\to\infty}\kappa^2\lambda^+_{\mathrm{var}}(\varepsilon,K,\kappa)\\
&\quad\le\sup_{\substack{f\in H^1_{\mathrm{per}}(Q)\\ \|f\|_2=1}}\biggl[\frac{\nu\gamma^2}{\rho}p\int_Q dx\,f^2(x)\int_Q dy\,f^2(y)\int_{\varepsilon\rho}^{K\rho}du\,p_G^{(Q)}(y-x,u)\\
&\qquad\qquad\qquad-\|\nabla_{\mathbb{R}^3}f\|_2^2\biggr]\\
&\quad=\mathcal{P}_p^{(Q)}(\varepsilon,K;\gamma,\rho,\nu).
\end{aligned}
\tag{7.52}
$$

Finally, let $Q\uparrow\mathbb{R}^3$ and use the following:

Lemma 7.4. *Let $\mathcal{P}_p(\varepsilon,K;\gamma,\rho,\nu)$ be as defined in (5.22). Then*

$$
\limsup_{Q\uparrow\mathbb{R}^3}\mathcal{P}_p^{(Q)}(\varepsilon,K;\gamma,\rho,\nu)\le\mathcal{P}_p(\varepsilon,K;\gamma,\rho,\nu).
\tag{7.53}
$$

The proof of Lemma 7.4 is deferred to Section 7.4. Combining (7.52) and (7.53), we have completed the proof of Lemma 5.5. □

We close this section by proving Lemma 7.3.

Proof of Lemma 7.3. Abbreviate

$$\nu_s = \frac{1}{p}\sum_{k=1}^{p}\mu_{k,s} \in \mathcal{M}_1(Q), \qquad 0 \le s \le (N+M+1)\delta\kappa. \tag{7.54}$$

Since $\mu \mapsto \widehat{S}_Q(\mu)$ is convex, we have

$$\begin{aligned}\text{lhs } (7.50) \le{}& \frac{p^2}{2}\frac{A}{\kappa}\int_0^{(N+M+1)\delta\kappa} ds \int_0^{(N+M+1)\delta\kappa} dt\, \mathbb{1}_{\{(L-1)\delta\kappa \le |t-s| \le (M+1)\delta\kappa\}}\\ &\times \int_Q \nu_s(dx)\int_Q \nu_t(dy) p_G^{(Q)}\Big(y-x, \frac{\rho}{\kappa}|t-s| + 2\delta\rho\Big)\\ &- p\int_0^{(N+M+1)\delta\kappa} ds \widehat{S}_Q(\nu_s),\end{aligned} \tag{7.55}$$

where we symmetrize the integrals w.r.t. $s$ and $t$. Let $B > 0$ be the size of $Q$, that is, $Q = [-B, B)^3$. Then $p_G^{(Q)}$ admits the Fourier representation

$$p_G^{(Q)}(x,t) = \frac{1}{(2B)^3}\sum_{q\in\mathbb{Z}^3} e^{-(\pi/B)^2|q|^2 t} e^{-i(\pi/B)q\cdot x}, \qquad x \in Q,\ t > 0. \tag{7.56}$$

Let

$$\widehat{\nu}_s(q) = \int_Q e^{i(\pi/B)q\cdot x}\nu_s(dx), \qquad q \in \mathbb{Z}^3. \tag{7.57}$$

Then we may rewrite

$$\begin{aligned}\text{rhs } (7.55) ={}& \frac{p^2}{2}\frac{A}{\kappa}\int_0^{(N+M+1)\delta\kappa} ds \int_0^{(N+M+1)\delta\kappa} dt\, \mathbb{1}_{\{(L-1)\delta\kappa \le |t-s| \le (M+1)\delta\kappa\}}\\ &\times \frac{1}{(2B)^3}\sum_{q\in\mathbb{Z}^3} e^{-(\pi/B)^2|q|^2[(\rho/\kappa)|t-s|+2\delta\rho]}\widehat{\nu}_s(q)\overline{\widehat{\nu}_t(q)}\\ &- p\int_0^{(N+M+1)\delta\kappa} ds\, \widehat{S}_Q(\nu_s).\end{aligned} \tag{7.58}$$

Since this expression is real-valued and

$$\mathrm{Re}(\widehat{\nu}_s(q)\overline{\widehat{\nu}_t(q)}) \le \tfrac12|\widehat{\nu}_s(q)|^2 + \tfrac12|\widehat{\nu}_t(q)|^2, \tag{7.59}$$

we find, after inserting (7.59) into (7.58) and afterward undoing the symmetrization w.r.t. $s$ and $t$, that

$$\begin{aligned}\text{rhs (7.58)} \le{}& p^2 \frac{A}{\kappa}\int_0^{(N+M+1)\delta\kappa} ds \int_{s+(L-1)\delta\kappa}^{s+(M+1)\delta\kappa} dt \\ &\times \frac{1}{(2B)^3}\sum_{q\in\mathbb{Z}^3} e^{-(\pi/B)^2|q|^2[(\rho/\kappa)(t-s)+2\delta\rho]}|\widehat{\nu}_s(q)|^2 \\ &- p\int_0^{(N+M+1)\delta\kappa} ds\, \widehat{S}_Q(\nu_s).\end{aligned} \tag{7.60}$$

Again using (7.56) and (7.57), we see that

$$\begin{aligned}&\text{rhs (7.60)}\\ &\quad = p\int_0^{(N+M+1)\delta\kappa} ds \\ &\qquad \times \left[\frac{A}{\kappa}p\int_Q \nu_s(dx)\int_Q \nu_s(dy)\int_{(L-1)\delta\kappa}^{(M+1)\delta\kappa} du\, p_G^{(Q)}\left(y-x,\frac{\rho}{\kappa}u+2\delta\rho\right)\right. \\ &\qquad\qquad \left. - \widehat{S}_Q(\nu_s)\right].\end{aligned} \tag{7.61}$$

Clearly, this expression does not exceed the right-hand side of (7.50). □

7.4. *Proofs of Lemmas 7.1, 7.2 and 7.4.*

Proof of Lemma 7.1. Let $X^{(\kappa)}$ be the scaled random walk on $\mathbb{Z}^3_\kappa$ [as in (7.1)], let $\tau^{(\kappa)}$ be the first time $X^{(\kappa)}$ exits $Q$, and let $\Xi^{(\kappa)}$ be the normalized occupation time measure of $X^{(\kappa)}$ [as in (7.13)]. Define the conditional probability measures

$$\mathbb{Q}_x^{(\kappa)}(\cdot) = \mathbb{P}_x^{X^{(\kappa)}}(\Xi^{(\kappa)} \in \cdot \,|\tau^{(\kappa)} > \delta\kappa,\ X^{(\kappa)}(\delta\kappa) \in \tfrac{1}{2}Q). \tag{7.62}$$

Let $\zeta_0$ denote the principal eigenvalue of the Laplacian $\Delta_Q$ with Dirichlet boundary condition in $L^2(Q)$. We will prove the following:

(a) uniformly in $x \in \frac{1}{2}Q$,

$$\lim_{\kappa\to\infty}\frac{1}{\delta\kappa}\log \mathbb{P}_x^{X^{(\kappa)}}\left(\tau^{(\kappa)} > \delta\kappa,\ X^{(\kappa)}(\delta\kappa) \in \frac{1}{2}Q\right) = \zeta_0; \tag{7.63}$$

(b) the family $(\mathbb{Q}_x^{(\kappa)})_{\kappa>0}$ satisfies the full large deviation principle on $\mathcal{M}_1(Q)$, uniformly in $x \in \frac{1}{2}Q$, with rate $\delta\kappa$ and with rate function $S_Q + \zeta_0$ [recall (7.22)].

As a consequence of (a) and (b), the family $(\widetilde{\mathbb{Q}}_x^{(\kappa)})_{\kappa>0}$ of sub-probability measures defined by

$$\widetilde{\mathbb{Q}}_x^{(\kappa)}(\cdot) = \mathbb{P}_x^{X^{(\kappa)}}(\Xi^{(\kappa)} \in \cdot\,,\ \tau^{(\kappa)} > \delta\kappa,\ X^{(\kappa)}(\delta\kappa) \in \tfrac{1}{2}Q) \tag{7.64}$$

satisfies the full large deviation principle on $\mathcal{M}_1(Q)$, uniformly in $x \in \frac{1}{2}Q$, with rate $\delta\kappa$ and rate function $S_Q$. The latter, in turn, implies Lemma 7.1.

The proof of assertions (a) and (b) is achieved as follows. Given a potential $V \in C_c^\infty(Q)$, let $\zeta_0(V)$ denote the principal eigenvalue of $\Delta_Q + V$ with Dirichlet boundary condition in $L^2(Q)$. It is well known that $V \mapsto \zeta_0(V)$ is Gateaux differentiable and that $S_Q$ has the following representation as a Legendre transform:

$$S_Q(\mu) = \sup_{V \in C_c^\infty(Q)} \left[\int_Q V\, d\mu - \zeta_0(V)\right], \qquad \mu \in C_c^\infty(Q)^*, \tag{7.65}$$

where $C_c^\infty(Q)^*$ is the algebraic dual of $C_c^\infty(Q)$ equipped with the weak* topology [(7.65) is dual to the Rayleigh–Ritz formula for $\zeta_0(V)$]. We may therefore apply a uniform (w.r.t. the starting point) version of Dawson and Gärtner [9], Theorem 3.4, to see that, in order to prove (a) and (b), it is enough to show that

$$\begin{aligned} \lim_{\kappa\to\infty} \frac{1}{\delta\kappa} \log \mathbb{E}_x^{X^{(\kappa)}} \biggl(&\exp\biggl[\int_0^{\delta\kappa} V(X^{(\kappa)}(s))\, ds\biggr] \\ &\times \mathbb{1}\biggl\{\tau^{(\kappa)} > \delta\kappa, X^{(\kappa)}(\delta\kappa) \in \frac{1}{2}Q\biggr\}\biggr) \\ = \zeta_0(V), \end{aligned} \tag{7.66}$$

uniformly in $x \in \frac{1}{2}Q$ for all $V \in C_c^\infty(Q)$. (An argument similar to that in [9], Section 3.5, shows that $S_Q(\mu) < \infty$, $\mu \in C_c^\infty(Q)^*$ imply $\mu \in \mathcal{M}_1(Q)$, which is needed for the application of [9], Theorem 3.4.) Note that assertion (a) coincides with (7.66) for $V = 0$.

Fix $V \in C_c^\infty(Q)$. Abbreviate

$$\begin{aligned} s_-^{(\kappa)}(t) = \log \inf_{x \in (1/2)Q} \mathbb{E}_x^{X^{(\kappa)}} \biggl(&\exp\biggl[\int_0^t V(X^{(\kappa)}(s))\, ds\biggr] \\ &\times \mathbb{1}\Bigl\{\tau^{(\kappa)} > t,\ X^{(\kappa)}(t) \in \tfrac{1}{2}Q\Bigr\}\biggr). \end{aligned} \tag{7.67}$$

Fix $T > 0$. For $t = \delta\kappa$, split the integral in the right-hand side of (7.67) into the sum of $\lfloor \delta\kappa/T \rfloor$ integrals over intervals of length $T_\kappa = \delta\kappa/\lfloor \delta\kappa/T \rfloor$. Then, using the Markov property of $X^{(\kappa)}$ at the splitting points, we get

$$s_-^{(\kappa)}(\delta\kappa) \ge \lfloor \delta\kappa/T \rfloor s_-^{(\kappa)}(T_\kappa). \tag{7.68}$$

Hence,

$$\liminf_{\kappa\to\infty}\frac{s_-^{(\kappa)}(\delta\kappa)}{\delta\kappa}\geq\frac{1}{T}\liminf_{\kappa\to\infty}s_-^{(\kappa)}(T_\kappa)$$

$$=\frac{1}{T}\log\inf_{x\in(1/2)Q}\mathbb{E}_x^W\biggl(\exp\biggl[\int_0^T V(W(s))\,ds\biggr]\times\mathbb{1}\biggl\{\tau>T,\ W(T)\in\frac{1}{2}Q\biggr\}\biggr),\tag{7.69}$$

where $W$ is Brownian motion on $\mathbb{R}^3$ with generator $\Delta_{\mathbb{R}^3}$ and $\tau$ denotes the first time $W$ exits $Q$. To derive the last line of (7.69) we use a uniform version of Donsker's invariance principle. It is well known that the right-hand side of (7.69) tends to $\zeta_0(V)$ as $T\to\infty$. Therefore, we arrive at the lower bound

$$\liminf_{\kappa\to\infty}\frac{s_-^{(\kappa)}(\delta\kappa)}{\delta\kappa}\geq\zeta_0(V).\tag{7.70}$$

To get the corresponding upper bound, abbreviate

$$s_+^{(\kappa)}(t)=\log\sup_{x\in Q}\mathbb{E}_x^{X^{(\kappa)}}\biggl(\exp\biggl[\int_0^t V(X^{(\kappa)}(s))\,ds\biggr]\mathbb{1}\{\tau^{(\kappa)}>t\}\biggr).\tag{7.71}$$

Then, in analogy with the above considerations, we obtain, through a superadditivity argument, that

$$\limsup_{\kappa\to\infty}\frac{s_+^{(\kappa)}(\delta\kappa)}{\delta\kappa}\leq\zeta_0(V).\tag{7.72}$$

We then combine (7.70) and (7.72) to get (7.66). □

Proof of Lemma 7.2. Let $X^{(\kappa,Q)}$ denote the random walk on $Q^{(\kappa)}=Q\cap\mathbb{Z}_\kappa^3$ obtained by wrapping $X^{(\kappa)}$ around $Q^{(\kappa)}$ [recall (7.26)]. Let

$$\widehat{\Xi}^{(\kappa)}(A)=\frac{1}{\delta\kappa}\int_0^{\delta\kappa}ds\,\mathbb{1}_A(X^{(\kappa,Q)}(s)),\qquad A\subset\mathbb{R}^3\text{ Borel},\tag{7.73}$$

and

$$\widehat{\mathbb{Q}}_x^{(\kappa)}(\cdot)=\mathbb{P}_x^{X^{(\kappa,Q)}}(\Xi^{(\kappa)}\in\cdot).\tag{7.74}$$

Then the analogue of (b) reads as follows:

(b′) The family $(\widehat{\mathbb{Q}}_x^{(\kappa)})_{\kappa>0}$ satisfies the full large deviation principle on $\mathcal{M}_1(Q)$, uniformly in $x\in Q$, with rate $\delta\kappa$ and rate function $\widehat{S}_Q$ [recall (7.45)].

The proof of assertion (b′) follows the same lines as the proofs of assertions (a) and (b) and is, in fact, even simpler. Using (b′) together with (7.73) and (7.74), we arrive at the assertion claimed in Lemma 7.2. □

Proof of Lemma 7.4. Let $Q = Q_B = [-B,B)^3$. Write $Q_B(q) = Q_B + q$, $q \in \mathbb{R}^3$. Let

$$\widehat{p}^{(Q_B)}(x,t) = \sum_{k\in\mathbb{Z}^3} p_G(x + 2Bk, t) \tag{7.75}$$

denote the $Q_B$-periodization of the Gaussian transition kernel $p_G$. Recall that $H^1_{\mathrm{per}}(Q_B)$ denotes the space of functions in $H^1(Q_B)$ with periodic boundary conditions.

Fix $B > 1$ and $f \in H^1_{\mathrm{per}}(Q_B)$ with $\|f\|_2 = 1$. Put $A = B - \sqrt{B}$. Let $\widehat{f}$ denote the $Q_B$-periodic extension of $f$ to $\mathbb{R}^3$. Then

$$\frac{1}{|Q_B|}\int_{Q_B} dq \int_{Q_A(q)} dx\, \widehat{f}^2(x) = \frac{|Q_A|}{|Q_B|} \tag{7.76}$$

and hence there exists $q \in Q_B$ (depending on $B, f$) such that

$$\int_{Q_A(q)} dx\, \widehat{f}^2(x) \ge \frac{|Q_A|}{|Q_B|}. \tag{7.77}$$

Let $h_B : \mathbb{R}^3 \to \mathbb{R}$ be a smooth function (depending on $B, q$) satisfying

$$0 \le h_B \le 1, \qquad h_B = \begin{cases} 1, & \text{on } Q_A(q), \\ 0, & \text{on } \mathbb{R}^3 \setminus Q_B(q). \end{cases} \tag{7.78}$$

We may assume that

$$D = \|\Delta(h_B(1-h_B)) + 2|\nabla h_B|^2\|_\infty < \infty \tag{7.79}$$

with $D$ not dependent on $B, q, f$. Define

$$f_B = \frac{h_B \widehat{f}}{\|h_B \widehat{f}\|_2}. \tag{7.80}$$

Then $f_B \in H^1(\mathbb{R}^3)$ and $\|f_B\|_2 = 1$. Moreover, by (7.77) and (7.78) we have

$$\frac{|Q_A|}{|Q_B|} \le \|h_B \widehat{f}\|_2^2 \le 1. \tag{7.81}$$

Hence, $\|h_B \widehat{f}\|_2 \to 1$ as $B \to \infty$.

Next, observe that

$$\|x - y + 2Bk\|_\infty \ge 2B(\|k\|_\infty - 1) + 2(B - A), \qquad x, y \in Q_A(q),\ k \in \mathbb{Z}^3 \setminus \{0\}. \tag{7.82}$$

Because

$$p_G(x,t) = (4\pi t)^{-3/2} \exp[-\|x\|^2/4t] \le (4\pi t)^{-3/2} \exp[-\|x\|_\infty^2/4t], \tag{7.83}$$

it follows from (7.82) that there exists $\delta_B$ (not depending on $q, f$), satisfying $\delta_B \to 0$ as $B \to \infty$, such that

$$\int_{\varepsilon\rho}^{K\rho} \widehat{p}_G^{(Q_B)}(x-y,t)\,dt \le \int_{\varepsilon\rho}^{K\rho} p_G(x-y,t)\,dt + \delta_B, \qquad x,y \in Q_A(q). \tag{7.84}$$

Moreover, from this it also follows that there exists a constant $C < \infty$ (not depending on $B \ge 1, q, f$) such that

$$\int_{\varepsilon\rho}^{K\rho} \widehat{p}_G^{(Q_B)}(x-y,t)\,dt \le C, \qquad x,y \in Q_A(q). \tag{7.85}$$

With the above estimates in place, we next derive an upper bound for

$$\int_{Q_B} dx \int_{Q_B} dy \int_{\varepsilon\rho}^{K\rho} dt\, \widehat{p}_G^{(Q_B)}(x-y,t) \widehat{f}^{\,2}(x) \widehat{f}^{\,2}(y). \tag{7.86}$$

Since $\widehat{f}$ is $Q_B$-periodic, we may replace the domain of integration $Q_B \times Q_B$ by $Q_B(q) \times Q_B(q)$. After that, we may split the integral into two parts: $Q_A(q) \times Q_A(q)$ and $[Q_B(q) \times Q_B(q)] \setminus [Q_A(q) \times Q_A(q)]$. The latter coincides with the union of $[Q_B(q) \setminus Q_A(q)] \times Q_B(q)$ and $Q_B(q) \times [Q_B(q) \setminus Q_A(q)]$. Therefore, using (7.77), (7.84) and (7.85), we obtain

$$\begin{aligned}
&\int_{Q_B} dx \int_{Q_B} dy \int_{\varepsilon\rho}^{K\rho} dt\, \widehat{p}_G^{(Q_B)}(x-y,t) \widehat{f}^{\,2}(x) \widehat{f}^{\,2}(y) \\
&\quad \le \int_{Q_A(q)} dx \int_{Q_A(q)} dy \int_{\varepsilon\rho}^{K\rho} dt\, p_G(x-y,t) \widehat{f}^{\,2}(x) \widehat{f}^{\,2}(y) \\
&\qquad + \delta_B + 2C \int_{Q_B(q)\setminus Q_A(q)} dx\, \widehat{f}^{\,2}(x) \\
&\quad \le \int_{Q_A(q)} dx \int_{Q_A(q)} dy \int_{\varepsilon\rho}^{K\rho} dt\, p_G(x-y,t) f_B^2(x) f_B^2(y) \\
&\qquad + \delta_B + 2C \frac{|Q_B \setminus Q_A|}{|Q_B|} \\
&\quad \le \frac{|Q_A|}{|Q_B|} \int_{Q_A(q)} dx \int_{Q_A(q)} dy \int_{\varepsilon\rho}^{K\rho} dt\, p_G(x-y,t) f^2(x) f^2(y) \\
&\qquad + \delta_B + 3C \frac{|Q_B \setminus Q_A|}{|Q_B|},
\end{aligned} \tag{7.87}$$

where, in the second inequality, we use the fact that $\widehat{f}^{\,2} = (h_B \widehat{f})^2 = \|h_B \widehat{f}\|_2^2 f_B^2 \le f_B^2$ on $Q_A(q)$.

Next, we derive a lower bound for $\|\nabla f\|_2^2$ in terms of $f_B$. First, estimate

$$\begin{aligned}\|\nabla f\|_2^2 &= \int_{Q_B(q)} dx\,|\nabla(h_B\widehat{f}\,) + \nabla((1-h_B)\widehat{f}\,)|^2 \\ &\geq \int_{Q_B(q)} dx|\nabla(h_B\widehat{f})|^2 + 2\int_{Q_B(q)} dx\,\nabla(h_B\widehat{f})\cdot\nabla((1-h_B)\widehat{f}\,).\end{aligned}\tag{7.88}$$

But

$$\nabla(h_B\widehat{f})\cdot\nabla((1-h_B)\widehat{f}\,) \geq (\widehat{f}\,\nabla(h_B(1-h_B)))\cdot\nabla\widehat{f} - |\nabla h_B|^2\widehat{f}^{\,2}\tag{7.89}$$

and integration by parts shows that

$$\begin{aligned}&\int_{Q_B(q)} dx\,(\widehat{f}\nabla(h_B(1-h_B)))\cdot\nabla\widehat{f} \\ &\qquad = -\tfrac{1}{2}\int_{Q_B(q)} dx\,\widehat{f}^{\,2}\Delta(h_B(1-h_B)).\end{aligned}\tag{7.90}$$

Hence, recalling the definition of $f_B$ and taking into account (7.77), (7.79) and (7.81), we obtain

$$\begin{aligned}\|\nabla f\|_2^2 &\geq \|h_B\widehat{f}\|_2^2\|\nabla f_B\|_2^2 \\ &\qquad - \int_{Q_B(q)\setminus Q_A(q)} dx\,\widehat{f}^{\,2}[\nabla(h_B(1-h_B)) + 2|\nabla h_B|^2] \\ &\geq \frac{|Q_A|}{|Q_B|}\|\nabla f_B\|_2^2 - \frac{|Q_B\setminus Q_A|}{|Q_B|}D.\end{aligned}\tag{7.91}$$

Combining (7.87) and (7.91) and abbreviating $\alpha = (\nu\gamma^2/\rho)p$, we arrive at

$$\begin{aligned}&\alpha\int_{Q_B} dx\int_{Q_B} dy\int_{\varepsilon\rho}^{K\rho} dt\,\widehat{p}_G(x-y,t)f^2(x)f^2(y) - \|\nabla f\|_2^2 \\ &\qquad \leq \frac{|Q_A|}{|Q_B|}\mathcal{P} + \alpha\delta_B + (3\alpha C + D)\frac{|Q_B\setminus Q_A|}{|Q_B|}.\end{aligned}\tag{7.92}$$

Since $C, D$ and $\delta_B$ do not depend on $f$, we conclude that [recalling (7.52)]

$$\begin{aligned}&\mathcal{P}_p^{(Q_B)}(\varepsilon, K;\gamma,\rho,\nu) \\ &\qquad \leq \frac{|Q_A|}{|Q_B|}\mathcal{P}_p(\varepsilon, K;\gamma,\rho,\nu) + \alpha\delta_B + (3\alpha C + D)\frac{|Q_B\setminus Q_A|}{|Q_B|}.\end{aligned}\tag{7.93}$$

Now let $B\to\infty$ and use the fact that $\delta_B\to 0$ and $|Q_A|/|Q_B|\to 1$, to arrive at the assertion claimed in (7.53). □

**8. Proofs of Lemmas 5.6–5.8.** In this section, we prove Lemmas 5.6–5.8, which deal with the terms that are asymptotically negligible as $\kappa\to\infty$.

8.1. *Proof of Lemma* 5.6.

PROOF OF LEMMA 5.6. Using the rough bound

$$p_{\rho/\kappa}(X(t)-X(s),t-s)\le p_{\rho/\kappa}(0,t-s)=p\left(0,\frac{\rho}{\kappa}(t-s)\right), \tag{8.1}$$

we conclude from (5.30) and (5.31) that

$$\kappa^2\lambda_{\mathrm{off}}^{+}(a,\kappa)\le\frac{\nu\gamma^2}{\rho}\kappa\int_{\rho a\kappa^2}^{\infty}dt\,p(0,t). \tag{8.2}$$

Because of (5.43) and (5.44), the expression in the right-hand side is bounded above by a constant times $a^{-(d-2)/2}\kappa^{-(d-3)}$. From this, the claims in (5.32) and (5.33) follow. □

8.2. *Proof of Lemma* 5.7. For the proof of Lemma 5.7, we need two more lemmas. Let $\mathcal{G}$ denote the Green operator acting on functions $V:\mathbb{Z}^d\to[0,\infty)$ as

$$(\mathcal{G}V)(x)=\sum_{y\in\mathbb{Z}^d}G(y-x)V(y),\qquad x\in\mathbb{Z}^d, \tag{8.3}$$

with $G(z)=\int_0^\infty dt\,p(z,t)$. Let $\|\cdot\|_\infty$ denote the supremum norm.

LEMMA 8.1. *For any* $V:\mathbb{Z}^d\to[0,\infty)$ *and* $x\in\mathbb{Z}^d$,

$$\mathbb{E}_x^X\left(\exp\left[\int_0^\infty dt\,V(X(t))\right]\right)\le(1-\|\mathcal{G}V\|_\infty)^{-1}, \tag{8.4}$$

*provided that*

$$\|\mathcal{G}V\|_\infty<1. \tag{8.5}$$

LEMMA 8.2. *For any* $\alpha,\beta>0$ *and* $a>0$,

$$\begin{aligned}&\mathbb{E}_{0,0}^{X,Y}\left(\exp\left[\alpha\int_0^T ds\int_s^{s+a\kappa^3}dt\,p_\beta(Y(t)-X(s),t-s)\right]\right)\\&\quad\le\mathbb{E}_0^X\left(\exp\left[\alpha\int_0^T ds\int_s^{s+a\kappa^3}dt\,p_\beta(X(s),t-s)\right]\right).\end{aligned} \tag{8.6}$$

Before giving the proofs of Lemmas 8.1 and 8.2, we first prove Lemma 5.7.

PROOF OF LEMMA 5.7. Using Lemma 8.2, we get from (5.34) that

$$\begin{aligned}\Lambda_{\mathrm{mix}}(T;a,\kappa)&\le\frac{1}{T}\log\mathbb{E}_0^X\left(\exp\left[\frac{\nu\gamma^2}{\kappa^2}\int_0^T ds\int_s^{s+a\kappa^3}dt\,p_{\rho/\kappa}(X(s),t-s)\right]\right)\\&\le\frac{1}{T}\log\mathbb{E}_0^X\left(\exp\left[\int_0^\infty ds\,V_{a,\kappa}(X(s))\right]\right),\end{aligned} \tag{8.7}$$

where

$$V_{a,\kappa}(x) = \frac{\nu\gamma^2}{\rho\kappa}\int_0^{\rho a\kappa^2} dt\, p(x,t), \qquad x \in \mathbb{Z}^d. \tag{8.8}$$

It follows from (5.43) and (5.44) that, as $\kappa \to \infty$,

$$\|\mathcal{G}V_{a,\kappa}\|_\infty = \frac{\nu\gamma^2}{\rho\kappa}\int_0^{\rho a\kappa^2} dt \int_t^\infty ds\, p(0,s) \tag{8.9}$$

tends to zero for $d \geq 4$ and $0 < a < \infty$ and tends to a constant times $a^{1/2}$ for $d = 3$. Hence, by Lemma 8.1, for large $\kappa$ the expectation in the right-hand side of (8.7) is finite for $0 < a < a_0$ with $a_0 = \infty$ for $d \geq 4$ and $a_0 > 0$ sufficiently small for $d = 3$. Thus, by letting $T \to \infty$ in (8.7), we conclude that

$$\lambda^+_{\text{mix}}(a,\kappa) = 0 \qquad \forall\, 0 < a < a_0,\ \kappa \geq \kappa_0(a). \tag{8.10}$$

This yields (5.37). To prove (5.36), simply note that for all $0 < a < \infty$,

$$\Lambda_{\text{mix}}(T;\infty,\kappa) \leq \Lambda_{\text{mix}}(T;a,\kappa) + \frac{\nu\gamma^2}{\rho\kappa}\int_{\rho a\kappa^2}^\infty dt\, p(0,t) \tag{8.11}$$

and hence

$$\kappa^2\lambda^+_{\text{mix}}(\infty,\kappa) \leq \frac{\nu\gamma^2}{\rho}\kappa\int_{\rho a\kappa^2}^\infty dt\, p(0,t) \qquad \forall\, 0 < a < a_0,\ \kappa \geq \kappa_0(a). \tag{8.12}$$

Now, proceed as with (8.2) to obtain the claimed assertion. □

### 8.3. *Proofs of Lemmas 8.1 and 8.2.*

Proof of Lemma 8.1. A Taylor expansion of the exponential function yields

$$\begin{aligned} &\mathbb{E}^X_x\left(\exp\left[\int_0^\infty dt\, V(X(t))\right]\right) \\ &\quad = \sum_{n=0}^\infty \int_0^\infty dt_1 \int_{t_1}^\infty dt_2 \cdots \int_{t_{n-1}}^\infty dt_n \\ &\qquad\qquad \times \mathbb{E}^X_x(V(X(t_1))V(X(t_2)) \times \cdots \times V(X(t_n))). \end{aligned} \tag{8.13}$$

But,

$$\int_0^\infty dt_1 \int_{t_1}^\infty dt_2 \cdots \int_{t_{n-1}}^\infty dt_n\, \mathbb{E}^X_x(V(X(t_1))V(X(t_2)) \times \cdots \times V(X(t_n)))$$

$$
\begin{aligned}
&= \sum_{y_1\in\mathbb{Z}^d} \int_0^\infty dt_1\, p(y_1-x,t_1)V(y_1) \\
&\quad \times \sum_{y_2\in\mathbb{Z}^d} \int_{t_1}^\infty dt_2\, p(y_2-y_1,t_2-t_1)V(y_2) \\
&\quad \times\cdots\times \sum_{y_n\in\mathbb{Z}^d} \int_{t_{n-1}}^\infty dt_n\, p(y_n-y_{n-1},t_n-t_{n-1})V(y_n) \qquad (8.14)\\
&= \sum_{y_1\in\mathbb{Z}^d} G(y_1-x)V(y_1) \sum_{y_2\in\mathbb{Z}^d} G(y_2-y_1)V(y_2) \\
&\quad \times\cdots\times \sum_{y_n\in\mathbb{Z}^d} G(y_n-y_{n-1})V(y_n) \\
&\le \|\mathcal{G}V\|_\infty^n.
\end{aligned}
$$

Substituting this into (8.13) and summing the geometric series, we arrive at the assertion claimed in (8.4). □

Proof of Lemma 8.2. Using the Fourier representation of the transition kernel [recalling (2.14)]

$$
p_\beta(x,t) = \oint dk\, e^{-\beta t\widehat{\varphi}(k)} e^{-ik\cdot x} \qquad (8.15)
$$

and expanding the exponential function in a Taylor series, we find that

$$
\begin{aligned}
&\mathbb{E}_{0,0}^{X,Y}\left(\exp\left[\alpha\int_0^T ds\int_s^{s+a\kappa^3} dt\, p_\beta(Y(t)-X(s),t-s)\right]\right) \\
&\quad = \sum_{n=0}^\infty \alpha^n \int_0^T ds_1 \int_{s_1}^T ds_2 \cdots \int_{s_{n-1}}^T ds_n \\
&\qquad \times \int_{s_1}^{s_1+a\kappa^3} dt_1 \int_{s_2}^{s_2+a\kappa^3} dt_2 \cdots \int_{s_n}^{s_n+a\kappa^3} dt_n \qquad (8.16)\\
&\qquad \times \oint dk_1 \oint dk_2 \cdots \oint dk_n \exp\left[-\beta\sum_{j=1}^n (t_j-s_j)\widehat{\varphi}(k_j)\right] \\
&\qquad \times \mathbb{E}_0^Y\left(\exp\left[-i\sum_{j=1}^n k_j\cdot Y(t_j)\right]\right)\mathbb{E}_0^X\left(\exp\left[i\sum_{j=1}^n k_j\cdot X(s_j)\right]\right).
\end{aligned}
$$

Here, to factorize the two expectations, we have used the fact that the random walks $X$ and $Y$ are independent. By symmetry of $X$ and $Y$, these two expectations are real-valued. An explicit computation shows that the second expectation is strictly positive. (Use the fact that the $s_i$ are ordered and that $X$ has independent increments so that the expectation factors into a product.) The first expectation is clearly less than or equal to 1. Hence, the above expression can be bounded from above by the same expression with $Y$ replaced by 0. This, in turn, yields (8.6). □

8.4. *Proof of Lemma* 5.8. We begin by noting two facts. First, define

$$\Lambda_{\mathrm{full}}(T;\kappa) = \frac{1}{T}\log \mathbb{E}_0^X\left(\exp\left[\frac{\nu\gamma^2}{\kappa^2}\int_0^T ds\int_s^\infty dt\, p_{\rho/\kappa}(X(t)-X(s),t-s)\right]\right) \tag{8.17}$$

and

$$\lambda_{\mathrm{full}}^+(\kappa) = \limsup_{T\to\infty}\Lambda_{\mathrm{full}}(T;\kappa). \tag{8.18}$$

By splitting the second integral in the right-hand side of (8.17) into a diagonal, a variational and an off-diagonal part (in accordance with Lemmas 5.4–5.6), applying Hölder's inequality to separate the parts [similarly as in (5.66)] and applying Lemmas 5.4–5.6, we find that

$$\limsup_{\kappa\to\infty}\kappa^2\lambda_{\mathrm{full}}^+(\kappa) \le \frac{\nu\gamma^2}{r_d}, \qquad \text{if } d\ge 4, \tag{8.19}$$

while

$$\limsup_{\kappa\to\infty}\kappa^2\lambda_{\mathrm{full}}^+(\kappa) \le \frac{\nu\gamma^2}{r_3} + \left(\frac{\nu\gamma^2}{\rho}\right)^2\mathcal{P}, \qquad \text{if } d=3. \tag{8.20}$$

Second, note that Lemma 6.3 for $k=0$ yields the bound

$$\mathbb{E}_0^X\left(\exp\left[\alpha\int_0^\infty dt\, p_{\rho/\kappa}(X(t),t)\right]\right) \le \exp\left[\frac{\alpha G_0(0)}{1-\alpha G_0(0)}\right] \le \exp\left[\frac{2\alpha}{r_d}\right], \tag{8.21}$$

provided that

$$0\le\alpha\le\frac{r_d}{2}. \tag{8.22}$$

Proof of Lemma 5.8. Using the rough bound (8.1), we have

$$\begin{aligned}
&\int_0^T ds\left(\int_s^\infty dt\, p_{\rho/\kappa}(X(t)-X(s),t-s)\right)\\
&\quad\times\left(\int_0^s du\, p_{\rho/\kappa}(X(s)-X(u),s-u)\right)\\
&\le \int_0^T ds\left(\int_s^{s+\kappa^{3/2}} dt\, p_{\rho/\kappa}(X(t)-X(s),t-s)\right)
\end{aligned}\tag{8.23}$$

$$\times \left( \int_{s-\kappa^{3/2}}^{s} du \, p_{\rho/\kappa}(X(s)-X(u), s-u) \right)$$
$$+ 2\left( \int_{\kappa^{3/2}}^{\infty} du \, p_{\rho/\kappa}(0,u) \right) \int_0^T ds \int_s^\infty dt \, p_{\rho/\kappa}(X(t)-X(s), t-s).$$

Substituting this into (5.38) and applying the Cauchy–Schwarz inequality, we find that

$$\Lambda_{\text{rem}}(T;\kappa) \le \Lambda_{\text{rem}}^{(1)}(T;\kappa) + \Lambda_{\text{rem}}^{(2)}(T;\kappa), \tag{8.24}$$

where

$$\begin{aligned} \Lambda_{\text{rem}}^{(1)}(T;\kappa) = \frac{1}{2T} \log \mathbb{E}_0^X \Bigg( \exp\Bigg[ \frac{2\nu\gamma^3}{\kappa^3} \int_0^T ds \\ \times \left( \int_s^{s+\kappa^{3/2}} dt \, p_{\rho/\kappa}(X(t)-X(s), t-s) \right) \\ \times \left( \int_{s-\kappa^{3/2}}^{s} du \, p_{\rho/\kappa}(X(s)-X(u), s-u) \right) \Bigg] \Bigg) \end{aligned} \tag{8.25}$$

and

$$\begin{aligned} \Lambda_{\text{rem}}^{(2)}(T;\kappa) = \frac{1}{2T} \log \mathbb{E}_0^X \Bigg( \exp\Bigg[ \frac{\nu\gamma^2}{\kappa^2} \left( \frac{4\gamma}{\kappa} \int_{\kappa^{3/2}}^{\infty} du \, p_{\rho/\kappa}(0,u) \right) \\ \times \int_0^T ds \int_s^\infty dt \, p_{\rho/\kappa}(X(t)-X(s), t-s) \Bigg] \Bigg). \end{aligned} \tag{8.26}$$

To prove Lemma 5.8, it will be enough to show that

$$\lim_{\kappa\to\infty} \kappa^2 \limsup_{T\to\infty} \Lambda_{\text{rem}}^{(i)}(T;\kappa) = 0, \qquad i = 1, 2. \tag{8.27}$$

Since for $d \ge 3$,

$$\frac{4\gamma}{\kappa} \int_{\kappa^{3/2}}^{\infty} du \, p_{\rho/\kappa}(0,u) \to 0 \qquad \text{as } \kappa \to \infty, \tag{8.28}$$

(8.27) for $i = 2$ follows from (8.17)–(8.20) with $\nu$ replaced by $\nu$ times the integral in (8.28). To prove (8.27) for $i = 1$, we split the integral in the right-hand side of (8.25) as follows:

$$\int_0^T ds = \left( \sum_{\substack{k=1 \\ \text{even}}}^{\lceil T/2\kappa^{3/2} \rceil} + \sum_{\substack{k=1 \\ \text{odd}}}^{\lceil T/2\kappa^{3/2} \rceil} \right) \int_{(k-1)2\kappa^{3/2}}^{k2\kappa^{3/2}} ds. \tag{8.29}$$

Note that the summands in each of the two sums are i.i.d. Hence, substituting (8.29) into (8.25) and applying the Cauchy–Schwarz inequality, we find

that

$$
\begin{aligned}
&\Lambda^{(1)}_{\mathrm{rem}}(T;\kappa)\\
&\quad\le \frac{\lceil T/2\kappa^{3/2}\rceil}{4T}\log \mathbb{E}^X_0\biggl(\exp\biggl[\frac{4\nu\gamma^3}{\kappa^3}\int_0^{2\kappa^{3/2}} ds \\
&\qquad\qquad \times\biggl(\int_s^{s+\kappa^{3/2}} dt\, p_{\rho/\kappa}(X(t)-X(s),t-s)\biggr)\\
&\qquad\qquad \times\biggl(\int_{s-\kappa^{3/2}}^{s} du\, p_{\rho/\kappa}(X(s)-X(u),s-u)\biggr)\biggr]\biggr).
\end{aligned}
\tag{8.30}
$$

Letting $T\to\infty$ and applying Jensen's inequality, we arrive at

$$
\begin{aligned}
&\limsup_{T\to\infty}\Lambda^{(1)}_{\mathrm{rem}}(T;\kappa)\\
&\quad\le \frac{1}{8\kappa^{3/2}}\log \mathbb{E}^{X,Y}_{0,0}\biggl(\exp\biggl[\frac{4\nu\gamma^3}{\kappa^3}2\kappa^{3/2}\biggl(\int_0^{\kappa^{3/2}} dt\, p_{\rho/\kappa}(X(t),t)\biggr)\\
&\qquad\qquad \times\biggl(\int_0^{\kappa^{3/2}} du\, p_{\rho/\kappa}(Y(u),u)\biggr)\biggr]\biggr),
\end{aligned}
\tag{8.31}
$$

where we use the fact that the increments of $X$ over the time intervals $[s,s+\kappa^{3/2}]$ and $[s-\kappa^{3/2},s]$ are independent in order to replace the expectation over the single random walk $X$ by an expectation over the two independent random walks $X,Y$. Since for $d\ge 3$,

$$
\begin{aligned}
\frac{4\nu\gamma^3}{\kappa^3}2\kappa^{3/2}\int_0^\infty du\, p_{\rho/\kappa}(Y(u),u) &\le \frac{8\nu\gamma^3}{\kappa^{3/2}}\int_0^\infty du\, p_{\rho/\kappa}(0,u)\\
&= \frac{8\nu\gamma^3}{r_d\rho\kappa^{1/2}}\to 0 \qquad \text{as } \kappa\to\infty,
\end{aligned}
\tag{8.32}
$$

we may apply (8.21) and (8.22) with $\alpha$ equal to the left-hand side of (8.32) to see that for large $\kappa$,

$$
\text{rhs (8.31)} \le \frac{1}{8\kappa^{3/2}}\log \mathbb{E}^Y_0\biggl(\exp\biggl[\frac{2}{r_d}\frac{8\nu\gamma^3}{\kappa^{3/2}}\int_0^\infty du\, p_{\rho/\kappa}(Y(u),u)\biggr]\biggr).
\tag{8.33}
$$

Finally, we may apply (8.21) and (8.22) once more, this time with $\alpha = 16\nu\gamma^3/r_d\kappa^{3/2}$, to obtain that for large $\kappa$,

$$
\limsup_{T\to\infty}\Lambda^{(1)}_{\mathrm{rem}}(T;\kappa)\le \frac{1}{8\kappa^{3/2}}\frac{2}{r_d}\biggl(\frac{16\nu\gamma^3}{r_d\kappa^{3/2}}\biggr) = \frac{4\nu\gamma^3}{r_d^2\kappa^3}.
\tag{8.34}
$$

This implies (8.27) for $i=1$. $\square$

**Acknowledgment.** The first author is grateful for hospitality at EURANDOM.

Institut für Mathematik
Technische Universität Berlin
Strasse des 17. Juni 136
D-10623 Berlin
Germany
E-mail: jg@math.tu-berlin.de

Mathematical Institute
Leiden University
P.O. Box 9512
2300 CA Leiden
The Netherlands
and
EURANDOM
P.O. Box 513
5600 MB Eindhoven
The Netherlands
E-mail: denholla@math.leidenuniv.nl